\documentclass[oneside,notitlepage,12pt]{article}

\usepackage{amssymb}
\usepackage[leqno]{amsmath}
\usepackage{amsfonts}
\usepackage{amsopn}
\usepackage{amstext}
\usepackage{amsthm}

\usepackage[all]{xy}
\newdir{ >}{{}*!/-9pt/@{>}}

\usepackage{verbatim}
\usepackage[colorlinks, backref]{hyperref}
\usepackage{makeidx}

\usepackage{calrsfs}
\usepackage{fourier-orns}
\usepackage{clock} 

\newcommand{\Em}{\mathbb M}

\providecommand{\cal}{\mathcal}
\renewcommand{\Bbb}{\mathbb}
\renewcommand{\frak}{\mathfrak}
\newenvironment{pf}{\begin{proof}}{\end{proof}}

\newcommand{\Aa}{{\Bbb{A}}}
\newcommand{\Aaa}{{\cal{A}}}

\newcommand{\Dee}{{\cal{D}}}
\newcommand{\Ef}{{\cal{F}}}

\newcommand{\Yu}{{\cal{U}}}

\newcommand{\Emm}{{\frak{M}}}

\newcommand{\Nat}{{\Bbb{N}}}
\newcommand{\Z}{{\mathbb Z}} 
\newcommand{\Qyu}{{\Bbb{Q}}}
\newcommand{\Err}{{\Bbb{R}}}

\newcommand{\lam}{{\lambda}}
\newcommand{\al}{\alpha}

\newcommand{\sig}{\sigma}
\newcommand{\eps}{\varepsilon}
\renewcommand{\phi}{\varphi}
\renewcommand{\rho}{\varrho}

\newcommand{\rest}{\restriction}
\newcommand{\unii}{\mathbb I}
\newcommand{\ntr}{{n\in\omega}}

\newcommand{\loe}{\leqslant}
\newcommand{\goe}{\geqslant}

\newcommand{\subs}{\subseteq}
\newcommand{\sups}{\supseteq}
\newcommand{\nnempty}{\ne\emptyset}
\newcommand{\argum}{\:\cdot\:}
\newcommand{\ovr}{\overline}
\renewcommand{\iff}{\Longleftrightarrow}

\newcommand{\cl}{\operatorname{cl}}

\newcommand{\w}{\operatorname{w}}
\newcommand{\dens}{\operatorname{dens}}

\newcommand{\conv}{\operatorname{conv}}

\newcommand{\G}{{\mathbb G}}

\newcommand{\id}[1]{{\operatorname{i\!d}_{#1}}} 
\newcommand{\cf}{\operatorname{cf}}
\newcommand{\dom}{\operatorname{dom}}
\newcommand{\Dom}{\operatorname{Dom}}

\newcommand{\suppt}{\operatorname{suppt}}

\newcommand{\liminv}{\varprojlim}

\newcommand{\oraz}{\qquad\text{and}\qquad}

\newcommand{\poset}{{\Bbb{P}}}

\newcommand{\Es}{{\cal{S}}}

\newcommand{\join}{\vee}

\newcommand{\Lev}{\operatorname{Lev}}
\newcommand{\Ht}{\operatorname{ht}}

\newcommand{\by}{/}

\newtheorem{tw}{Theorem}[section]
\newtheorem{wn}[tw]{Corollary}
\newtheorem{lm}[tw]{Lemma}
\newtheorem{prop}[tw]{Proposition}
\newtheorem{claim}[tw]{Claim}
\theoremstyle{definition}

\newtheorem{ex}[tw]{Example}

\theoremstyle{remark}

\newcommand{\set}[1]{\{#1\}}
\newcommand{\setof}[2]{\{#1\colon #2\}}

\newcommand{\sett}[2]{\{#1\}_{#2}}
\newcommand{\sn}[1]{\{#1\}} 
\newcommand{\dn}[2]{\{#1,#2\}} 
\newcommand{\pair}[2]{\langle #1, #2 \rangle} 
\newcommand{\triple}[3]{\langle #1, #2, #3 \rangle} 
\newcommand{\map}[3]{#1\colon #2 \to #3} 
\newcommand{\img}[2]{#1[#2]} 
\newcommand{\inv}[2]{{#1}^{-1}[#2]} 

\newcommand{\Cantor}{2^\omega}

\newcommand{\dpower}[2]{[#1]^{#2}}
\newcommand{\truj}[1]{{\dpower{#1}{\loe2}}}

\newcommand{\fra}{Fra\"iss\'e}
\newcommand{\jon}{J\'onsson}
\newcommand{\frajon}{\fra-\jon}
\newcommand{\U}{\mathbb U}

\newcommand{\Ama}{{\mathcal A}} 

\providecommand{\nat}{\omega}

\newcommand{\ciag}[1]{{\sett{{#1}_n}{\ntr}}}

\newcommand{\dualcat}[1]{{\ensuremath{{#1}^{\operatorname{op}}}}}

\newcommand{\Arr}{\operatorname{Arr}}
\newcommand{\iso}{\approx}
\newcommand{\uzup}[2]{{{\operatorname{Seq}}_{<{#1}}{\left(#2\right)}}}
\newcommand{\uzupiso}[2]{{{\operatorname{Seq}}^{\rm{iso}}_{<{#1}}{\left(#2\right)}}}
\newcommand{\uzuple}[2]{{{\operatorname{Seq}}_{\loe{#1}}{\left(#2\right)}}}

\newcommand{\ciagi}[1]{\sig{#1}}

\newcommand{\norm}[1]{\|#1\|}

\newcommand{\abs}[1]{|#1|}

\newcommand{\fK}{{\mathfrak{K}}}
\newcommand{\fL}{{\mathfrak{L}}}
\newcommand{\fM}{{\mathfrak{M}}}
\newcommand{\fC}{{\mathfrak{C}}}
\newcommand{\fS}{{\mathfrak{S}}}
\newcommand{\fT}{{\mathfrak{T}}}
\newcommand{\fU}{{\mathfrak{U}}}
\newcommand{\fB}{{\mathfrak{B}}}
\newcommand{\fF}{{\mathfrak{F}}}

\newcommand{\bS}{{\mathbb{S}}}

\newcommand{\cmp}{\circ} 

\newcommand{\komp}{\ensuremath{\mathfrak C\mathfrak o\mathfrak m\mathfrak p}} 
\newcommand{\rp}[1]{{\ddag{#1}}} 

\newcommand{\sets}{\ensuremath{\mathfrak S\mathfrak e\mathfrak t}} 
\newcommand{\banach}{\ensuremath{\mathfrak B}} 
\newcommand{\sban}{\ensuremath{\banach_{\aleph_0}}} 
\newcommand{\sbaniso}{\ensuremath{\banach_{\aleph_0}^{\operatorname{iso}}}} 
\newcommand{\LO}{\ensuremath{\mathfrak L\!\mathfrak O}}
\newcommand{\Quna}{\ensuremath{\Qyu_{\omega_1}}} 

\newcommand{\Pel}{\mathbf P}

\newcommand{\cont}{\ensuremath{\mathfrak c}}

\newcommand{\trees}{\ensuremath{\mathfrak{T}_2}}

\newcommand{\wek}[1]{{\vec{#1}}}

\newcommand{\C}{{\ensuremath\mathcal C}} 

\newcommand{\separator}{\begin{center} \leafright \leafright \leafright \decotwo \leafleft \leafleft\leafleft
\end{center}}

\newcommand{\fun}[3]{{\mathfrak{f}\left(#1,#2, #3 \right)}}

\newcommand{\qlim}{{\mathcal L}}
\newcommand{\cocones}[1]{\Delta\left(#1\right)}
\newcommand{\coconesiso}[1]{\Delta^{\rm{iso}}\left(#1\right)}

\title{\fra\ sequences: category-theoretic approach to universal homogeneous structures}
\author{
{\sc Wies{\l}aw Kubi\'s}\footnote{Research partially supported by the GA\v CR grant P
201/12/0290 (Czech Republic)}\\ \\
{\small Institute of Mathematics, Academy of Sciences of the Czech Republic}\\
{\small \v{Z}itn\'a 25, 115 67 Praha 1, Czech Republic}\\
{\small\texttt{kubis@math.cas.cz}}\\
and \\
{\small Institute of Mathematics,}
{\small Jan Kochanowski University}\\
{\small \'Swietokrzyska 15, 25-406 Kielce, Poland}
}

\date{\today}

\makeindex

\begin{document}

\maketitle

\begin{abstract}
We develop a category-theoretic framework for universal homogeneous objects, with some applications in the theory of Banach spaces, linear orderings, and in the topology of compact  Hausdorff spaces.

\ 


\end{abstract}

\tableofcontents

\section{Introduction}\label{rognowgrow}

When studying a given class of structures, it is an interesting and important issue to find a ``special" object that is universal for the class, namely, every other object ``embeds" into this special one. A better property of this special object would be some sort of homogeneity, namely, that every isomorphism between two ``small" substructures extends to an automorphism of the special object (as we shall see later, the exact meaning of ``small" will depend on the class under consideration).
Probably one of the earliest results in this spirit was Cantor's theorem on the uniqueness of the set of rational numbers among all countable linear orders.
More precisely, $\Qyu$ is the only (up to isomorphism) countable linear order which contains all other countable linear orders and every isomorphism between two finite subsets extends to an automorphism of $\Qyu$.
A typical statement of Cantor's theorem is different: the extension property is equivalent to saying that $\Qyu$ is order dense and has no end-points.
The argument used in the proof is now called the ``back-and-forth" method: an automorphism extending a given finite isomorphism is constructed inductively, interchanging the domain and the co-domain at each step.

It was Urysohn who found an analogue of Cantor's back-and-forth argument\footnote{
It seems that the back-and-forth argument was actually developed by Huntington (1904) and Hausdorff (1914), see \url{http://en.wikipedia.org/wiki/Cantor_back-and-forth_method}.}
for metric spaces. Namely, he found a complete separable metric space $\U$ which has very similar properties to the rationals, now considering isometries between finite subsets. 
Urysohn's work \cite{Urysohn} was actually forgotten for many years; it received some significant attention in the end of the 20th century due to problems in topological dynamics; perhaps the most representative work is~\cite{KPT}.

One of the most important works on this subject, independent of Urysohn, was done by Roland \fra~\cite{F1} in 1954. This is a model-theoretic approach to the back-and-forth argument, which can be partially applied also to the case of Urysohn's metric space.
Roughly speaking, \fra\ considers a class $\fK$ of finite (or, at least, finitely generated) models of some first order language.
The class should have the amalgamation property, that is, each two embeddings of the same model $Z\in \fK$ can be extended to a further embedding into a bigger model in $\fK$. In other words, given two embeddings $\map fZX$, $\map gZY$, there should exist a model $W\in \fK$ and embeddings $\map {f'}XW$ and $\map {g'}YW$ satisfying $f' \cmp f = g' \cmp g$.
If such a class has only countably many isomorphic types and each two models embed into a common one (the joint embedding property), then there exists a countable model $\U$ that can be represented as the union of a chain of models from the class $\fK$, contains isomorphic copies of all models in $\fK$ and has the following strong homogeneity property: every isomorphism between submodels of $\U$ which are in $\fK$ extends to an automorphism of $\U$. Furthermore, the model $\U$ is unique up to isomorphism.
It is often called the \emph{\fra\ limit} of the class $\fK$.

This line of investigation was further continued by \jon~\cite{Jon} and Morley \& Vaught~\cite{MorleyVaught} (see also Yasuhara~\cite{Yasuhara}), where uncountable classes of models were studied.
Besides the amalgamation property, typical cardinal-arithmetic assumption $\kappa = 2^{<\kappa}$ is needed for the existence of the universal homogeneous structure of cardinality $\kappa$.
One has to mention a curious independent work of Trnkov\'a~\cite{Trnkova} with metamathematical results on universal categories, in the setting of Bernays-G\"odel set theory.
One of the main tools is the amalgamation property for certain classes of categories, treated just as first-order structures.

All the authors cited above assume that the class of structures is closed under unions of chains of length less than the size of the universal homogeneous structure (in the last case, the union of any chain of small categories is a category).
One of the objectives of this work is to relax this assumption and to make the theory general enough for capturing new cases and obtaining new examples of universal homogeneous objects.

We believe that category theory is the proper language for \fra-\jon\ limits. 
In fact, this has already been confirmed by the works of Droste \& G\"obel~\cite{DrGoe93, DrGoe92}, where the authors consider some categories of first-order structures with special types of embeddings, obtaining new applications in algebra and theoretical computer science.
In this context, one has to mention a recent work of Pech \& Pech~\cite{Pech} where the authors, based on the results of Droste \& G\"obel, develop the theory of \fra\ limits in comma categories, leading to universal homomorphisms and universal retractions.
Finally, Irwin \& Solecki~\cite{IrSo} presented a variant of \fra\ theory with reversed arrows, i.e., epimorphisms of finite structures instead of embeddings.
By this way, they obtained an interesting new characterization of the \emph{pseudo-arc}---a certain connected compact metric space, never associated to \fra\ limits before.

There is no doubt that universal structures with strong homogeneity properties can be discovered or identified in various areas of pure mathematics, theoretical computer science and even mathematical physics (see \cite{Droste}).
Model-theoretic \fra\ limits are nowadays important objects of study in combinatorics, permutation group theory and topological dynamics, see Macpherson's survey~\cite{Macpherson} for more information and further references. 
Category theory brings much more freedom for dealing with \fra\ limits, offering the possibility of constructing new objects from old and eliminating superfluous assumptions.
It has been demonstrated in \cite{DrGoe93, DrGoe92, Pech} that category-theoretic approach brings new important examples of universal homogeneous objects and the work \cite{IrSo} shows that one of the simplest constructions in category theory, namely, passing to the opposite category, leads to new and somewhat surprising examples.
Actually, one of the author's construction of a universal pre-image for a certain class of compact linearly ordered spaces~\cite{K_classR} turns out to be the \fra\ limit of a category whose arrows are increasing quotient maps.

Summarizing, category-theoretic approach offers powerful tools for constructing new universal homogeneous objects as well as identifying existing objects, discovering their homogeneity properties.

Here a notational issue has to be pointed out.
Namely, in category theory the notion of a ``universal object" is totally different from the notion of a ``universal structure" in model theory and related areas.
In order to avoid this confusion, we shall replace the adjective ``universal" by ``cofinal" in the latter case.
Namely, an object $U$ is defined to be \emph{cofinal} for a class $\fK$ if every object from $\fK$ embeds into $U$, where ``embedding" will be just an arrow of a category under consideration.

\separator

As mentioned above, the aim of this work is extending the theory of \frajon\ limits in the framework of pure category theory.
The key point of our approach is dealing with sequences instead of (co-)limits, where a \emph{sequence} is nothing but a functor from an ordinal into the category.
The crucial notion is that of a \emph{\fra\ sequence}, a sequence which is supposed to ``converge" to a universal (cofinal, in our terminology) homogeneous object in a bigger category, where homogeneity is meant with respect to the original category.
In other words, we deal with the base category of ``small" objects and we use sequences for encoding the category of ``large" objects.
This approach is similar in spirit to the idea of forcing in set theory, where one deals with approximations of a generic object, often working in the base model instead of its forcing extension.

The paper is organized as follows.
Section~\ref{SectPreLims} contains basic notation, including a very short treatment of categories of sequences.
Section~\ref{SecFraMejnone} presents the main ideas of this work: \fra\ sequences, their existence and properties.
Section~\ref{SectFraJonxqlims} deals with discontinuity: We describe a way of dealing with \fra\ sequences that are continuous with respect to some quasi-limiting operator.
We also introduce the notion of an amalgamation structure, specializing the result on cofinality of a \fra\ sequence in the category of sequences.
Section~\ref{SectExmplsFirstSeries} contains selected examples of classes of categories with \fra\ sequences, demonstrating the ideas of previous sections.
In particular, we discuss categories of functors (diagrams), monoids, and certain categories of trees.
We also present an example of a category with many incomparable \fra\ sequences.
Section~\ref{reterpeteairs} deals with categories of embedding-projection pairs, aiming at applications to Banach space theory and topology of compact spaces.

Finally, Section~\ref{SectAppsEPetsoon} contains applications, the first four of them are equivalent to the Continuum Hypothesis.
One of the results there, which actually could fit into the model-theoretic framework, is the existence of a unique Banach space of density $\aleph_1$ that is isometrically universal for this class of spaces and isometrically homogeneous for separable Banach spaces\footnote{There is a recent work \cite{ACCGM} elaborating the method of constructing Gurari\u\i-like Banach spaces using amalgamations. The article~\cite{ACCGM} was however partially inspired by our result, with a reference to an early draft of this work.}\label{Ffootone}.
It turns out that such a Banach space has not been discovered before.
Another result is a Banach space $\Pel$, complementably universal for the class of spaces with monotone Schauder bases of length $\loe \omega_1$.
More precisely, every Banach space from this class is linearly isometric to a 1-complemented subspace of $\Pel$.
Yet another result shows the existence of a Banach space of density $\aleph_1$ that is complementably universal for the class of Banach spaces with the so-called \emph{projectional resolution of the identity}. This class of Banach spaces had been studied extensively by several authors, in particular, in connections with renorming theory, see the monographs~\cite{DGZ} and~\cite{Fabian} for further references.
We also present a dual result in topology, for the class of \emph{Valdivia compact spaces} that can be described in the language of retractions onto metrizable compacta.
We show that there exists a Valdivia compact $K$ of weight $\aleph_1$ such that all other Valdivia compacta of the same weight embed as retracts into $K$.
Both results are straightforward applications of Section~\ref{reterpeteairs}, by considering two natural categories of embedding-projection pairs.
Another result, now without any cardinal-arithmetic assumptions, is the existence of a sequence of continuous functions on the Cantor set which is universal both for quotient maps and for topological embeddings.
This is obtained by adapting the ideas of embedding-projection pairs to the category of finite nonempty sets.
We finally describe, again without any cardinal-arithmetic assumptions, a linearly ordered set of size $\aleph_1$ which could serve as a natural generalization of the rationals, with homogeneity property involving increasing embeddings and projections at the same time.
Its topological counterpart identifies as a \fra\ limit the author's example~\cite{K_classR} of a universal pre-image for the class of linearly ordered Valdivia compacta, thus showing its homogeneity and uniqueness with respect to these properties.

\section{Preliminaries}\label{SectPreLims}

Categories will usually be denoted by letters $\fK$, $\fL$, $\fM$, etc.
\index{category}\index{functor}
All functors considered here are assumed to be covariant, unless otherwise specified.
Let $\fK$ be a category. We shall write ``$a\in\fK$" for ``$a$ is an object of $\fK$". Given $a,b\in\fK$, we shall denote by $\fK(a,b)$ the set of all $\fK$-morphisms from $a$ to $b$. The composition of two compatible arrows $f$ and $g$ will be denoted by $g\cmp f$.
\index{subcategory}\index{subcategory!-- full}
A \emph{subcategory} of $\fK$ is a category $\fL$ such that each object of $\fL$ is an object of $\fK$ and each arrow of $\fL$ is an arrow of $\fK$ (with the same domain and co-domain). We write $\fL\subs \fK$.
Recall that a subcategory $\fL$ of $\fK$ is \emph{full} if $\fL(a,b)=\fK(a,b)$ for every objects $a,b\in\fL$.
We say that $\fL$ is \emph{cofinal} in $\fK$
if for every object $x\in\fK$ there exists an object $y\in\fL$ such that $\fK(x,y)\nnempty$.
\index{category!-- cofinal}\index{cofinal subcategory}
The opposite category to $\fK$ will be denoted by $\dualcat \fK$. That is, the objects of $\dualcat\fK$ are the objects of $\fK$ and all arrows are reversed, i.e. $\dualcat\fK(a,b)=\fK(b,a)$ for every $a,b\in \fK$.
\index{category!-- ordered}
Recall that a category $\fK$ is \emph{ordered} if $|\fK(x,y)| \loe 1$ and $\fK(x,y) \nnempty \ne \fK(y,x)$ implies $x = y$ for every $\fK$-objects $x,y$.
Removing the last condition we get the notion of a \emph{quasi-ordered category}.
Every (not necessarily ordered) category induces a partial order $\loe$ on the objects of $\fK$ defined by the formula $\fK(x,y)\nnempty$ iff $x \loe y$.
Every partially ordered set $\pair P \loe$ can be viewed as an ordered category $\fK_P$ with $P$ the class of objects and the class of arrows defined by $\fK_P(x,y) = \sn {\pair xy}$ whenever $x\loe y$ and $\fK_P(x,y) = \emptyset$ otherwise.
In particular, ordinals treated as well ordered sets are important examples of ordered categories.

Let $\fK$ be a category. We say that $\fK$ has the \emph{amalgamation property}
\index{amalgamation property}
if for every $a,b,c\in\fK$ and for every morphisms $f\in\fK(a,b)$, $g\in\fK(a,c)$ there exist $d\in\fK$ and morphisms $f'\in\fK(b,d)$ and $g'\in\fK(c,d)$ such that $f'\cmp f = g'\cmp g$. If, additionally, for every arrows $f''$, $g''$ such that $f''\cmp f=g''\cmp g$ there exists a unique arrow $h$ satisfying $h\cmp f' = f''$ and $h\cmp g' = g''$ then the pair $\pair{f'}{g'}$ is a \emph{pushout} of $\pair fg$.
\index{pushout}
Reversing the arrows, we define the \emph{reversed amalgamation} and the \emph{pullback}.
\index{pullback}
\index{reversed amalgamation}\index{amalgamation!-- reversed}
We say that $\fK$ is \emph{directed}
if for every $a,b\in\fK$ there exists $g\in\fK$ such that both sets $\fK(a,g)$, $\fK(b,g)$ are nonempty.
In model theory, where the arrows are embeddings, this is usually called the \emph{joint embedding property}.
\index{joint embedding property}\index{category!-- directed}

\index{sequence}\index{inductive sequence}\index{sequence!-- inductive}
Fix a category $\fK$ and fix an ordinal $\delta>0$. An \emph{inductive $\delta$-sequence}
in $\fK$ is 
formally a covariant functor from $\delta$ (treated as a poset category) into $\fK$. In other words, it could be described as a pair of the form $\pair {\sett{a_\xi}{\xi<\delta}} {\sett{a^\eta_\xi}{\xi<\eta<\delta}}$, where $\delta$ is an ordinal, $\setof{a_\xi}{\xi<\delta}\subs\fK$ and $a^\eta_\xi\in\fK(a_\xi,a_\eta)$ are such that $a^\rho_\eta \cmp a_\xi^\eta = a^\rho_\xi$ for every $\xi<\eta<\rho<\delta$. We shall denote such a sequence shortly by $\wek a$. The ordinal $\delta$ is the \emph{length} of $\wek a$.
\index{sequence!-- length of}

\index{category!-- $\kappa$-complete}
Let $\kappa$ be an infinite cardinal. A category $\fK$ is \emph{$\kappa$-complete}
if all inductive sequences of length $<\kappa$ have co-limits in $\fK$. Every category is $\aleph_0$-complete, since the co-limit of a finite sequence is its last object.
\index{category!-- $\kappa$-bounded}
A category $\fK$ is \emph{$\kappa$-bounded} if for every inductive sequence $\wek x$ in $\fK$ of length $\lam<\kappa$ there exist $y\in\fK$ and a co-cone of arrows $\sett{y_\al}{\al<\lam}$ such that $\map {y_\al}{x_\al}y$ and $y_\beta \cmp x_\al^\beta = y_\al$ for every $\al<\beta<\lam$. 
Obviously, every $\kappa$-complete category is $\kappa$-bounded. We shall write ``$\sig$-complete" and ``$\sig$-bounded" for ``$\aleph_1$-complete" and ``$\aleph_1$-bounded" respectively.

\index{$\Dom(\argum)$}
\index{dominating family}\index{family!-- dominating}
We shall need the following notion concerning families of arrows. Fix a family of arrows $\Ef$ in a given category $\fK$. We shall write $\Dom(\Ef)$ for the set $\setof{\dom(f)}{f\in\Ef}$.
We say that $\Ef$ is \emph{dominating} in $\fK$ if the family of objects $\Dom(\Ef)$ is cofinal in $\Ef$ and moreover for every $a\in\Dom(\Ef)$ and for every arrow $\map fax$ in $\fK$ there exists an arrow $g$ in $\fK$ such that $g\cmp f\in \Ef$.

For all undefined category-theoretic notions we refer to Mac Lane~\cite{MacLane} or Johnstone~\cite{Johnstone}.

\subsection{Categories of sequences}

Fix a category $\fK$ and denote by $\uzup\kappa\fK$ the class of all sequences in $\fK$ which have length $<\kappa$. We shall write $\uzuple\kappa\fK$ instead of $\uzup{\kappa^+}\fK$ and $\ciagi\fK$ instead of $\uzuple{\aleph_0}\fK$.
\index{$\uzup\kappa{\argum}$}\index{$\uzuple\kappa{\argum}$}\index{$\ciagi{\argum}$}
We would like to turn $\uzup\kappa\fK$ into a category in such a way that an arrow from a sequence $\wek a$ into a sequence $\wek b$ induces an arrow from $\lim\wek a$ into $\lim \wek b$, whenever $\fK$ is embedded into a category in which sequences $\wek a$, $\wek b$ have co-limits.

\index{transformation of sequences}\index{sequence!-- transformation of}
Fix two sequences $\wek a$ and $\wek b$ in a given category $\fK$. Let $\lam=\dom(\wek a)$, $\rho=\dom(\wek b)$. A \emph{transformation} from $\wek a$ to $\wek b$ is, by definition, a natural transformation from $\wek a$ into $\wek b\cmp \phi$, where $\map\phi\lam\rho$ is an order preserving map (i.e. a covariant functor from $\lam$ into $\rho$, treated as ordered categories).

In order to define an arrow from $\wek a$ to $\wek b$ we need to identify some transformations. Fix two natural transformations $\map F{\wek a}{\wek b\cmp \phi}$ and $\map G{\wek a}{\wek b\cmp \psi}$. We shall say that $F$ and $G$ are \emph{equivalent} if the following conditions hold:
\begin{enumerate}
	\item[(1)] For every $\al$ there exists $\beta\goe \al$ such that $\phi(\al)\loe\psi(\beta)$ and $b_{\phi(\al)}^{\psi(\beta)} \cmp F(\al) = G(\beta) \cmp a_\al^\beta$.
	\item[(2)] For every $\al$ there is $\beta\goe\al$ such that $\psi(\al)\loe\phi(\beta)$ and $b_{\psi(\al)}^{\phi(\beta)} \cmp G(\al) = F(\beta) \cmp a_\al^\beta$.
\end{enumerate}
It is rather clear that this defines an equivalence relation, which is actually a congruence on the category of transformations.
Every equivalence class of this relation will be called an \emph{arrow} (or \emph{morphism}) from $\wek a$ to $\wek b$.
It is easy to check that this indeed defines a category structure on all sequences in $\fK$.
Formally, this is the quotient category with respect to the equivalence relation described above.
The identity arrow of $\wek a$ is the equivalence class of the identity natural transformation $\map {\id{\wek a}}{\wek a}{\wek a}$.
\index{sequence!-- identity transformation of}
\index{morphism of sequences}\index{sequence!-- morphism of}
\index{arrow of sequences}\index{sequence!-- arrow of}

Categories of sequences are special cases (or rather ``parts") of more general categories called \emph{Ind-completions}, see Chapter VI of Johnstone's monograph \cite{Johnstone}.

We shall later need the following two facts.

\begin{lm}\label{LmDigprz}
Let $\fK$ be a category, $\kappa \goe \aleph_0$ a regular cardinal, and let $\wek{\wek x}$ be a $\kappa$-sequence in $\uzuple \kappa \fK$ whose elements are sequences of length $\kappa$.
Then $\wek{\wek x}$ has co-limit in $\uzuple \kappa \fK$.
\end{lm}

\begin{pf}
Refining inductively each $\wek x_\al$ ($\al < \kappa$) to a cofinal subsequence, we may assume that all bonding maps are natural transformations.
Now look at $\wek{\wek x}$ as a functor from $\kappa \times \kappa$ into $\fK$ and let $\wek y$ be the diagonal sequence.
Then $\wek y$ is easily seen to be the co-limit of $\wek{\wek x}$.
\end{pf}

\begin{prop}
Assume $\fK$ is a $\kappa$-complete category, where $\kappa = \cf \kappa \goe \aleph_0$.
Then $\uzuple \kappa \fK$ is $\kappa^+$-complete.
\end{prop}

\begin{pf}
The fact that $\fK$ is $\kappa$-complete implies that so is $\uzuple \kappa \fK$.
Finally, $\kappa^+$-comp\-lete\-ness follows from Lemma~\ref{LmDigprz}.
\end{pf}

\subsection{Partial orders and trees}

By a \emph{tree}
\index{tree}
we mean a partially ordered set $\pair T\loe$ which is a meet semilattice, i.e. every two elements of $T$ have the greatest lower bound, and for every $t\in T$ the interval $\setof{x\in T}{x<t}$ is well ordered. Every tree $T$ has a single minimal element $0_T$, called the \emph{root} of $T$.
\index{tree!-- root of}
\index{tree!-- immediate successor}
An \emph{immediate successor} of $t\in T$ is an element $s>t$ such that no $x\in T$ satisfies $t<x<s$.
The set of all immediate successors of $t$ will be denoted by $t^+$.
\index{subtree}
A \emph{subtree} of a tree $T$ is a subset $S\subs T$ which is again a semilattice, possibly with a different infimum.
\index{tree!-- binary}
A tree $T$ is \emph{binary} if every $t\in T$ has at most two immediate successors. We shall denote by $\max T$ the set of all maximal elements of $T$. A tree $T$ is \emph{bounded}
\index{tree!-- bounded}\index{bounded tree}
if for every $x\in T$ there is $t\in \max T$ such that $x\loe t$. Recall that an \emph{initial segment} of a poset $\pair T\loe$
\index{tree!-- initial segment}\index{initial segment}\index{initial subtree}
is a subset $A$ of $T$ satisfying $\setof{x\in T}{x\loe t}\subs A$ for every $t\in A$. 
Given a tree $T$ and $t\in T$, we define the \emph{level} of $t$ in $T$ as the order type of the interval $[0,t)=\setof{s\in T}{s < t}$.
We denote it by $\Lev_T(t)$.
The \emph{height} of $T$, $\Ht(T)$ is the minimal ordinal $\al$ such that no $t\in T$ satisfies $\Lev_T(t) = \al$.
\index{tree!-- level}\index{tree!-- height}

For all undefined set-theoretic notions we refer to Jech~\cite{Jech} or Kunen~\cite{Kunen}.

\section{\fra\ sequences}\label{SecFraMejnone}

Below we introduce the key notion of this work.

\index{\fra\ sequence}\index{sequence!-- \fra}
Let $\fK$ be a category and let $\kappa$ be a cardinal. A \emph{\fra\ sequence} of length $\kappa$ in $\fK$ (briefly: a \emph{$\kappa$-\fra\ sequence}) is an inductive sequence $\wek u$ satisfying the following conditions:
\begin{enumerate}
	\item[(U)] For every $x\in \fK$ there exists $\xi<\kappa$ such that $\fK(x,u_\xi)\nnempty$.
	\item[(A)] For every $\xi<\kappa$ and for every arrow $f\in\fK(u_\xi,y)$, where $y\in\fK$, there exist $\eta\goe\xi$ and $g\in\fK(y,u_\eta)$ such that $u_\xi^\eta= g\cmp f$.
\end{enumerate}
An inductive sequence satisfying (U) will be called \emph{$\fK$-cofinal}.
\index{sequence!-- cofinal}
\index{sequence!-- $\fK$-cofinal}
More generally, a collection $\Yu$ of objects of $\fK$ is \emph{$\fK$-cofinal} if for every $x\in \fK$ there is $u\in\Yu$ such that $\fK(x,u)\nnempty$.
\index{family!-- cofinal}
Condition (A) will be called \emph{amalgamation property}.
\index{sequence!-- amalgamation}

\subsection{Basic properties}

Let $\wek v$ be a $\kappa$-sequence in a category $\fK$. We say that $\wek v$ has the \emph{extension property}
\index{sequence!-- extension property}
if the following holds:
\begin{itemize}
	\item[(E)] For every arrows $\map fab$, $\map g a {v_\al}$ in $\fK$, where $\al<\kappa$, there exist $\beta\goe\al$ and an arrow $\map hb{v_\beta}$ such that $v_\al^\beta \cmp g = h \cmp f$.
\end{itemize}
Clearly, this condition implies (A).

\begin{prop}\label{krakrakra}
Let $\wek u$ be a $\kappa$-\fra\ sequence in a category $\fK$. Then $\fK$ is directed. Moreover, the following conditions are equivalent:
\begin{enumerate}
	\item[(a)] $\wek u$ has the extension property.
	\item[(b)] $\fK$ has the amalgamation property.
\end{enumerate}
\end{prop}

\begin{pf}
The first statement is trivial.
Assume (a) and fix arrows $\map fzx$, $\map gzy$. Using (U), find $\map hx{u_\al}$, $\al<\kappa$. Using (E), find $\beta\goe\al$ and $\map ky{u_\beta}$ such that $k\cmp g = u_\al^\beta \cmp h\cmp f$. Thus (b) holds.

Finally, assume (b) and fix arrows $\map fab$ and $\map ga{u_\al}$, $\al<\kappa$. Using (b), find arrows $\map{f_1}bw$ and $\map{g_1}{u_\al}w$ so that $f_1\cmp f = g_1\cmp g$. Using (A) for the sequence $\wek u$, we find $\beta\goe\al$ and $\map hw{u_\beta}$ so that $h\cmp g_1 = u_\al^\beta$. Thus $(h\cmp f_1)\cmp f=h\cmp g_1\cmp g=u_\al^\beta\cmp g$, which shows that (a) holds.
\end{pf}

\begin{prop}\label{laskacervena}
Assume $\fK$ is a directed category. Then every sequence in $\fK$ satisfying condition (A) is \fra.
\end{prop}

\begin{pf} 
Let $\wek u$ be a sequence in $\fK$ satisfying (A). Fix $x\in\fK$.
Since $\fK$ is directed, there exist $w\in\fK$ and arrows $\map f{u_0}w$, $\map gxw$. Using (A), we find an arrow $\map hw{u_\xi}$ such that $h\cmp f=u_0^\xi$. Thus $\fK(x,u_\xi)\nnempty$, which shows (U).
\end{pf}

\begin{prop}\label{laskanebeska}
Let $\fK$ be a category, let $\wek u$ be an inductive sequence of length $\kappa$ in $\fK$ and let $S\subs \kappa$ be unbounded in $\kappa$.
\begin{enumerate}
	\item[(a)] If\/ $\wek u$ is a \fra\ sequence in $\fK$ then $\wek u\rest S$ is \fra\ in $\fK$.
	\item[(b)] If\/ $\fK$ has the amalgamation property and $\wek u\rest S$ is a \fra\ sequence in $\fK$ then so is $\wek u$.
\end{enumerate}
\end{prop}

\begin{pf}
Assume $\wek u$ is a \fra\ sequence. Then $\wek u\rest S$ clearly satisfies (U). In order to check (A), fix $\map f{u_\xi}y$ with $\xi\in S$. Then $u_\xi^\eta=g\cmp f$ for some arrow $g$ and for some $\eta\goe\xi$. Since $S$ is unbounded in $\kappa$, there is $\al\in S$ such that $\al\goe\eta$. Then $u_\xi^\al = u_\eta^\al\cmp g\cmp f$, which shows that $\wek u\rest S$ satisfies (A).

Now assume $\wek u\rest S$ is a \fra\ sequence. Clearly, $\wek u$ satisfies (U). Fix $\map f{u_\xi}y$, $\xi<\kappa$. Find $\al\in S$ with $\al\goe\xi$. Using the amalgamation property of $\fK$, find $\map {f'}{u_\al}z$ such that the diagram
$$\xymatrix{
{u_\al} \ar[r]^{f'} & z\\
{u_\xi} \ar[u]^{u_\xi^\al}\ar[r]^f & y \ar[u]_g
}$$
commutes for some arrow $g$ in $\fK$. Now, using (A) for $\wek u\rest S$, we can find $\beta\in S$ such that $\beta\goe\al$ and $h\cmp f' = u_\al^\beta$ holds for some $\map hz{u_\beta}$. This shows that $\wek u$ satisfies (A).
\end{pf}

Recall that we consider sequences up to equivalence with respect to the relation defined in Section~\ref{SectPreLims}.
Thus, we need to show that being a \fra\ sequence does not depend on the representation.
This is indeed true, assuming amalgamations:

\begin{prop}\label{Prlaskaamodr}
Let $\fK$ be a category with the amalgamation property and let $\map {\wek g}{\wek u}{\wek v}$, $\map {\wek h}{\wek v}{\wek u}$ be transformations of sequences such that $\wek g \cmp \wek h$ is equivalent to the identity of $\wek v$.
If $\wek u$ is a \fra\ sequence then so is $\wek v$.
\end{prop}

\begin{pf}
Condition (U) is obvious.
Fix a $\fK$-arrow $\map e {v_\al} y$ and assume $\map {h_\al} {v_\al} {u_{\al'}}$.
Using the amalgamation property, find $\map k {u_{\al'}} w$ and $\map \ell y w$ such that $\ell \cmp e = k \cmp h_\al$.
Using the fact that $\wek u$ is \fra, find $\beta \goe \al'$ and a $\fK$-arrow $f$ such that $f \cmp k = u_{\al'}^\beta$.
Finally, if $\beta'\goe \al$ is such that $\map {g_\beta}{u_\beta}{v_{\beta'}}$, then $v_\al^{\beta'} = g_\beta \cmp u_{\al'}^\beta \cmp h_\al = g_\beta \cmp f \cmp \ell \cmp e.$
This shows (A).
\end{pf}

A \fra\ sequence can possibly have finite length. In that case, by Proposition \ref{laskanebeska}(a), there is also a \fra\ sequence of length one---it is an object $u$ which is cofinal in $\fK$ and which satisfies the following version of (A): given $f\in\fK(u,x)$, where $x\in \fK$, there exists $g\in \fK(x,u)$ such that $g\cmp f=\id u$. We shall call $u$ a \emph{\fra\ object} in $\fK$.
\index{\fra\ object}
Given a \fra\ object $u$, the sequence $u\to u \to \dots$, where each arrow is identity, is a \fra\ sequence of length $\omega$. Thus, it follows from Theorem \ref{jednorodnost} below that a possible \fra\ object is unique, up to isomorphism. Below we give a direct proof of this fact.

\begin{prop}\label{asfasfqwrQWR}
Assume $u,v$ are \fra\ objects in a category $\fK$. Then $u\iso v$. If moreover all arrows in $\fK$ are monomorphisms then every arrow $\map fux$ is an isomorphism.
\end{prop}

\begin{pf}
Applying (U) for $v$, we find a morphism $\map{f_0} uv$ which, using (A) for $u$, has a left inverse $\map{g_0}vu$, i.e. $g_0\cmp f_0 = \id u$. Now, using (A) for $v$, we obtain an arrow $\map {f_1}uv$ such that $f_1\cmp g_0 = \id v$. Observe that
$$f_1 = f_1\cmp \id u = f_1\cmp (g_0\cmp f_0) = (f_1\cmp g_0)\cmp f_0 = \id v\cmp f_0 = f_0.$$
Hence $f_0\cmp g_0=\id v$, which shows that $f_0$ is an isomorphism. 

Finally, let $\map fux$ be a morphism in $\fK$. Again by (A), $f$ has a left inverse $\map gxu$. Assuming $g$ is a monomorphism, we deduce that $f\cmp g=\id x$, because $g\cmp (f\cmp g) = g\cmp \id x$. Thus $f$ is an isomorphism.
\end{pf}

We finish this subsection with the following stability result which is well-known for model-theoretic \fra\ limits.

\begin{tw}\label{ThmStbFrcsest}
Let $\fK$ be a category, $\kappa \goe \aleph_0$ a regular cardinal, and let $\map {\wek{\wek u}} \kappa {\uzuple \kappa \fK}$ be a sequence in $\uzuple \kappa \fK$ such that each $\wek u_\al$ is a $\kappa$-\fra\ sequence in $\fK$.
Then the co-limit of $\wek{\wek u}$ in $\uzuple \kappa \fK$ is a \fra\ sequence in $\fK$.
\end{tw}

Recall that by Lemma~\ref{LmDigprz}, the co-limit of $\wek {\wek u}$ indeed exists.

\begin{pf}
Since every cofinal subsequence of a \fra\ sequence is again \fra, we can inductively refine each $\wek u_\al$ so that the bonding maps become natural transformations. 
We look at $\wek{\wek u}$ as a ``two-dimensional sequence" in $\fK$ and we denote by $u_{\al,\beta}$ the $\beta$th object of the sequence $\wek u_\al$.
Now the co-limit of $\wek{\wek u}$ is the diagonal sequence $\wek v$, where $v_\al = u_{\al,\al}$ for $\al < \kappa$.
It is clear that $\wek v$ satisfies (U).
Fix a $\fK$-arrow $\map f {v_\xi} y$.
Since $\wek u_\xi$ is \fra, there are $\eta > \xi$ and a $\fK$-arrow $\map g y{u_{\xi,\eta}}$ such that $g \cmp f = u_{\xi,\xi}^{\xi,\eta}$.
Finally, let $h = u_{\xi,\eta}^{\eta,\eta} \cmp g$, where $u_{\xi,\eta}^{\eta,\eta}$ is the suitable component of the natural transformation from $\wek u_\xi$ to $\wek u_\eta$.
Then $h \cmp f = v_\xi^\eta$, which shows that $\wek v$ satisfies (A).
\end{pf}

\subsection{The existence}

We present below a simple yet useful criterion for the existence of a \fra\ sequence. In case of sequences of length $\loe\aleph_1$, this criterion becomes a characterization.

\index{dominating family of arrows}
Recall that a family of arrows $\Ef$ is \emph{dominating} in $\fK$ if it satisfies the following two conditions.
\begin{enumerate}
	\item[(D1)] The family $\Dom(\Ef)$ is cofinal in $\fK$, i.e. for every $x\in \fK$ there is $a\in\Dom(\Ef)$ such that $\fK(x,a)\nnempty$.
	\item[(D2)] Given $a\in\Dom(\Ef)$ and $\map f a y$ in $\fK$, there exist $\map g y b$ in $\fK$ such that $g\cmp f\in \Ef$.
\end{enumerate}

\begin{tw}[Existence]\label{jofjaiopf}
Let $\kappa$ be an infinite regular cardinal and let $\fK$ be a $\kappa$-bounded directed category with the amalgamation property. Assume further that $\Ef\subs\Arr(\fK)$ is dominating in $\fK$ and $|\Ef|\loe\kappa$.
Then there exists a \fra\ sequence of length $\kappa$ in $\fK$.
\end{tw}
\index{Existence Theorem}

\begin{pf}
Let $\Dom(\Ef)=\sett{a_\al}{\al<\kappa}$ and enumerate $\Ef\times \kappa$ as $\sett{\pair{f_\al}{i_\al}}{\al<\kappa}$ so that for each $p\in\Ef\times\kappa$ the set $\setof{\al<\kappa}{p=\pair{f_\al}{i_\al}}$ is unbounded in $\kappa$. We shall construct inductively a $\kappa$-\fra\ sequence $\wek u$, so that the following conditions are satisfied:
\begin{enumerate}
  \item[(i)] $u^\al_\eta\cmp u^\eta_\xi=u^\al_\xi$ for every $\xi<\eta<\al$ (i.e. $\wek u$ is indeed a sequence).
	\item[(ii)] $u_{\al}\in\Dom(\Ef)$ and $\fK(a_\al,u_{\al})\nnempty$.
	\item[(iii)] Given $\xi<\al$, if $\xi=i_\al$ and $\dom(f_\al)=u_\xi$ then there exists an arrow $h$ in $\fK$ such that $h\cmp f_\al=u_\xi^{\al}$.
\end{enumerate}
We start with $u_0=a_0$.
Assume that $\beta<\kappa$ is such that $u_\xi$ and $u_\xi^\eta$ have been constructed for all $\xi<\eta<\beta$. Using the fact that $\fK$ is $\kappa$-bounded, find $v\in\fK$ and $\map{j_\al}{u_\al}v$ such that $j_\xi = j_\eta\cmp u_\xi^\eta$ holds for every $\xi<\eta<\beta$. Using the fact that $\fK$ is directed, we may ensure that $\fK(a_\beta,v)\nnempty$.
Now, if $\map{f_\beta}{u_\xi}y$ and $\xi=i_\beta<\beta$ then, using amalgamation, we may find arrows $\map h{v}{w}$ and $\map gyw$ so that $g\cmp f_\beta = h\cmp j_\xi$ holds. Using (D1), we may further assume that $w\in \Dom(\Ef)$. Finally, set $u_\beta:=w$ and $u^\beta_\eta:=h\cmp j_\eta$ for $\eta<\beta$. It is clear that conditions (i) -- (iii) are satisfied.

It follows that the construction can be carried out. It remains to check that $\map{\wek u}\kappa\fK$ is a \fra\ sequence.
Condition (i) says that $\wek u$ is indeed an inductive sequence. Conditions (D1) and (ii) imply (U). In order to justify (A), fix $\xi<\kappa$ and $f\in \fK(u_\xi,x)$, where $x\in \fK$. We need to find $\al>\xi$ and an arrow $g$ so that $g\cmp f=u_\xi^\al$.
Since $u_\xi\in\Dom(\Ef)$, using (D2), we can find $g\in\Ef$ such that $g = k\cmp f$ for some arrow $k$. Now find 
$\al>\xi$ such that $f_\al=g$ and $i_\al=\xi$.
By (iii), $h\cmp g=u_\xi^{\al+1}$ for some arrow $h$. Hence $(h\cmp k)\cmp f=u_\xi^{\al+1}$, which completes the proof.
\end{pf}

The existence of a \fra\ sequence can also be proved by using a Baire category argument. We shall demonstrate it below for the countable case, where the classical Baire Category Theorem is invoked.

\begin{wn}
Let $\fK$ be a directed category with the amalgamation property. Assume that $\fK$ is dominated by a countable family of arrows. Then $\fK$ has an $\omega$-\fra\ sequence.
\end{wn}

\begin{pf}
Without loss of generality, we may assume that a countable dominating family $\Ef$ is a subcategory of $\fK$.
An $\omega$-sequence $\wek x$ in $\Ef$ may be regarded as a function from $\Delta = \setof{\pair mn}{m\loe n}$ into $\Ef$ satisfying the obvious conditions. Thus, the set $S$ of all $\omega$-sequences in $\Ef$ is a closed subspace of the Polish space $\Ef^\Delta$, endowed with the product topology.
Given an object $x$ in $\Ef$, let $U_x$ be the set of all $\wek x\in S$ for which there exists an arrow $x \to \wek x$.
Clearly, $U_x$ is open and dense in $S$.
Given an arrow $\map fab$ in $\Ef$ and $\ntr$, let
$$V_{f,n} = \setof{\wek x \in S}{x_n = a \implies (\exists\;m > n)(\exists\;g)\; g \cmp f = x_n^m}.$$
Again, $V_{f,n}$ is open and dense in $S$ (for the density one needs to use amalgamations).
Using the Baire Category Theorem, we can find a sequence $\wek u \in S$ that belongs to all the sets defined above.
It is easy to check, using the fact that $\Ef$ is dominating in $\fK$, that $\wek u$ is a \fra\ sequence in $\fK$.
\end{pf}

The argument above can be repeated in the general case, however one has to use the topological space $S^\kappa$, where $\kappa$ is a regular cardinal and the topology is generated by open sets of the form $V_s = \setof{x}{s\subs x}$, where $s\in S^{<\kappa}$. It is a standard and easy fact that in this topology, the intersection of $\kappa$ many open dense sets is dense in $S^\kappa$.

We shall sometimes need the following variant of Theorem~\ref{jofjaiopf}, involving a functor:

\begin{prop}\label{strangewef}
Let $\fK$ be a directed category with the amalgamation property. Let $\kappa\goe\aleph_0$ be a regular cardinal and let $\map \Phi\fK\fL$ be a functor satisfying the following condition:
\begin{enumerate}
	\item[$(\mathcal C)$] For every $\wek x\in\uzup\kappa\fK$ there exist $a\in\fK$ and a co-cone of arrows $\Ef$ in $\fK$ such that $\pair{\Phi(a)}{\img \Phi\Ef}$ is the co-limit of $\img \Phi{\wek x}$ in $\fL$.
\end{enumerate}
If $\fK$ has a \fra\ sequence of length $\kappa$ then there exists a $\kappa$-\fra\ sequence $\wek u$ in $\fK$ such that $\img \Phi{\wek u}$ is continuous in $\fL$.
\end{prop}

\begin{pf}
Note that the assumptions on the existence of a $\kappa$-\fra\ sequence imply that $\fK$ is dominated by at most $\kappa$ many arrows.
Repeat the construction from the proof of Theorem \ref{jofjaiopf}, taking care that $\Phi(u_\delta)$ be the co-limit of $\img \Phi{\wek u\rest\delta}$ at every limit stage $\delta$.
\end{pf}

Selected examples of \fra\ sequences will be described in Sections~\ref{SectExmplsFirstSeries} and \ref{SectAppsEPetsoon}.
Below we illustrate the main concepts in the case of ordered categories.

\begin{ex}
Let $\pair P \loe$ be a partially ordered (or, more generally, quasi-ordered) set, treated as a category.
Observe that the amalgamation property is equivalent to directedness.
A sequence in $P$ is nothing but a well ordered chain.
Note that a subset of $P$ is dominating iff it is cofinal.
In particular, a sequence in $P$ is \fra\ iff it is cofinal in $P$.
Finally, it is well known and easy to show by a simple transfinite induction that a directed poset with no maximal element has a well ordered cofinal chain of type $\kappa$ (a $\kappa$-\fra\ sequence in our terminology) if and only if it is $\kappa$-bounded (i.e. every subset of cardinality $<\kappa$ has an upper bound in $P$) and has any cofinal set of cardinality $\kappa$.
Such a cardinal $\kappa$ is necessarily regular. 
\end{ex}

\subsection{Cofinality}

Below we discuss the crucial property of a \fra\ sequence: cofinality in the category of sequences.

\begin{tw}[Cofinality for $\omega$-sequences]\label{saofjapiwf}
Assume $\wek u$ is a \fra\ sequence in a category $\fK$ with the amalgamation property. Then for every countable inductive sequence $\wek x$ in $\fK$ there exists a morphism of sequences $\map F{\wek x}{\wek u}$.
\end{tw}
\index{cofinality for $\omega$-sequences}

\begin{pf}
We use the extension property (property (E)) of the sequence $\wek u$, which is equivalent to the amalgamation property of $\fK$ (Proposition \ref{krakrakra}). Let $\wek x$ be an $\omega$-sequence in $\fK$.
Using (U), find an arrow $\map {f_0}{x_0}{u_{\al_0}}$. Now assume that arrows $f_0,\dots,f_{n-1}$ have been defined so that $\map{f_m}{x_m}{u_{\al_m}}$ and the diagram
$$\xymatrix{
x_\ell \ar[r]^{f_\ell} & u_{\al_\ell} \\
x_k \ar[u]^{x_k^\ell} \ar[r]^{f_k} & u_{\al_k} \ar[u]_{u_{\al_k}^{\al_\ell}}
}$$
commutes for every $k<\ell<n$ (in particular $\al_0\loe\al_1\loe\dots\loe \al_{n-1}$).
Using (E), find $\al_n\goe\al_{n-1}$ and an arrow $\map f{x_n}{u_{\al_n}}$ so that $f\cmp x^n_{n-1} = u_{\al_{n-1}}^{\al_n}\cmp f_{n-1}$ and define $f_n:=f$. Given $m<n-1$, by the induction hypothesis, we get
$$f_n\cmp x_m^n=f_n\cmp x_{n-1}^n\cmp x_m^{n-1} = u_{\al_{n-1}}^{\al_n}\cmp f_{n-1}\cmp x_m^{n-1} = u_{\al_{n-1}}^{\al_n}\cmp u_{\al_m}^{\al_{n-1}}\cmp f_m = u_{\al_{m}}^{\al_n}\cmp f_m.$$
Finally, setting $F=\sett{f_n}{\ntr}$, we obtain the required morphism $\map F{\wek x}{\wek u}$.
\end{pf}

The proof above can be easily extended to uncountable sequences, assuming continuity:

\begin{tw}\label{continussdfsdfarfqarw}
Let $\fK$ be a category with the amalgamation property and let $\wek u$ be a \fra\ sequence of regular length $\kappa$ in $\fK$. Then for every continuous sequence $\wek x\in\uzuple\kappa\fK$ there exists an arrow of sequences $\map F{\wek x}{\wek u}$.
\end{tw}

\begin{pf}
We repeat the construction from the proof of Theorem \ref{saofjapiwf}. In the case of a limit ordinal $\delta$, we let $\al_\delta$ to be the supremum of $\setof{\al_\xi}{\xi<\delta}$ and we define $f_\delta$ to be the unique arrow satisfying $f_\delta\cmp x_\al^\delta = f_\al$ for every $\al<\delta$. This is possible, because $x_\delta$ together with the co-cone of arrows $\sett{x_\xi^\delta}{\xi<\delta}$ is, by assumption, the co-limit of $\wek x \rest \delta$. Thus, the construction from the proof of Theorem \ref{saofjapiwf} can be carried out, obtaining the desired arrow $\map F{\wek x}{\wek u}$.
\end{pf}

We shall see later that an uncountable \fra\ sequence may not be cofinal for $\omega_1$-sequences. From Theorem \ref{saofjapiwf} we immediately get the following characterization of the existence of a \fra\ sequence of length $\omega_1$.

\begin{wn}
Let $\fK$ be a category with the amalgamation property. There exists a \fra\ sequence of length $\omega_1$ in $\fK$ if and only if $\fK$ is $\sig$-bounded and dominated by a family consisting of at most $\aleph_1$ arrows.
\end{wn}

\begin{pf}
The ``if" part is a special case of Theorem \ref{jofjaiopf}.
Let $\wek u$ be an $\omega_1$-\fra\ sequence in $\fK$. Then $\fK$ is directed and the family $\setof{u_\al^\beta}{\al\loe\beta<\omega_1}$ is dominating in $\fK$. Fix $\wek x\in\ciagi\fK$. Theorem \ref{saofjapiwf} says that there exists an arrow of sequences $\map{\wek f}{\wek x}{\wek u}$, so some $u_\al$ provides a bound for $\wek x$. Thus, every countable sequence is bounded in $\fK$.
\end{pf}

\subsection{The back-and-forth principle}

\index{back-and-forth principle}
Fix a category $\fK$. Let $\wek u$, $\wek v$ be \fra\ sequences in $\fK$. We shall say that $\pair{\wek u}{\wek v}$ satisfies the \emph{back-and-forth principle} if for every $\al$ below the length of $\wek u$, for every arrow $\map{f}{u_\al}{\wek v}$ there exists an isomorphism of sequences $\map{h}{\wek u}{\wek v}$ such that $h\cmp u^\infty_\al = f$, i.e. the following diagram commutes:
$$\xymatrix{
\wek u\ar[r]^h & \wek v\\
u_\al\ar[u]^{u^\infty_\al}\ar[ur]_f & 
}$$
Since there exists at least one arrow $\map{f}{u_0}{\wek v}$, this implies that $\wek u\iso \wek v$. 
We shall say that $\wek u$ satisfies the \emph{back-and-forth principle} if $\pair{\wek u}{\wek u}$ does. The following simple statement shows the importance of this property.

\begin{prop}[Homogeneity]
Assume $\fK$ is a category with the amalgamation property and $\wek u$ is a \fra\ sequence in $\fK$ which satisfies the back-and-forth principle. Then for every arrows $\map fab$, $\map ia{\wek u}$ and $\map jb{\wek u}$ there exists an isomorphism of sequences $\map h{\wek u}{\wek u}$ such that $h \cmp i = j\cmp f$, i.e. the following diagram commutes.
$$\xymatrix{
\wek u\ar[r]^h & \wek u\\
a\ar[u]^{i}\ar[r]_f & b\ar[u]_j
}$$
\end{prop}

\begin{pf}
Let $\al$ be such that $i = u_\al^\infty \cmp i_\al$ for some $\map {i_\al}a{u_\al}$.
Let $\wek v$ be the sequence obtained from $\wek u$ by cutting all objects with indices below $\al$ and adding the arrow $\map {i_\al}a{u_\al}$ at the beginning. That is,
$$\xymatrix{
a\ar[rr]^{i_\al} & & u_\al\ar[rr]^{u_\al^{\al+1}} & & u_{\al+1}\ar[rr] & & \dots
}$$
Since $\fK$ has the amalgamation property, by Proposition \ref{Prlaskaamodr}, $\wek v$ is \fra. The back-and-forth principle applied to $j\cmp f$ yields the required isomorphism $h$.
\end{pf}

It turns out that countable \fra\ sequences always satisfy the back-and-forth principle. We shall see in Section \ref{drzewa} that this not true for sequences of length $\omega_1$.

\begin{tw}[Uniqueness]\label{jednorodnost}
Assume that $\wek u$, $\wek v$ are $\omega$-\fra\ sequences in a given category $\fK$. Assume further that $k,\ell<\omega$ and $\map f{u_k}{v_\ell}$ is an arrow in $\fK$. Then there exists an isomorphism $\map F{\wek u}{\wek v}$ in $\ciagi\fK$ such that the diagram
$$\xymatrix{
\wek u\ar[r]^F & \wek v\\
u_k\ar[u]^{u_k^\infty}\ar[r]^f & v_\ell\ar[u]_{v_\ell^\infty}
}$$
commutes. In particular ${\wek u}\iso {\wek v}$.
\end{tw}

Notice that in the statement above we do not need to assume that the given category has the amalgamation property.

\begin{pf}
We construct inductively arrows $\map{f_n}{u_{k_n}}{v_{\ell_n}}$, $\map{g_n}{v_{\ell_n}}{u_{k_{n+1}}}$, where $k_0\loe\ell_0 < k_1\loe \ell_1< \dots$ and for each $\ntr$ the diagram
$$\xymatrix{
u_{k_n} \ar[dr]^{f_n}\ar[rr]^{u_{k_n}^{k_{n+1}}} & & u_{k_{n+1}} \ar[rd]^{f_{n+1}}\ar[rr]^{u_{k_{n+1}}^{k_{n+2}}} & & u_{k_{n+2}}\ar@{-->}[r] &\\
& v_{\ell_n} \ar[ur]_{g_n}\ar[rr]_{v_{\ell_n}^{\ell_{n+1}}} & & v_{\ell_{n+1}} \ar[ur]_{g_{n+1}}\ar@{-->}[rr] & &
}$$
commutes.

We start with $f_0:=f$, $k_0:=k$, $\ell_0:=\ell$, possibly replacing $f$ by some arrow of the form $j_\ell^m\cmp f$ to ensure that $k_0\loe\ell_0$. Using property (A) of the sequence $\wek u$, find $k_1>k_0$ and $\map{g_1}{v_0}{u_{k_1}}$ such that $g_1\cmp f = u_0^{k_1}$.
Assume that $f_m,g_m$ have already been constructed for $m\loe n$. Using the amalgamation of $\wek v$, find $\ell_{n+1}\goe k_{n+1}$ and an arrow $\map{f_{n+1}}{u_{k_{n+1}}}{v_{\ell_{n+1}}}$ such that $f_{n+1}\cmp g_n = v_{\ell_n}^{\ell_{n+1}}$. Now, using the amalgamation of $\wek u$, we find $k_{n+2}>\ell_{n+1}$ and an arrow $\map{g_{n+1}}{v_{\ell_{n+1}}}{u_{k_{n+2}}}$ such that $g_{n+1}\cmp f_{n+1} = u_{k_{n+1}}^{k_{n+2}}$.
By the induction hypothesis, $g_n\cmp f_n = u_{k_n}^{k_{n+1}}$, therefore the above diagram commutes. This finishes the construction.

Finally, set $F = \sett{f_n}{\ntr}$ and $G = \sett{g_n}{\ntr}$. Then $\map F{\wek u}{\wek v}$, $\map G{\wek v}{\wek u}$ are morphisms of sequences and by a simple induction we show that
\begin{equation}
g_n\cmp v_{\ell_m}^{\ell_n}\cmp f_m = u_{k_m}^{k_{n+1}}\oraz f_n\cmp u_{k_{m+1}}^{k_n}\cmp g_m = v_{\ell_m}^{\ell_n}
\tag{*}\end{equation}
holds for every $m<n<\nat$. This shows that $F\cmp G=\id{\wek v}$ and $G\cmp F=\id{\wek u}$, therefore $F$ is an isomorphism. The equality $v_0^\infty\cmp f = F\cmp u_0^\infty$ means that $v_{0}^{\ell_n}\cmp f = f_n\cmp u_0^{k_n}$ should hold for every $\ntr$. Fix $n>0$. Applying (*) twice (with $m=0$ and $m=n-1$ respectively), we get
$$f_n\cmp u^{k_n}_0 = f_n\cmp g_{n-1}\cmp v_0^{\ell_{n-1}}\cmp f = v^{\ell_n}_{\ell_{n-1}}\cmp v^{\ell_{n-1}}_0\cmp f = v^{\ell_n}_0\cmp f.$$
Thus $v_0^\infty\cmp f = F\cmp u_0^\infty$.

Finally, notice that, by property (U) of the sequence $\wek v$, for some $\ell<\omega$ there exists an arrow $\map f{u_0}{v_\ell}$, so applying the first part we see that ${\wek u}\iso {\wek v}$.
\end{pf}

\subsection{\fra\ sequences and functors}

We now discuss a possible extension of the notion of a \fra\ sequence, involving a functor.
The motivation is two-fold: First, we may want to have a sequence that is not \fra\ in a given category $\fK$, yet resembles the \fra\ property when moving to a different, richer category $\fL$.
Second, given a \fra\ sequence, we may be unable to show its cofinality, however we may succeed after transferring the sequence to a different category. Of course, this makes sense only if the functor we are using preserves relevant information.

\index{\fra\ sequence!-- over a functor}
Fix a functor $\map \Phi \fK \fL$.
We say that a sequence $\map {\wek u} \kappa \fL$ is \emph{\fra\ over} $\Phi$ if the following conditions are satisfied.
\begin{enumerate}
	\item[($U_\Phi$)] For every $x\in \fL$ there is $\al<\kappa$ such that $\fL(x,\Phi(u_\al))\nnempty$.
	\item[($A_\Phi$)] Given $\al<\kappa$ and $\map f{u_\al}y$ in $\fK$, there are $\beta\goe\al$ and an arrow $\map gy{\Phi(u_\beta)}$ such that $\Phi(u_\al^\beta) = g \cmp \Phi(f)$.
\end{enumerate}
In case of the identity functor (where $\fK=\fL$) this is just the usual notion of a \fra\ sequence.

\begin{prop}
Let $\map\Phi\fK\fL$ be a covariant functor and let $\wek u$ be a \fra\ sequence in $\fK$. If $\img\Phi\fK$ is cofinal in $\fL$ then $\wek u$ is \fra\ over $\Phi$.
\end{prop}

Aiming at the problem of existence and cofinality, we now introduce the following natural concept.
\index{amalgamation property!-- mixed}
\index{mixed amalgamation property}
Let $\map \Phi\fK\fL$ be a covariant functor. 
We say that $\Phi$ has the {\em mixed amalgamation property}
if for every objects $a,b,c$ in $\fK$, for every arrow $\map fab$ in $\fK$ and for every arrow $\map g{\Phi(a)}{\Phi(c)}$ in $\fL$, there exist $d\in \fK$ and arrows $\map{f'} c d$, $\map{g'}{\Phi(b)}{\Phi(d)}$ in $\fK$ and $\fL$ respectively, such that the following diagram
$$\xymatrix{
\Phi(b)\ar[r]^{g'} & \Phi(d)\\
\Phi(a)\ar[u]^{\Phi(f)}\ar[r]^{g} & \Phi(c)\ar[u]_{\Phi(f')}
}$$
is commutative.
The usual amalgamation can be rephrased as the amalgamation property of the identity functor.
Given a category $\fL$ and its subcategory $\fK$, we say that $\pair{\fK}{\fL}$ has the \emph{mixed amalgamation property} if the inclusion functor $\fK \subs \fL$ has this property.

In order to formulate the existence result, we need to adapt the notion of a dominating family.
Let $\map \Phi \fK \fL$ be a functor.
A family $\Ef \subs \fK$ will be called \emph{$\Phi$-dominating} if
\begin{enumerate}
	\item[(D$1_\Phi$)] For every $\fK$-object $x$, there is $a \in \Dom(\Ef)$ such that $\fL(\Phi(x), \Phi(a)) \nnempty$.
	\item[(D$2_\Phi$)] Given $a \in \Dom(\Ef)$ and a $\fK$-arrow $\map f a y$, there exist $\map h a b$ in $\Ef$ and an $\fL$-arrow $\map g {\Phi(y)}{\Phi(b)}$ such that $\Phi(h) = g \cmp \Phi(f)$.
\end{enumerate}

\begin{tw}\label{ThmtwIstnienief}

Let $\kappa$ be an infinite regular cardinal and let $\map \Phi \fK \fL$ be a functor satisfying the following conditions:
\begin{enumerate}
	\item[$(F_1)$] For every $\fK$-objects $x,y$ there exists a $\fK$-object $z$ such that $\fK(x,z)\nnempty \ne \fL(\Phi(y), \Phi(z))$.
	\item[$(F_2)$] $\Phi$ has the mixed amalgamation property.
	\item[$(F_3)$] $\fK$ is $\kappa$-bounded and has a $\Phi$-dominating family consisting of $\loe \kappa$ arrows.
\end{enumerate}
Then there exists a $\kappa$-\fra\ sequence over $\Phi$.
If, additionally, $\Phi$ satisfies condition ($\mathcal C$) of Proposition~\ref{strangewef} then there exists a $\kappa$-\fra\ sequence over $\Phi$ such that $\img \Phi {\wek u}$ is continuous.
\end{tw}

\begin{pf}
This is actually a repetition of the arguments from the proof of Theorem~\ref{jofjaiopf} combined with Proposition~\ref{strangewef}.
The only difference is that we need to use mixed amalgamations.
More precisely, having defined $u_\xi$ and $u_\xi^\eta$ for $\xi < \eta \loe \beta$, as in the proof of Theorem~\ref{jofjaiopf}, given a $\fK$-arrow $\map f{u_\al} y$
that belongs to a fixed $\Phi$-dominating family $\Ef$, we use the mixed amalgamation property followed by condition (D$2_\Phi$) in order to get a $\fK$-object $u_{\beta+1}$ and a $\fK$-arrow $\map {u_\beta^{\beta+1}}{u_\beta}{u_{\beta+1}}$ such that $u_\beta^{\beta+1} \in \Ef$ and $\Phi(u_\beta^{\beta+1}) = g \cmp \Phi(f)$ for some $g \in \fL$.

If $\Phi$ satisfies ($\mathcal C$) then we use this condition at limit steps in order to make the sequence $\Phi$-continuous.
\end{pf}

The mixed amalgamation property is also the main tool for showing cofinality of a \fra\ sequence ``transferred" to a different, possibly much bigger, category.
As we shall see in Subsection~\ref{SubXtreesfinerrw}, there exist quite natural pairs of categories $\fK \subs \fL$ such that $\fK$ has neither the amalgamation property nor a \fra\ sequence, yet there exists a \fra\ sequence over the inclusion $\fK \subs \fL$.

\begin{tw}\label{TmixedamsCof}
Let $\kappa$ be an infinite regular cardinal and let $\map\Phi\fK\fL$ be a functor with the mixed amalgamation property.
Assume further that $\map{\wek u}\kappa \fK$ is a \fra\ sequence over $\Phi$.
Then for every sequence $\wek x \in\uzuple{\kappa}\fK$ such that $\img\Phi{\wek x}$ is continuous in $\fL$, there exists an arrow $\map F{\img\Phi{\wek x}}{\img\Phi{\wek u}}$ in $\uzuple{\kappa}\fL$.
\end{tw}

\begin{pf}
Fix a sequence $\wek x$ in $\fK$, of length $\lam = \cf \lam \loe \kappa$, and assume that $\img \Phi{\wek x}$ is continuous in $\fL$.
We define inductively arrows $\map{f_\al}{\img\Phi{x_\al}}{\img\Phi{u_{\theta(\al)}}}$ so that $\wek f=\sett{f_\al}{\al < \lam}$ will be an arrow of sequences in $\fL$.
We start with an arbitrary arrow $\map{f_0}{\img\Phi{x_0}}{\img\Phi{u_{\theta(0)}}}$, using the fact that $\wek u$ is \fra\ over $\Phi$.
Fix $\al>0$ and suppose $f_\xi$ have already been defined for all $\xi<\al$ so that
$$f_\eta \cmp \img\Phi{x_\xi^\eta} = \img\Phi{u_{\theta(\xi)}^{\theta(\eta)}} \cmp f_\xi$$
holds for every $\xi<\eta<\al$.

If $\al$ is a limit ordinal, define $\theta(\al)=\sup_{\xi<\al}\theta(\xi)$ and, using the continuity of $\img\Phi{\wek x}$, let $\map{f_\al}{\img\Phi{x_\al}}{\img\Phi{u_{\theta(\al)}}}$ be the unique arrow satisfying $f_\al \cmp \img\Phi{x^\al_\xi} = \img\Phi{u_{\theta(\xi)}^{\theta(\al)}} \cmp f_\xi$ for every $\xi<\al$.

Now assume that $\al=\beta+1$.
Using the mixed amalgamation property of $\Phi$, we find an arrow $\map g{u_{\theta(\beta)}}w$ in $\fK$ and an arrow $\map h{\img\Phi{x_\al}}{\img\Phi w}$ in $\fL$ such that $h \cmp \img\Phi{j_\beta^\al} = \img\Phi{g} \cmp f_\beta$.
Since $\wek u$ is \fra\ over $\Phi$, there are $\gamma\goe\theta(\beta)$ and $\map k{\img\Phi w}{\img\Phi{u_\gamma}}$ such that 
$\img\Phi{u_{\theta(\beta)}^\gamma} = k \cmp \img\Phi g$.
Finally, set $\theta(\al):=\gamma$ and $f_\al:=k \cmp h$.

This finishes the construction and completes the proof.
\end{pf}

\section{\frajon\ categories}\label{SectFraJonxqlims}

In this section we describe a way of dealing with discontinuous \fra\ sequences. Namely, it may happen that ``short" sequences do not have co-limits, however there is a (sometimes quite natural) operator of quasi-limit (see the definition below).
If this operator is good enough, the back-and-forth principle holds.
If moreover the quasi-limit commutes with certain amalgamations, it is possible to prove cofinality of a quasi-continuous \fra\ sequence.
The results of this section are motivated by applications, described later.

\subsection{Quasi-limiting operators}

What is really needed for the back-and-forth principle?
It is rather clear that the existence of co-limits for ``short" sequences is sufficient, because then one can restrict attention to continuous \fra\ sequences.
Looking for specific applications, we are forced to relax this condition.
So, let us analyze co-limits for a moment.

Given a category $\fK$, let $\cocones\fK$ be the category of all diagrams of the form $\sett{f_\xi}{\xi<\lam}$ where $\lam$ is an ordinal and all $f_\xi$s have the same co-domain.
In other words, $\cocones \fK$ is the category of all functors from posets of the form $\lam \cup \sn\infty$ where $\lam$ has the discrete order (different elements are incomparable) and $\infty > \xi$ for every $\xi \in \lam$.
The arrows are natural transformations.
Assume all sequences in $\fK$ of length $<\kappa$ have co-limits in $\fK$.
Then we have an operator $\map \lim {\uzup\kappa \fK}{\cocones\fK}$ which assigns to each sequence $\map{\wek x} \lam \fK$ its co-limiting co-cone $\lim\wek x$.
This is indeed a functor, because of the universality of the co-limit.
In the back-and-forth argument, we actually deal with isomorphisms at each limit step, so it turns out that only some properties of $\lim$ are needed.

\index{quasi-limiting operator}
Denote by $\uzupiso \kappa \fK$ the category $\uzup \kappa \fK$ restricted to isomorphisms of sequences only.
A \emph{quasi-limiting operator} for $\uzup \kappa \fK$ is a functor $\map \qlim {\uzupiso \kappa \fK}{\cocones \fK}$ assigning to each sequence $\wek x$ a co-cone $\qlim {\wek x}$ for $\wek x$.
More precisely, given an isomorphism of sequences $\map {\wek h}{\wek x}{\wek y}$ there exists an isomorphism $\map {h_\infty}{\qlim{\wek x}}{\qlim{\wek y}}$ commuting with $\wek h$ in the usual sense.
\index{quasi-limiting operator!-- partial}
We say that $\qlim$ is a \emph{partial quasi-limiting operator} if it is defined on a subcategory of $\uzupiso \kappa \fK$ with the property that if $\qlim \wek x$ is defined then so is $\qlim \wek y$ whenever $\wek y \iso \wek x$.
An example is the co-limiting operator, defined only for sequences that have co-limits.

Typical examples of quasi-limiting operators come from co-limits, ``computed" in a bigger category containing $\fK$.

\begin{prop}\label{Pkwasilimitz}
Let $\map F \fK \fL$ be a faithful functor, $\kappa = \cf \kappa > \aleph_0$, and assume the following conditions are satisfied:
\begin{enumerate}
	\item[$(\circledast_1)$] For every sequence $\wek x$ in $\fK$ of length $\lam < \kappa$ there exists a 
co-cone $\qlim \wek x := \sett{x_\xi^\infty}{\xi < \lam} \subs \fK$ such that $\sett{F(x_\xi^\infty)}{\xi < \lam}$ is the co-limiting co-cone of $F(\wek x)$ in $\fL$.
	\item[$(\circledast_2)$] For every $\fK$-objects $a,b$, for every $\fL$-isomorphism $\map h {F(a)}{F(b)}$ there exists a $\fK$-isomorphism $\map g a b$ such that $h = F(g)$.
\end{enumerate}
Then $\qlim$ is a quasi-limiting operator on $\uzup \kappa \fK$ satisfying
$F(\qlim \wek x) = \lim F(\wek x)$
for every $\wek x \in \uzup \kappa \fK$.
\end{prop}

\begin{pf}
Fix an isomorphism of sequences $\map {\wek h}{\wek x}{\wek y}$ in $\fK$.
Let the co-cones $\qlim \wek x$, $\qlim \wek y$ have co-domains $\ovr x$ and $\ovr y$, respectively.
There is an isomorphism $\map {k}{\ovr x}{\ovr y}$ in $\fL$, commuting with $F(\wek h)$.
That is, for every $\al, \beta$ such that $\map {h_\al}{x_\al}{x_\beta}$, we have a commuting diagram
$$\xymatrix{
F(x_\al) \ar[d]_{F(h_\al)} \ar[rrr]^{F(x_\al^\infty)} & & & F(\ovr x) \ar[d]^k \\
F(y_\beta) \ar[rrr]_{F(y_\beta^\infty)} & & & F(\ovr y)
}$$
commutes.
Let $\map g {\ovr x}{\ovr y}$ be an isomorphism such that $F(g) = k$.
Then the diagram
$$\xymatrix{
x_\al \ar[d]_{h_\al} \ar[rrr]^{x_\al^\infty} & & & \ovr x \ar[d]^g \\
y_\beta \ar[rrr]_{y_\beta^\infty} & & & \ovr y
}$$
is commutative, because $F$ is faithful.
This shows that $\qlim$ is indeed a functor from $\uzupiso \kappa \fK$ into the co-cones of $\fK$.
The faithfulness of $F$ guarantees that $\qlim$ is well defined, i.e. the co-cone in condition $(\circledast_1)$ is unique, up to isomorphism.
\end{pf}

The next example shows that quasi-limiting operators may have nothing to do with co-limits.

\begin{ex}
Consider $\omega_1$ as a poset category.
For each limit ordinal $\delta < \omega_1$ add its copy $\delta^+$, incomparable with $\delta$ and such that 
$\xi < \delta^+ < \eta$ whenever $\xi < \delta < \eta$.
Let $\fK = \omega_1 \cup \setof{\delta^+}{\delta \in \lim(\omega_1)}$.
It is clear that no strictly increasing sequence has a co-limit in $\fK$.
On the other hand, there exists a quasi-limiting operator $\map \qlim {\uzup {\omega_1} \fK}{\coconesiso \fK}$ defined by
$$\qlim{\wek x} = \sup_{\ntr}P(x_n),$$
where $P$ is the canonical projection of $\fK$ onto $\omega_1$.

This example also shows that a quasi-limiting operator may not be unique, since the formula $\qlim^+{\wek x} = (\qlim{\wek x})^+$ defines another one.
\end{ex}

Note that in the example above, as in every poset category, the back-and-forth principle holds for every two \fra\ sequences.
A \fra\ sequence in a poset is a cofinal chain whose length is a regular cardinal.
If such a chain exists, then it is equivalent to any other cofinal chain.
Summarizing: an infinite poset $\poset$ has a \fra\ sequence if and only if it is $\kappa$-directed (i.e. every subset of cardinality $< \kappa$ has an upper bound) for some cardinal $\kappa$, and dominated by $\kappa$ many elements, i.e. $\poset$ has a cofinal subset of cardinality $\kappa$.
Of course, such a cardinal $\kappa$ must be regular.

\subsection{The back-and-forth principle revisited}

\index{sequence!-- $\qlim$-continuous}
\index{$\qlim$-continuous!-- sequence}
Let $\fK$ be a category and let $\qlim$ be a partial quasi-limiting operator on $\uzup \kappa \fK$.
Let $\rho < \kappa$ be an ordinal.
We say that a sequence $\map {\wek x} \rho \fK$ is $\qlim$-continuous if for every limit ordinal $\delta < \rho$ it holds that
$$x_\delta = \qlim (\wek x \rest \delta)$$
with the co-cone $\sett{x_\xi^\delta}{\xi < \delta}$.
This means in particular that $\qlim$ is defined on all subsequences of $\wek x$.

The following fact is rather trivial, knowing the proof of Theorem~\ref{jofjaiopf}.

\begin{prop}\label{Pstrangitmil}
Let $\fK$ be a directed category with the amalgamation property. Let $\kappa\goe\aleph_0$ be a regular cardinal and let $\qlim$ be a quasi-limiting operator on $\uzup \kappa \fK$.
If $\fK$ has a \fra\ sequence of length $\kappa$ then $\fK$ also has an $\qlim$-continuous one.
\end{prop}

\begin{tw}\label{Tfjogosdg}
Let $\kappa > \aleph_0$ be a regular cardinal and let $\fK$ be a category with a partial quasi-limiting operator $\qlim$ on $\uzup \kappa\fK$.
Then every two $\qlim$-continuous $\kappa$-\fra\ sequences in $\fK$ satisfy the back-and-forth principle.

\end{tw}

\begin{pf}
Let $\wek u$, $\wek v$ be $\qlim$-continuous \fra\ sequences of length $\kappa$ in $\fK$ and fix an arrow $\map f{u_{0}}{\wek v}$.
We construct inductively arrows $\map{f_\al}{u_{\psi(\al)}}{v_{\phi(\al)}}$ and $\map{g_\al}{v_{\phi(\al)}}{u_{\psi(\al+1)}}$, where $\map \phi\kappa\kappa$, $\map \psi\kappa\kappa$ are increasing functions, so that the following conditions are satisfied.
\begin{enumerate}
	\item[(i)] $\al\loe \psi(\al)\loe \phi(\al) < \psi(\al+1)$.
	\item[(ii)] $g_\al\cmp f_\al = u_{\psi(\al)}^{\psi(\al+1)}$ and $f_{\al+1}\cmp g_\al = v_{\phi(\al)}^{\phi(\al+1)}$.
	\item[(iii)] $\xi<\eta\implies v_{\phi(\xi)}^{\phi(\eta)}\cmp f_\xi = f_\eta\cmp u_{\psi(\xi)}^{\psi(\eta)}$.
\end{enumerate}
We start with $\psi(0)=0$, $\phi(0)=\al$ and $f_0=f$, where $0<\al<\kappa$ is such that $f$ is equivalent to $\map {f_0}{u_0}{v_\al}$.

Fix $\beta>0$ and assume that $\sett{f_\al}{\al<\beta}$ and $\sett{g_\al}{\al<\beta}$ have already been defined, so in particular $\phi(\al)$, $\phi(\al)$, $\psi(\al)$ and $\psi(\al+1)$ have been defined for every $\al<\beta$.

Suppose first that $\beta$ is a successor ordinal, say $\beta=\al+1$. Using the fact that $\wek v$ is \fra, find $\phi(\beta)\goe\psi(\al)$ and $\map{f_\beta}{u_{\psi(\beta)}}{v_{\phi(\beta)}}$ so that the second part of (ii) holds.
Using the fact that $\wek u$ is \fra, find $\psi(\beta+1) > \phi(\beta)$ and $\map{g_\beta}{v_{\phi(\beta)}}{u_{\psi(\beta+1)}}$ so that the first part of (ii) holds. Clearly, (i) is satisfied and (iii) follows easily from the inductive hypothesis and (ii).

Assume now that $\beta$ is a limit ordinal. Let $\rho := \sup_{\al<\beta}\phi(\al) = \sup_{\al<\beta}\psi(\al)$ and define $\phi(\beta):=\psi(\beta):=\rho$, $\psi(\beta+1):=\rho+1$. By the induction hypothesis $\rho\goe\beta$, therefore (i) holds.
Observe that $\rho$ is a limit ordinal.

First, observe that $F = \sett{f_\xi}{\xi<\rho}$ and $G = \sett{g_\xi}{\xi<\rho}$ witness an isomorphism between $\wek u\rest \rho$ and $\wek v \rest \rho$, which extends to an isomorphism of their $\qlim$-limits $\map h{u_\rho}{v_\rho}$.
More specifically,
\begin{equation}
h \cmp u^\rho_{\psi(\xi)} = v^\rho_{\phi(\xi)} \cmp f_\xi
\oraz
h^{-1} \cmp v^\rho_{\phi(\xi)} = u^\rho_{\psi(\xi+1)} \cmp g_\xi
\tag{\textxswup}\label{eqsdngoisjegf}
\end{equation}
hold for every $\xi < \rho$.

Define $f_\beta := h$ and $g_\beta := u^{\rho+1}_\rho \cmp h^{-1}$.
Then, $g_\beta \cmp f_\beta = u^{\rho+1}_\rho$, i.e. the first part of (ii) holds for $\beta$ in place of $\al$ (the second part is taken care of in the next, successor, inductive step).
The first equation in (\ref{eqsdngoisjegf}) shows (iii).
This finishes the inductive construction.

We claim that $\wek f = \sett{f_\al}{\al<\kappa}$ is an isomorphism from $\wek u$ to $\wek v$ which witnesses the back-and-forth principle. For this aim, we first observe that $\wek f$ is an arrow of sequences, because of (iii).
Define $\wek g = \sett{g_\al}{\al<\kappa}$.
We claim that $g$ is an arrow of sequences, inverse to $\wek f$.

Using (ii) and (iii), we see that
\begin{align*}
g_\eta \cmp v_{\phi(\xi)}^{\phi(\eta)} &= g_\eta \cmp v_{\phi(\xi+1)}^{\phi(\eta)} \cmp v_{\phi(\xi)}^{\phi(\xi+1)} = g_\eta \cmp v_{\phi(\xi+1)}^{\phi(\eta)} \cmp f_{\xi+1} \cmp g_\xi\\
&= g_\eta \cmp f_\eta \cmp u_{\psi(\xi+1)}^{\psi(\eta)} \cmp g_\xi = u_{\psi(\eta)}^{\psi(\eta+1)} \cmp u_{\psi(\xi+1)}^{\psi(\eta)} \cmp g_\xi\\
&= u_{\psi(\xi+1)}^{\psi(\eta+1)} \cmp g_\xi
\end{align*}
for every $\xi<\eta$, which shows that $\wek g=\sett{g_\al}{\al<\kappa}$ is indeed an arrow of sequences.
Furthermore, again by (ii), we deduce that $\wek g$ is the inverse of $\wek f$, which shows that $\wek f$ is an isomorphism.

Finally, $\wek f \cmp u_0^\infty = f$, by the first stage of the construction.
\end{pf}

As a corollary, we obtain the result proved in \cite{DrGoe92}, stated there in a different language.

\begin{wn}\label{oghowehe}
Let $\fK$ be a category. Every two continuous \fra\ sequences in $\fK$ of the same regular length satisfy the back-and-forth principle.
\end{wn}

\begin{tw}
Let $\kappa>\aleph_0$ be a regular cardinal and let $\fK$ be a full, cofinal and relatively $\kappa$-complete subcategory of a category $\fL$ with the amalgamation property. Then
\begin{enumerate}
	\item[(a)] There exists at most one, up to isomorphism, $\kappa$-\fra\ sequence in $\fK$.
	\item[(b)] A possible $\kappa$-\fra\ sequence in $\fK$ is cofinal for $\uzuple\kappa\fK$ and satisfies the back-and-forth principle.
\end{enumerate}
\end{tw}

\begin{pf}
Let $\wek u$ be a $\kappa$-\fra\ sequence in $\fK$. 
Since $\fK$ is full and cofinal in $\fL$, $\wek u$ is \fra\ in $\fL$.
Let us make this sequence continuous in $\fL$, by adding co-limits at each limit ordinal below $\kappa$. Denote the extended sequence by $\tilde u$. Since $\fL$ has the amalgamation property, the sequence $\tilde u$ is \fra\ in $\fL$ (by Proposition \ref{laskanebeska}(b)). Given another $\kappa$-\fra\ sequence $\wek v$, let us extend it in the same way, obtaining a \fra\ sequence $\tilde v$ in $\fL$.
By Corollary \ref{oghowehe}, the two extended sequences $\tilde u$, $\tilde v$ satisfy the back-and-forth principle. Since $\fK$ is a full subcategory of $\fL$, we conclude that $\wek u$, $\wek v$ satisfy the back-and-forth principle in $\fK$.
This shows (a) and the first part of (b). The second part of (b) follows from Theorem \ref{continussdfsdfarfqarw} and the fact that each $\kappa$-sequence in $\fK$ can be completed to a continuous $\kappa$-sequence in $\fL$.
\end{pf}

\begin{wn}
Let $\kappa$ be a regular cardinal and let $\fK$ be a directed category with the amalgamation property and with a dominating family consisting of at most $\kappa$ arrows. Assume further that either $\fK$ is $\kappa$-complete or $\fK$ is $\kappa$-bounded and $\uzup\kappa\fK$ has the amalgamation property.
Then
\begin{enumerate}
	\item[(a)] There exists a unique, up to isomorphism, $\kappa$-\fra\ sequence $\wek u$ in $\fK$.
	\item[(b)] The sequence $\wek u$ is cofinal for $\uzuple\kappa\fK$ and satisfies the back-and-forth principle.
\end{enumerate}
\end{wn}

\subsection{Amalgamation structures}

We now describe a general way of obtaining cofinality, using a quasi-limiting operator and amalgamations.
As we shall see later, it may happen that sequences of amalgamations do not converge and any inductive procedure of mapping a given sequence into a \fra\ sequence fails at the first limit step.
It is possible, however, that some special amalgamations survive limit steps.
For this aim, we consider the following concepts.

\index{amalgamation structure}
\index{category!-- amalgamation structure}
Fix a category $\fK$ with the amalgamation property.
An \emph{amalgamation structure} on $\fK$ is a class $\Ama$ of commuting squares of $\fK$ (formally: functors from the 4-element Boolean algebra into $\fK$) satisfying the following conditions:
\begin{enumerate}
	\item[(A1)] Every square with two parallel identities is in $\Ama$.
	\item[(A2)] Given squares
$$\xymatrix{
b \ar[r]^{g'} & d \\
c \ar[u]^g \ar[r]_f & a \ar[u]_{f'}
}
\qquad
\xymatrix{
d \ar[r]^{f''} & v \\
a \ar[u]^{f'} \ar[r]_h & u \ar[u]_{h'}
}$$
in $\Ama$, the composition square
$$\xymatrix{
b \ar[r]^{f'' \cmp g'} & v \\
c \ar[u]^g \ar[r]_{h \cmp f} & u \ar[u]_{h'}
}$$
is again in $\Ama$.
	\item[(A3)] Given $\fK$-arrows $\map fca$, $\map gcb$, there exist $\fK$-arrows $\map {f'}av$ and $\map {g'}bv$ such that the square
$$\xymatrix{
b \ar[r]^{g'} & v \\
c \ar[r]_f \ar[u]^g & a \ar[u]_{f'}
}$$
is in $\Ama$.
\end{enumerate}

Clearly, the class $\Ama$ of all commuting squares in $\fK$ is an amalgamation structure as long as $\fK$ has the amalgamation property.
In Section~\ref{reterpeteairs} we shall consider an important type of amalgamation structures in categories of embedding-projection pairs.

\index{amalgamation structure!-- $\qlim$-continuous}
Now assume $\qlim$ is a quasi-limiting operator on $\uzup \kappa\fK$.
	We say that an amalgamation structure $\Ama$ is \emph{$\qlim$-continuous} if for every limit ordinal $\rho < \kappa$, for every $\qlim$-continuous sequences $\map {\wek x}{\rho+1}\fK$ and $\map {\wek y}{\rho+1}\fK$, for every sequence of $\fK$-arrows $\sett{\map {f_\xi}{x_\xi}{y_\xi}}{\xi < \rho}$ such that for each $\xi < \eta < \rho$ the square
$$\xymatrix{
y_\xi \ar[rr]^{y_\xi^\eta} & & y_\eta \\
x_\xi \ar[u]^{f_\xi} \ar[rr]_{x_\xi^\eta} & & x_\eta \ar[u]_{f_\eta}
}$$
is in $\Ama$, there exists $\map {f_\rho}{x_\rho}{y_\rho}$ such that for each $\xi < \rho$ the square
$$\xymatrix{
y_\xi \ar[rr]^{y_\xi^\rho} & & y_\rho \\
x_\xi \ar[u]^{f_\xi} \ar[rr]_{x_\xi^\rho} & & x_\rho \ar[u]_{f_\rho}
}$$
is in $\Ama$.

The existence of an $\qlim$-continuous amalgamation structure is sufficient for proving cofinality, as we shall see in a moment.

We now introduce the main notion of this section.
For defining one of the conditions, we need to consider posets $\truj \lam$ of at most two-element subsets of an ordinal $\lam$.
This poset can be regarded as the set of all ordered pairs $\pair \xi \eta$ in $\lam \times \lam$ with $\xi \loe \eta$.
The ordering is coordinate-wise.
Given a functor $\map X{\truj \lam}\fK$, we shall denote by $X_\Delta$ its \emph{diagonal sequence} $X \cmp \delta$, where $\map \delta \lam \fK$ is defined by $\delta(\xi) = \dn \xi \xi$.
Furthermore, we shall denote by $X^\al$ the $\al$th \emph{vertical sequence} $X \cmp v_\al$, where $\map {v_\al}\lam {\truj \lam }$ is defined by $v_\al(\xi) = \dn{\al}{\al + \xi}$.

\index{\frajon\ category}
\index{category!-- \frajon}
A \emph{\frajon\ category} is a triple of the form $\triple \fK \qlim \Ama$, where $\qlim$ is a quasi-limiting operator on $\uzup \kappa \fK$ for some infinite regular cardinal $\kappa$, $\Ama$ is an $\qlim$-continuous amalgamation structure on $\fK$ and the following condition is satisfied.
\begin{enumerate}
	\item[($\Delta$)] Given an infinite cardinal $\lam < \kappa$, given a functor $\map X{\truj \lam}\fK$ such that for every $\al_0 < \al_1 \loe \beta_0 < \beta_1$ the square
$$\xymatrix{
X_{\{\al_0,\beta_1\}} \ar[rr]^{X_{\{\al_0,\beta_1\}}^{\{\al_1,\beta_1\}}} & & X_{\{\al_1,\beta_1\}} \\
& & \\
X_{\{\al_0,\beta_0\}} \ar[uu]^{X_{\{\al_0,\beta_0\}}^{\{\al_0,\beta_1\}}} \ar[rr]_{X_{\{\al_0,\beta_0\}}^{\{\al_1,\beta_0\}}} & & X_{\{\al_1,\beta_0\}} \ar[uu]_{X_{\{\al_1,\beta_0\}}^{\{\al_1,\beta_1\}}}
}$$
is in $\Ama$,
it holds that
$$\qlim(X_\Delta) = \qlim(\qlim X),$$
where $\qlim X$ is the sequence obtained from the $\qlim$-continuity of $\Ama$:
$$\xymatrix{
\qlim X^0 \ar[r] & \qlim X^1 \ar[r] & \dots \ar[r] & \qlim X^\al \ar[r] & \dots
}$$
and the corresponding arrows and co-cones are natural.
\end{enumerate}
We shall call condition ($\Delta$) the \emph{diagonalization property}.
\index{diagonalization property}

Note that the diagonalization property is obvious when one deals with co-limits, even after moving by a faithful functor to a different category.
Namely:

\begin{prop}\label{Prhgiues}
Let $\map F \fK \fL$ be a faithful functor satisfying conditions $(\circledast_1)$, $(\circledast_2)$ of Proposition~\ref{Pkwasilimitz} for some $\kappa = \cf \kappa > \aleph_0$ and let $\qlim$ be the quasi-limiting operator such that $F(\qlim \wek x) = \lim F(\wek x)$ for every $\wek x \in \uzup \kappa \fK$.
Then for every $\qlim$-continuous amalgamation structure $\Aaa$, $\triple \fK \qlim \Aaa$ has the diagonalization property, i.e., $\triple \fK \qlim \Aaa$ is a \frajon\ category.
\end{prop}

\begin{pf}
Assume $\Aaa$ is an $\qlim$-continuous amalgamation structure on $\fK$.
Fix a suitable functor $\map X {\dpower \lam 2}\fK$ in which all squares are in $\Aaa$.
Note that the usual co-limit has the diagonalization property, therefore we have $\lim F(X_\Delta) = \lim(\lim F(X))$, with the notation above.
Since $F$ is faithful and satisfies $(\circledast_2)$, we conclude that $\qlim(X_\Delta) = \qlim(\qlim X)$.
\end{pf}

We now give a simple example showing that ($\Delta$) does not follow from the other assumptions.

\begin{ex}
Let $P_0 = \setof{\pair m n \in \omega \times (\omega+1)}{m \loe n}$ treated as a poset category with the coordinate-wise ordering.
Let $P = P_0 \cup \dn a b$ where $a, b$ are incomparable and $\pair m n < a$ precisely for $m \loe n < \omega$, while $\pair m n < b$ for every $m \loe n \loe \omega$.
We define a quasi-limiting operator $\qlim$ on $P$ as follows.
Given a strictly increasing sequence $\wek x$ we define $\qlim \wek x = \sup_{P_0}\wek x$ provided that $\wek x \subs k \times (\omega + 1)$ for some $k$.
Otherwise, we have two possibilities: either $\wek x \subs \omega \times \omega$ and then we set $\qlim \wek x = a$, or else $\wek x = \ciag x$, where $x_n \in \omega \times \sn \omega$ for all but finitely many $\ntr$.
In this case we set $\qlim \wek x = b$.
It is clear that this  defines a quasi-limiting operator and $P_0 \subs P$ witnesses the failure of ($\Delta$).
\end{ex}

\subsection{Cofinality revisited}

\begin{tw}\label{Tkofnltyrew}
Let $\triple \fK \qlim \Ama$ be a $\kappa$-\frajon\ category, where $\kappa$ is an infinite regular cardinal.
Assume $\wek u$ is an $\qlim$-continuous $\kappa$-\fra\ sequence in $\fK$.
Then for every $\qlim$-continuous sequence $\wek x \in \uzuple \kappa \fK$ there exists an arrow of sequences
$$\map {\wek \i} {\wek x}{\wek u}.$$
\end{tw}

\begin{pf}
Fix an $\qlim$-continuous sequence $\map{\wek x}\lam \fK$, where $\lam \loe \kappa$.
If $\lam < \kappa$ then we may use $\qlim$ extending the sequence to $\lam + 1$ and we use condition (U) of the definition of a \fra\ sequence for obtaining $\wek \i$.
Thus, we may assume that $\lam = \kappa$.

We shall define inductively a functor $\map I{\truj \kappa} \fK$
such that $I^0 = \wek x$ and $I_\Delta$ is a cofinal subsequence $\wek u \cmp \phi$ of $\wek u$ (so $\map \phi \kappa \kappa$ will be an increasing map).
Using property (U) of the \fra\ sequence, we may assume that $x_0 = I_{\dn00} = u_0$ and we set $\phi(0) = 0$.
Suppose $I$ has been constructed for all pairs $\dn \xi \eta$ with $\xi \loe \eta < \beta$.

Assume first that $\beta =  \al + 1$.
Using $\Ama$-amalgamations, we construct $I_{\dn \xi \beta}$ and the corresponding arrows by induction on $\xi < \al$.
We start by using an $\Ama$-amalgamation, knowing that $I_{\dn 0\beta}$ must be equal to $x_\beta$.
At limit steps we use the fact that $\Ama$ is $\qlim$-continuous.
We finish this procedure at step $\al$.
Finally, using condition (A) of the \fra\ sequence, we find $\phi(\beta) > \phi(\al)$ and an arrow $\map g{I_{\dn \al \beta}}{u_{\phi(\beta)}}$ such that
$$g \cmp I_{\dn \al \al}^{\dn \al \beta} = u_{\phi(\al)}^{\phi(\beta)}.$$
Note that, by the inductive hypothesis, $I_{\dn \al \al} = u_{\phi(\al)}$.
We set
$$I_{\dn \beta \beta} = u_{\phi(\beta)} \oraz I_{\dn \al \beta}^{\dn \beta \beta} = g.$$

Assume now that $\beta$ is a limit ordinal.
We set $\phi(\beta) = \sup_{\xi < \beta}\phi(\xi)$.
For each $\xi < \beta$ let $I_{\dn \xi\beta}$ be the $\qlim$-limit of the sequence that looks as follows.
\begin{equation}
\xymatrix{
I_{\dn \xi \xi} \ar[r] & I_{\dn \xi{\xi+1}} \ar[r] & \dots \ar[r] & I_{\dn \xi \al} \ar[r] & I_{\dn \xi{\al + 1}} \ar[r] & \dots
}
\tag{\decosix}\label{eqwirdf}
\end{equation}
Since $\Ama$ is $\qlim$-continuous, for each $\al < \beta$ there is an arrow from $I_{\dn \al \beta}$ to $I_{\dn {\al + 1} \beta}$ compatible with all arrows below.
Applying the diagonalization property, we obtain $I_{\dn \beta\beta} = u_{\phi(\beta)}$ as the $\qlim$-limit of the sequence (\ref{eqwirdf}).

This finishes the inductive construction.
Finally, since $I_{\dn \al \al} = u_{\phi(\al)}$, the collection of arrows $\wek \i = \sett{I_{\dn 0 \al}^{\dn \al \al}}{\al < \kappa}$ is the required arrow from $\wek x$ to $\wek u$.
\end{pf}

We shall now discuss a class of amalgamation structures for which the proof above can be simplified, avoiding the diagonalization.

\index{amalgamation structure!-- stable}
An amalgamation structure $\Ama$ will be called \emph{stable} if it contains all squares of the form
$$\xymatrix{
b \ar[r]^g & c \\
a \ar[r]_{\id a} \ar[u]^f & a \ar[u]_h
}$$
It turns out that amalgamation structures considered in Section~\ref{reterpeteairs} have this property.

\begin{lm}
Let $\qlim$ be a quasi-limiting operator on $\uzup \kappa \fK$. 
Then every stable $\qlim$-continuous amalgamation structure in $\fK$ has the diagonalization property.
\end{lm}

\begin{pf}
Fix a functor $\map X {\dpower \lam{\loe2}} \fK$ as in condition ($\Delta$).
Given $\beta_0 < \beta < \lam$, the natural composition of the following two diagrams is in $\Ama$:
$$\xymatrix{
X_{\dn{0}{\beta_1}} \ar[rr]^{X_{\dn{0}{\beta_1}}^{\dn{\beta_0}{\beta_1}}} & & X_{\dn{\beta_0}{\beta_1}} \\
& & \\
X_{\dn{0}{\beta_0}} \ar[rr]_{X_{\dn{\beta_0}{\beta_0}}^{\dn{\beta_0}{\beta_0}}} \ar[uu]^{X_{\dn{0}{\beta_0}}^{\dn{0}{\beta_1}}} & & X_{\dn{\beta_0}{\beta_0}} \ar[uu]_{X_{\dn{\beta_0}{\beta_0}}^{\dn{\beta_0}{\beta_1}}}
}
\qquad \qquad
\xymatrix{
X_{\dn{\beta_0}{\beta_1}} \ar[rr]^{X_{\dn{\beta_0}{\beta_1}}^{\dn{\beta_1}{\beta_1}}} & & X_{\dn{\beta_1}{\beta_1}} \\
& & \\
X_{\dn{\beta_0}{\beta_0}} \ar[uu]^{X_{\dn{\beta_0}{\beta_0}}^{\dn{\beta_0}{\beta_1}}} \ar[uurr]_{X_{\dn{\beta_0}{\beta_0}}^{\dn{\beta_1}{\beta_1}}}
}
$$
This is because $\Ama$ is stable and the triangle above can be replaced by a square with an identity arrow at the bottom.
Now the functor $X$ can be reduces to a sequence of squares in $\Ama$ which, by $\qlim$-continuity, converges to an arrow from $\map F_0{\qlim X^0}{\qlim(X_\Delta)}$.
The same holds when replacing $0$ by a fixed positive ordinal $< \lam$, obtaining a co-cone of arrows $\sett{F_\xi}{\xi < \lam}$ for the sequence $\sett{\qlim X^\xi}{\xi < \lam}$.
In particular, $\qlim(\qlim X) = \qlim(X_\Delta)$.
\end{pf}

It is usually much easier to check that a given amalgamation structures is stable, rather than checking the diagonalization property directly (unless $\qlim = \lim$).

Notice that, in the case of stable amalgamation structures, the proof of Theorem~\ref{Tkofnltyrew} can be simplified by constructing an arrow of sequences directly.
\index{arrow of sequences!-- admissible}
More important here is the property of this arrow.
Namely, given an amalgamation structure $\Ama$, an arrow of sequences $\map F {\wek x}{\wek y}$ will be called \emph{$\Ama$-admissible} if it is equivalent to a natural transformation from $\wek x$ to $\wek y \cmp \phi$ (where, as usual, $\phi$ is an order preserving map between the corresponding ordinals) in which all squares belong to $\Ama$.
More precisely, 
all squares of the form
$$\xymatrix{
x_\al \ar[d]_{f_\al} \ar[rr]^{x_\al^\beta} & & x_\beta \ar[d]^{f_\beta} \\
y_{\phi(\al)} \ar[rr]_{y_{\phi(\al)}^{\phi(\beta)}} & & y_{\phi(\beta)}
}$$
are required to be in $\Ama$.
Summarizing, an easy adaptation of the proof of Theorem~\ref{reterpeteairs} gives:

\begin{tw}\label{Thmeotbnoerge}
Assume $\triple \fK \qlim \Ama$ is a $\kappa$-\frajon\ category, where $\kappa = \cf \kappa \goe \aleph_0$ and $\Ama$ is a stable amalgamation structure.
Assume $\wek u$ is an $\qlim$-continuous $\kappa$-\fra\ sequence in $\fK$.
Then for every $\qlim$-continuous sequence $\wek x$ of length $\loe \kappa$ there exists an $\Ama$-admissible arrow $\map {\wek \i} {\wek x} {\wek u}$.
\end{tw}

Clearly, all $\Ama$-admissible arrows provide a subcategory of the category of sequences and the statement above brings more information than Theorem~\ref{reterpeteairs}, even in the case of countable sequences.
We shall use it in Section~\ref{reterpeteairs}, aiming at specific applications.

\section{Some examples}\label{SectExmplsFirstSeries}

This section contains the first series of examples concerning \fra\ sequences, illustrating our ideas and results.

Perhaps one of the simplest is the category of sets $\sets$. If we allow all possible functions as morphisms, then the singleton is trivially a \fra\ object. If we deal with monics only, at the same time restricting to sets of cardinality $<\kappa$, we see that any $\kappa$-chain whose union has cardinality $\kappa$ is \fra.
This, however, is trivial and fits into the classical \frajon\ theory.
On the other hand, starting from the category of finite sets, one can get non-trivial \fra\ sequences. This will be demonstrated in Section~\ref{reterpeteairs}.
Below we present some examples that do not fit completely into the model-theoretic framework.

\subsection{Reversed \fra\ limits}

Here we describe briefly  the theory of inverse limits of models, the so-called \emph{projective \fra} theory, developed by Irwin \& Solecki~\cite{IrSo}.
Fix a first-order language $L$ and consider some class $\Emm$ of $L$-models. The classical \frajon\ theory deals with the category $\Emm^e$
of all embeddings of $L$-models.
\index{embedding of models}
Recall that $\map f M N$ is an \emph{embedding} if it is a one-to-one homomorphism such that for every $n$-ary relation symbol $R \in L$, the following equivalence
$$R^M(f(x_1), \dots, f(x_n)) \iff R^N(x_1, \dots, x_n)$$
holds for every $x_1, \dots, x_n \in M$.
\index{quotient map of models}
In the reversed \frajon\ theory one deals with \emph{quotient maps} $\map f M N$ that are, by definition, surjective homomorphisms satisfying the following formula
$$R^N(y_0, \dots, y_{n-1}) \implies (\exists\; x_0, \dots, x_{n-1} \in M)\; R^M(x_0,\dots, x_{n-1}) \land (\forall\; i < n)\; f(x_i) = y_i$$
for each $n$-ary relation symbol $R \in L$ and for each $y_0, \dots, y_{n-1} \in N$.
The two notions are dual in some sense, due to the following easy fact.

\begin{prop}
Let $\map f A B$ and $\map g B C$ be homomorphisms of models of a fixed first-order language. Then
\begin{enumerate}
	\item[(a)] $g \cmp f$ is an embedding $\implies$ $f$ is an embedding.
	\item[(b)] $g \cmp f$ is a quotient map $\implies$ $g$ is a quotient map.
\end{enumerate}
\end{prop}

The following result from \cite{IrSo} is a  direct application of the results from Section~\ref{SecFraMejnone}.
It can also be derived from the results of \cite{DrGoe92}, however our approach is more direct and explains why the topology is needed.

\begin{tw}\label{Twenrgoero}
Let $\Emm$ be a countable class of finite models of a fixed first-order language $L$. Suppose $\Emm$ satisfies the following conditions:
\begin{enumerate}
	\item[(J)] For every $a,b \in \Emm$ there exist $c\in \Emm$ and quotient maps $\map f c a$ and $\map g c b$.
	\item[(A)] Given quotient maps $\map f c a$ and $\map g c b$ with $a,b,c \in \Emm$, there exist $w \in \Emm$ and quotient maps $f'$, $g'$ for which the diagram
$$\xymatrix{
b \ar[d]_g & w \ar[l]_{g'} \ar[d]^{f'} \\
c & a \ar[l]^f
}$$
commutes.
\end{enumerate}
Then there exists a unique (up to a topological isomorphism) topological $L$-model $\Em$ satisfying the following conditions:
\begin{enumerate}
	\item[(1)] $\Em$ is the inverse limit of a sequence of models from $\Emm$ with quotient maps.
	\item[(2)] Every model from $\Emm$ is a continuous quotient of $\Em$.
	\item[(3)] Given continuous quotients $\map p \Em a$ and $\map q \Em b$ with $a, b\in \Emm$, given a quotient map $\map f b a$, there exists a topological isomorphism $\map h \Em \Em$ for which the diagram
$$\xymatrix{
\Em \ar[d]_p & \Em \ar[l]_{h} \ar[d]^{q} \\
a & b \ar[l]^f
}$$
is commutative.
\end{enumerate}
Furthermore, every $L$-model which is the inverse limit of a sequence of models from $\Emm$ with quotient maps is a continuous quotient of $\Em$.
\end{tw}

\begin{pf}
Let $\fK$ be the category whose objects are elements of $\Emm$ and an arrow from $a \in \Emm$ into $b \in \Emm$ is a quotient map $\map f b a$.
It is obvious that conditions (J) and (A) translate to directedness and the amalgamation property.
Since $\Emm$ is countable, the category $\fK$ is countable.
Now observe that the category $\sig \fK$ can be naturally identified with the category of compact $L$-models that are inverse limits of sequences of models from $\Emm$.
In fact, the topology becomes natural here, because given two sequences $\wek x$ and $\wek y$ in $\fK$ and taking $X$ and $Y$ to be their inverse limits in the category of sets, one can easily check that precisely the continuous quotient maps $\map f Y X$ correspond to $\ciagi \fK$-arrows from $\wek x$ to $\wek y$.
Summarizing: the existence and properties of $\Em$ follow directly from Theorems~\ref{jofjaiopf}, \ref{saofjapiwf} and \ref{jednorodnost}.
\end{pf}

The main problem to get the reversed \fra\ limit is, as usual, showing that the class of models in question has the reversed amalgamation property, i.e., satisfies condition (A) above.
In fact, in \cite{IrSo} this property for the class of linear graphs is crucial and non-trivial.
It turns out however, that for many classes of models the reversed amalgamation can be proved easily.

\begin{lm}\label{Lnrjner}
Let $\Emm$ be a class of models of a first-order language $L$. Assume $\Emm$ is closed under finite products and substructures.
Then for every quotient maps $\map f c a$ and $\map g c b$ such that $a,b,c \in \Emm$ there exist quotient maps $\map {f'} a w$ and $\map {g'} b w$ for which the diagram
$$\xymatrix{
b \ar[d]_g & w \ar[l]_{g'} \ar[d]^{f'} \\
c & a \ar[l]^f
}$$
is commutative.
\end{lm}

\begin{pf}
Define
$$w = \setof{\pair x y \in a \times b}{f(x) = g(y)}.$$
It is straight to check that $w$ is a submodel of the product $a \times b \in \Emm$.
We let $f'$ and $g'$ to be the canonical projections.
It is obvious that they are homomorphisms.
It remains to check that $f'$, $g'$ are quotient maps.
Fix an $n$-ary relation symbol $R$ in $L$ and fix $x_0, \dots, x_{n-1} \in a$ such that $R(x_0, \dots, x_{n-1})$ holds.
Since $R(f(x_0), \dots, f(x_0))$ holds and $g$ is a quotient map, there exist $b_0, \dots, b_{n-1} \in b$ such that $R(b_0, \dots, b_{n-1})$ holds and $g(b_i) = f(x_i)$ for $i < n$.
Let $p_i = \pair {x_i}{b_i}$ for $i < n$. 
Note that $\sett{p_i}{i < n} \subs w$ and $f'(p_i) = x_i$ for $i < n$.
Finally. $R(p_0, \dots, p_{n-1})$ holds.
This shows that $f'$ is a quotient map.
The same argument shows that $g'$ is a quotient map.
\end{pf}

The statement above is particularly useful for some classes of algebras, like groups, semilattices, rings, etc.
Recall that a topological algebra is \emph{pro-finite} if it is the inverse limit of a system of finite algebras with epimorphisms.
Lemma~\ref{Lnrjner} and the results above give the following:

\begin{wn}
Let $\Emm$ be a class of pro-finite algebras, stable under subalgebras and finite products. Assume that $\Emm$ contains countably many isomorphic types of finite algebras.
Then there exists a unique, up to a topological isomorphism, second countable pro-finite algebra $\Aa$ satisfying conditions (2) and (3) from Theorem~\ref{Twenrgoero}.
Furthermore, every second countable pro-finite algebra from $\Emm$ is a continuous quotient of $\Aa$.
\end{wn}

\subsection{Monoids}

\index{semigroup}\index{monoid}
Recall that a \emph{semigroup} is a structure of the form $\triple S \cmp 1$, where $\cmp$ is an associative binary operation on $S$, and $1\in S$ is such that $1 \cmp s = s = s \cmp 1$ holds for every $s \in S$.
A semigroup $\fS = \triple S \cmp 1$ can be viewed as a category with a single object $S$, whose arrows are all elements of $S$ and $1$ is the identity arrow.
Such a category is called a \emph{monoid}.
A sequence in a monoid $\fS$ is naturally represented as $\ciag a$, where $a_n \in S$ for every $\ntr$.
This corresponds to the formal sequence
$$\xymatrix{
S \ar[r]^{a_0} & S \ar[r]^{a_1} & S \ar[r] & \dots \ar[r] & S \ar[r]^{a_n} & S \ar[r] & \dots
}$$
In particular, the bonding arrow from the $n$th to the $m$th element of the sequence is the composition $a_{m-1} \cmp \dots \cmp a_n$.

Note that every monoid is automatically directed.
The amalgamation property can be rephrased as
$$(\forall\; f, g)(\exists\; f', g')\; f' \cmp f = g' \cmp g.$$
This is automatically fulfilled when $S$ is commutative.
In particular, every countable commutative monoid has a \fra\ sequence of length $\omega$.

It turns out that, even the set of natural numbers provides a monoid with many non-isomorphic $\omega$-sequences.
Namely, consider $\Nat^+$, the set of all positive integers, as a multiplicative monoid.
Formally, a sequence in $\Nat^+$ is a function $\map \phi \Delta \Nat^+$, where
$\Delta = \setof{\pair mn \in \nat \times \nat}{m \loe n}$,
satisfying $\phi(n,n) = 1$ and $\phi(n,k) \cdot \phi(m,n) = \phi(m,k)$ for every $m < n < k$.
Clearly, the relevant values are $\phi(n,n+1)$, therefore $\ciagi{\Nat^+}$ can be naturally identified with sequences of positive integers.
Now, it is an easy exercise that a sequence $\ciag u$ is \fra\ if and only if for every prime $p \in \Nat^+$, for every $\ntr$, there is $m>n$ such that $p$ divides $u_n u_{n+1} \cdot \dots \cdot u_m$.
Since there are infinitely many prime numbers, it is easy to see that $\ciagi \Nat^+$ has continuum many pairwise incomparable sequences.
This is not the case for the additive monoid $\Nat$, where there are only two isomorphic types in $\ciagi \Nat$, namely, the equivalence class of the constant 0 sequence and the equivalence class of the constant 1 sequence.
The latter one is \fra\ in $\Nat$.

Another class of monoids are semilattices.
Namely, fix a join semilattice $\triple L \join 0$ and let $\loe$ be the partial order induced by $\join$, i.e. $x \loe y$ iff $x\join y = y$.
Since this is a commutative monoid, it has the amalgamation property.
Let $\map {\wek u} \kappa L$ be a sequence in $L$.
Given $\al < \beta < \kappa$ we have that $u_0^\beta = u_\al^\beta \join u_0^\al$, therefore $u_0^\al \loe u_0^\beta$.
Assuming that $\wek u$ is \fra\ and using condition (A), we see that for every $x \in L$ there are $\al < \kappa$ and $y\in L$ such that $u_0^\al = y \join x$.
In other words, if $\wek u$ is \fra, then $\sett{u_0^\al}{\al < \kappa}$ is an increasing cofinal well ordered chain in $\pair L \loe$.
Conversely, given an increasing cofinal sequence $\sett{a_\al}{\al < \kappa}$ in $\pair L \loe$, the functor $\map {\wek u}\kappa L$, defined by $u_\al^\al = 0$ and $u_\al^\beta = a_\beta$ for $\al < \beta < \kappa$, is a \fra\ sequence in $L$.

An important class of monoids are transformation semigroups.
We briefly discuss the particular case of such monoids induced by \fra\ sequences.
Fix a category $\fK$ such that $\ciagi\fK$ has the amalgamation property and assume $\wek u$ is a countable \fra\ sequence in $\fK$.
Consider $M_{\wek u} = \ciagi \fK(\wek u, \wek u)$.
By Theorem~\ref{ThmStbFrcsest}, $M_{\wek u}$ is a $\sig$-complete monoid.
It turns out that $M_{\wek u}$ ``encodes" all information about the category $\sig\fK$ relevant for the existence of an $\omega_1$-\fra\ sequence:

\begin{tw}
Let $\fK$ be a category with an $\omega$-\fra\ sequence $\wek u$ and such that $\ciagi \fK$ has the amalgamation property.
Then $M_{\wek u}$ has the amalgamation property and is cofinal in $\ciagi \fK$.
Furthermore, $\ciagi \fK$ has an $\omega_1$-\fra\ sequence if and only if the monoid $M_{\wek u}$ has an $\omega_1$-\fra\ sequence.
\end{tw}

Note that, by the first part, a \fra\ sequence in $M_{\wek u}$ is \fra\ in $\ciagi \fK$.

\begin{pf}
The fact that $M_{\wek u}$ has the amalgamation property and is cofinal in $\ciagi \fK$ follows directly from Theorem~\ref{saofjapiwf}.
It remains to prove that if $\ciagi\fK$ has an $\omega_1$-\fra\ sequence then so does $M_{\wek u}$.
For this aim, it suffices to show that $M_{\wek u}$ is dominated by $\loe \aleph_1$ many arrows.
Fix a \fra\ sequence $\map {\omega_1} {\wek {\wek v}} {\ciagi \fK}$.
Cutting the sequence if necessary, we may assume that there is a $\ciagi\fK$-arrow $\map H {\wek u}{\wek v_0}$.
Using Theorem~\ref{saofjapiwf}, for each $\al<\omega_1$ choose an arrow $\map {F_\al}{\wek v_\al}{\wek u}$.
Let $G_\al$ denote the bonding arrow from $\wek v_0$ to $\wek v_\al$.
Define
$$\Ef = \setof{ F_\al \cmp G_\al \cmp H }{\al < \omega_1}.$$
We claim that $\Ef$ is dominating in $M_{\wek u}$.
Fix a $\ciagi \fK$-arrow $\map f {\wek u} {\wek u}$.
Using the amalgamation property, we find $\ciagi \fK$-arrows $H'$, $f'$ such that $f' \cmp f = H' \cmp H$.
Using the fact that $\wek{\wek v}$ is \fra, we find $\al > 0$ and a $\ciagi \fK$-arrow $g$ such that $g \cmp H' = G_\al$.
Finally, $(F_\al \cmp g \cmp f') \cmp f = F_\al \cmp G_\al \cmp H \in \Ef$.
\end{pf}

It is clear that the result above can be generalized to uncountable sequences, adding suitable assumptions concerning completeness.

\subsection{Diagrams: the role of pushouts}

We shall describe a quite general procedure of building a category with the amalgamation property.
First of all, let us note that typical categories with monics do not admit pushouts. For instance, this is the case with the category of finite sets with one-to-one maps.
In order to capture these situations, we define the following natural notion.

Let $\fK \subs \fL$ be two categories. We say that $\fK$ \emph{has pushouts in} $\fL$ (or, that \emph{$\pair \fK \fL$ has pushouts})
if for every $\fK$-arrows $\map fca$, $\map gcb$ there exist $\fK$-arrows $\map {f'}aw$ and $\map {g'}bw$ such that
$$\xymatrix{
b \ar[r]^{g'} & w \\
c \ar[r]_f \ar[u]^g & a \ar[u]_{f'}
}$$
is a pushout in $\fL$.

Now, fix a small category $\bS$ and let $\fun \bS \fK \fL$ be the category whose objects are covariant functors from $\bS$ to $\fL$ taking objects from $\fK$ and the arrows are natural transformations into $\fK$.
More precisely, $x$ is an object of iff $x$ is a functor from $\bS$ into $\fL$ such that $x(s)$ is an object of $\fK$ for every object $s$ of $\bS$.
Furthermore, $\map fxy$ is an $\fun \bS \fK \fL$-arrow iff it is a natural transformation from $x$ to $y$ whose all components (arrows) are in $\fK$.
In other words, for every object $s$ of $\bS$, the arrow $f(s)$ belongs to $\fK$.

\begin{lm}\label{Luniwdiagramsewrg}
Let $\fK \subs \fL$ be a pair of categories such that $\fK$ has pushouts in $\fL$ and let $\bS$ be a small category. Then $\fun \bS \fK \fL$ has pushouts in $\fun \bS \fL \fL$.
\end{lm}

\begin{pf}
Fix $a,b,c\in\fun \bS \fK \fL$ and fix natural transformations $\map fca$, $\map gcb$.
We are going to define $w\in \fun \bS \fK \fL$. Fix $p\in\bS$. Let $\pair{f'(p)}{g'(p)}$ be the pushout of $\pair{f(p)}{g(p)}$ in $\fL$, and let $w(p)$ be the common co-domain of $f'(p)$ and $g'(p)$.
Now fix an arrow $\map jpq$ in $\bS$. Using the property of a pushout, there is a unique arrow $w(j)$ making the following diagram commutative.
$$\xymatrix{
c(q)\ar[rrrr]^{g(q)}\ar[rd]_{f(q)} & & & & b(q)\ar[rd]^{g'(q)} & \\
& a(q)\ar[rrrr]^{f'(q)} & & & & w(q) \\
& & & & & \\
c(p)\ar[uuu]^{c(j)}\ar[rd]_{f(p)}\ar[rrrr]^{g(p)} & & & & b(p)\ar[uuu]_{b(j)}\ar[dr]^{g'(p)} & \\
&  a(p)\ar[uuu]^{a(j)}\ar[rrrr]^{f'(p)} & & & & w(p)\ar@{.>}[uuu]_{w(j)}
}$$
Indeed, $\pair{f'(q) \cmp a(j)}{g'(q) \cmp b(j)}$ is an amalgamation of $\pair{f(p)}{g(p)}$.
By uniqueness, we have that $w(j\cmp k) = w(j)\cmp w(k)$, whenever $j,k$ are compatible arrows in $\bS$.

Thus we have defined a functor $\map w\bS\fK$. Further, $p\mapsto f'(p)$ and $p\mapsto g'(p)$ are natural transformations from $a$ to $w$ and from $b$ to $w$ respectively.

It is clear that $f' \cmp f = g' \cmp g$.

It remains to check that $\pair{f'}{g'}$ is a pushout of $\pair fg$ in $\fun \bS \fL \fL$. 
For this aim, fix $v\in \fun \bS \fL \fL$ and natural transformations $\map {f''}av$, $\map {g''}bv$ such that $f''\cmp f = g''\cmp g$.
By the fact that $\pair{f'}{g'}$ is a pushout of $\pair fg$ in $\fL$, for each $p\in\bS$ there is a unique $\fL$-arrow $\map{h(p)}{w(p)}{v(p)}$ satisfying $h(p)\cmp f'(p) = f''(p)$ and $h(p)\cmp g'(p) = g''(p)$.
This defines uniquely a map $\map h\bS \Arr \fL$ that satisfies $h\cmp f' = f''$ and $h\cmp g' = g''$. It remains to check that $h$ is indeed a natural transformation.

Fix an arrow $\map jpq$ in $\bS$ and let $k = h(q)\cmp w(j)$, $\ell = v(j)\cmp h(p)$. We need to show that $k=\ell$.
Notice that $k\cmp f'(p) = f''(q)\cmp a(j)$ and $k\cmp g'(p) = g''(q)\cmp b(j)$. Also, $\ell\cmp f'(p) = f''(q)\cmp a(j)$ and $\ell\cmp g'(p) = g''(q)\cmp b(j)$.
It follows that $k=\ell$, because $\pair{f'(p)}{g'(p)}$ is a pushout of $\pair{f(p)}{g(p)}$ in $\fL$.
\end{pf}

The interesting part of the lemma above is that $\fun \bS \fK \fL$ has the amalgamation property, whenever $\fK$ has pushouts in $\fL$.

We now demonstrate a possible use of Lemma~\ref{Luniwdiagramsewrg} in the context of \fra\ model-theoretic structures.

\index{mixed pushout property}
Given categories $\fK \subs \fL$, we say that $\pair \fK \fL$ has the \emph{mixed pushout property} if for every $\fK$-arrow $\map ica$, for every $\fL$-arrow $\map fcb$, there exist $j\in \fK$ and $g\in \fL$ such that
$$\xymatrix{
b \ar[r]^j & w \\
c \ar[r]_i \ar[u]^f& a \ar[u]_g
}$$
is a pushout in $\fL$.
In case $\Emm$ is a class of models of some first-order language, it is natural to consider $\fK(\Emm)$ to be the category of all embeddings between models of $\Emm$ and $\fL(\Emm)$ to be the category of all homomorphisms between these models.
We then say that $\Emm$ has the \emph{mixed pushout property} if so does $\pair{\fK(\Emm)}{\fL(\Emm)}$.
We denote by $\ovr \Emm$ the class of all (countable) models that are unions of $\omega$-chains of models from $\Emm$.

\begin{tw}\label{ThmUnivHomBX}
Let $\Emm$ be a countable \fra\ class of finitely generated models, with the mixed pushout property.
Let $W$ denote the \fra\ limit of $\Emm$.
Then there exists a unique (up to isomorphism) homomorphism $\map L W W$ satisfying the following conditions.
\begin{enumerate}
	\item[(a)] For every $X,Y\in \ovr\Emm$, for every homomorphism $\map F X Y$ there exist embeddings $\map {I_X} X W$ and $\map {I_Y} Y W$ such that the square
$$\xymatrix{
W \ar[r]^L & W \\
X \ar[u]^{I_X} \ar[r]_F & Y \ar[u]_{I_Y}
}$$
is commutative.
	\item[(b)] Given finitely generated substructures $x_0,x_1,y_0,y_1$ of $W$ such that $\img L{x_i} \subs y_i$ for $i < 2$, given isomorphisms $\map {h_i}{x_i}{y_i}$ for $i < 2$ such that $L \cmp h_0 = h_1 \cmp L$, there exist automorphisms $\map {H_i} W W$ extending $h_i$ for $i < 2$, and such that $L \cmp H_0 = H_1 \cmp L$.
\end{enumerate}
\end{tw}

\begin{pf}
Consider the category $\fC = \fun 2 {\fK(\Emm)}{\fL(\Emm)}$, where $2$ denotes the two-element poset category.
The assumptions above, combined with Lemma~\ref{Luniwdiagramsewrg}, give a \fra\ sequence $\ciag \ell$ which translated back to $\fL(\Emm)$ looks as follows.
$$\xymatrix{
u_0 \ar[d]_{\ell_0} \ar[r] & u_1 \ar[d]_{\ell_1} \ar[r] & \dots \ar[r] & u_n \ar[d]_{\ell_n} \ar[r] & \dots \\
v_0 \ar[r] & v_1 \ar[r] & \dots \ar[r] & v_n \ar[r] & \dots \\
}$$
The horizontal arrows are embeddings and the vertical arrows are homomorphisms of models.
Without loss of generality, we may assume that the horizontal arrows are inclusions.
Now let $U = \bigcup_{\ntr}u_n$, $V = \bigcup_{\ntr}v_n$ and define $L = \bigcup_{\ntr}\ell_n$.
It is clear that $L$ satisfies conditions (a) and (b) above.
It remains to check that both $U$ and $V$ are isomorphic to the $\Emm$-\fra\ limit $W$.
For this aim, it suffices to show that both $\ciag u$ and $\ciag v$ are \fra\ sequences in $\fK(\Emm)$.

It is clear that both sequences satisfy (U), because identity arrows can be viewed as objects of $\fC$.

Fix $\ntr$ and fix an embedding $\map j {u_n}x$, where $x\in \Emm$.
Using the mixed pushout property, we can find an embedding $\map k {v_n}y$ and a homomorphism $\map f x y$ such that $k \cmp \ell_n = f \cmp j$.
Now $f$ can be regarded as an object of $\fC$ and the pair $\pair jk$ is an arrow of $\fC$. Using the \fra\ property of the sequence $\ciag \ell$, we find $m \goe n$ and embeddings $\map {j'}x {u_m}$, $\map {k'}y {v_m}$ such that $\ell_m \cmp j' = k' \cmp f$ and $j' \cmp j$, $k' \cmp k$ are inclusions.
In particular, this shows that $\ciag u$ is \fra\ and therefore $U = W$.

Now fix an embedding $\map i {v_n}z$.
Let $g = i \cmp \ell_n$.
The pair $\pair {\id{u_n}} i$ can be viewed as a $\fC$-arrow into $g$, therefore by the \fra\ property of $\ciag \ell$ we again find $m \goe n$ and embeddings $\map e{u_n}{u_m}$, $\map {i'}z {u_m}$ such that, among other conditions, the composition $i' \cmp i$ is the inclusion $v_n \subs v_m$.
This shows that $\ciag v$ is \fra\ and consequently $V = W$.
\end{pf}

Note that the \emph{mixed} pushout property is needed for concluding that the domain of the universal homogeneous homomorphism is the \fra\ limit of the class in question. Relaxing this by assuming only that embeddings have pushouts among all homomorphisms, one still gets a homomorphism with properties (a) and (b) above, whose co-domain is the \fra\ limit of the considered class of models.
It is clear how to formulate a similar statement for reversed \fra\ limits, where pushouts are replaced by pullbacks.

It might be interesting to compare Theorem~\ref{ThmUnivHomBX} with the recent work of Pech \& Pech \cite{Pech}, where the authors consider comma categories, in particular obtaining universal homomorphisms as \fra\ limits, constructed from categories with mixed pushouts.

\subsection{Universal homogeneous categories}

It is tempting to ask for the existence of a \fra\ sequence in ``categories of all categories".
Of course, one needs to be consistent with the axioms of set theory, therefore we should only talk about classes of small categories.
Actually, a small category is nothing but a directed graph endowed with an additional ``composition" operation, satisfying the obvious axioms.
Arrows between small categories are just covariant functors.
In order to avoid trivialities, we should consider injective functors (possibly onto full subcategories).
Summarizing, we can ask for a \fra\ sequence in a directed class of categories that has the amalgamation property with respect to injective functors.
\index{functor!-- injective}
Recall that a functor is \emph{injective} if it is one-to-one on the class of arrows (in particular, such a functor must be one-to-one on the class of objects).

Actually, this question had already been addressed by Trnkov\'a~\cite{Trnkova} in 1966. Her work was set up within the Bernays-G\"odel axiom system for set theory.
In this framework, every category is the union of a chain of small categories and the amalgamation property (which appears to be quite non-trivial---proved by Trnkov\'a in her earlier work~\cite{TrnkovaAmalgams}) allows concluding the existence of a ``\fra\ chain" of categories.
We can add more information to the results of \cite{Trnkova} by adding homogeneity.
Namely, every injective functor between small subcategories extends to a bijective functor of the big universal category under consideration.
Since we work in the usual Zermelo-Fraenkel set theory (with the axiom of choice), we discuss briefly one obvious class of categories.
Other examples (like additive, concrete categories, and so on) can be easily adapted from Trnkov\'a's work~\cite{Trnkova}, restricting the size of categories and adding some cardinal-arithmetic assumptions.

A natural example is the class $\fK$ of all \emph{finite categories}.
A category $\fK$ is \emph{finite} if it has finitely many isomorphic types and each hom-set of $\fK$ is finite.
It is rather clear that this class has the amalgamation property.
It also has the initial object (the empty category), therefore it is directed.
Obviously, $\fK$ is dominated by a countable family of injective functors, say, between finite categories whose objects are natural numbers and arrows are taken from a prescribed countable set.
By Theorem~\ref{jofjaiopf}, $\fK$ has a unique countable \fra\ sequence.

Let $\fU$ be the category obtained as the ``co-limit" of the \fra\ sequence in $\fK$.
Then $\fU$ is a category with countably many isomorphic types of objects and each hom-set of $\fU$ is countable.
It is an easy exercise to see that every other category with these properties is equivalent to a subcategory of $\fU$.
Finally, $\fU$ has the following homogeneity property:
Given an injective functor $\map F \fS \fT$, between finite subcategories of $\fU$, there exists a bijective functor $\map H \fU \fU$ extending $F$.

One can consider a similar construction using injective functors $\map F \fS \fT$ such that $\img F \fS$ is a full subcategory of $\fT$.
The \fra\ limit would be a countable category $\fU_f$ whose all hom-sets are finite, all other countable categories with this property would be equivalent to full subcategories of $\fU_f$ and the obvious variant of homogeneity would hold.
We leave the details to interested readers.

\subsection{Binary trees: failure of the back-and-forth principle}\label{drzewa}

We describe below the announced example of a category of trees which has many pairwise non-equivalent $\omega_1$-\fra\ sequences.

Recall that a tree $T$ is \emph{bounded} if for every $x\in T$ there is $t\in \max T$ such that $x\loe t$.
A subset $A$ of $T$ is \emph{closed} if $\sup C\in A$ for every chain $C\subs A$.
This is equivalent to saying that $A$ is closed with respect to the \emph{interval topology} on $T$ generated by intervals of the form $[0,t]$ and $(s,t]$, where $s<t$.
\index{tree!-- closed subset}\index{tree!-- interval topology}

We define the category $\trees$ as follows.
\index{category!-- \trees}\index{\trees}
The objects of \trees\ are nonempty countable bounded binary trees. An arrow from $T\in\trees$ into $S\in\trees$ is a tree (i.e. strictly order preserving) embedding $\map fTS$ such that $\img fT$ is a closed initial segment of $S$.

A tree $T$ is \emph{healthy}
\index{tree!-- healthy}\index{healthy tree}
if every element of $T\setminus\max(T)$ has at least two immediate successors and for every $t\in T$ and $\al<\Ht(T)$ there exists $s\goe t$ such that $\Lev_T(s)\goe\al$. An example of a healthy tree of height $\omega_1$ is
$$T = \setof{x\in 2^{<\omega_1}}{ |\setof{\al}{x(\al)=1}|<\aleph_0}.$$
Note that all levels of $T$ are countable. Setting $T_\al=\setof{x\in T}{\dom(x)\subs\al+1}$, we obtain an inductive sequence $\sett{T_\al}{\al<\omega_1}$ in the category \trees. More generally, if $S$ is any binary tree of height $\omega_1$ whose all levels are countable then, setting
$$S_\al=\setof{x\in S}{\text{the order type of $[0,x)$ is }\loe\al}$$ we obtain an inductive sequence $\wek S$ in \trees, where each $S_\al^\beta$ is the inclusion, which is an arrow in \trees. We shall say that $\wek S$ is the \emph{natural decomposition} of $S$.
\index{tree!-- natural decomposition of}

\begin{lm}\label{woehjqt}
Let $V\in\trees$ be a healthy tree of height $\al+1$. Then every $T\in \trees$ with $\Ht(T)\loe\al+1$ is isomorphic to a closed initial segment of $V$.
\end{lm}

\begin{pf}
Denote by $\fM$ the class of all nonempty bounded binary trees $T$ of height $\loe\al+1$ such that $\max(T)$ is finite. Given such a tree $T$, write $\max(T)=\{w_0, \dots, w_{m-1}\}$ and define inductively $T_0:=[0,w_0]$ and $T_k := [0,w_k]\setminus(T_0\cup\dots\cup T_{k-1})$. Then $\Dee=\{T_0,\dots,T_{m-1}\}$ is a natural decomposition of $T$ into connected chains, induced by the enumeration of $\max(T)$. 

Let $\fL$ be a category whose objects are pairs $\pair T\Dee$, where $T\in\fM$ and $\Dee$ is a natural decomposition into connected chains induced by an enumeration of $\max(T)$, as described above. Given $\pair T\Dee, \pair S\Ef\in\fL$, an arrow in $\fL$ is a tree embedding $\map fTS$ with the following properties: 
\begin{enumerate}
  \item[(a)] For each $D\in \Dee$ there is $F\in\Ef$ such that $\img fD$ is an initial segment of $F$.
  \item[(b)] For each $F\in \Ef$ there exists at most one $D$ such that $\img fD\subs F$.
\end{enumerate}
It is clear that these properties are preserved under the usual composition, so $\fL$ is indeed a category.
Now consider the given healthy tree $V$ with $\max(V)=\ciag e$ and define $V^n = \bigcup_{i\loe n}[0,e_i]$. 
Let $\Dee^n$ be the decomposition of $V^n$ induced by the enumeration $\{e_0,\dots,e_n\}$ of $\max(V^n)$.
Then $\pair{V^n}{\Dee^n}\in \fL$ and $\wek V = \sett{\pair{V^n}{\Dee^n}}{\ntr}$ is an inductive sequence in $\fL$. We claim that:
\begin{enumerate}
	\item[(1)] $\wek V$ is a \fra\ sequence in $\fL$ which has property (E).
	\item[(2)] If $\wek T$ is an inductive sequence in $\fL$ and $\wek f=\ciag f$ is an embedding of $\wek T$ into $\wek V$ then the embedding $\map{f}TV$ induced by $\wek f$ has the property that $\img fT$ is a closed initial segment of $V$.
\end{enumerate}
We first show (2): Fix $y\in V\setminus \img fT$. Find $m$ and $D\in\Dee_m$ such that $y\in D$. 
Let $\Es_n$ be the natural decomposition of $T_n$, $\Es=\bigcup_{\ntr}\Es_n$. Then there is at most one $S\in \Es$ such that $\img fS\subs D$. Moreover $\img fS$ is closed in $D$, so $D\setminus \img fS$ is a neighborhood of $y$ disjoint from $\img fT$.

For the proof of (1), fix $\pair T\Ef\in\fL$ and assume $\map fT{V_n}$ is an arrow in $\fL$. Let $T\subs T'$ and let $\Ef'\sups \Ef$ so that the inclusion $T\subs T'$ is an arrow between $\pair T\Ef$ and $\pair{T'}{\Ef'}$. Without loss of generality, we may assume that $\Ef'=\Ef\cup\sn A$, i.e. $T'$ differs from $T$ by only one new branch. Let $a=\min A$. Then, by the definition of the natural decomposition, $a$ has an immediate predecessor $c\in T$ (note that $0\in T$ so $a>0$). Find $F\in\Ef$ such that $c\in F$. Then $f(c)$ has exactly two immediate successors in $V$ and at most one belongs to $\img fT$, since $c$ has only two immediate successors in $T$. Let $d\in V\setminus \img fT$ be an immediate of $f(c)$. Find a big enough $m>n$ so that there exists $D\in \Dee_m$ with $d\in D$. Since each maximal element of $V$ has height $\al$, $D$ is a cofinal branch in $V$ and hence $A$ can be (uniquely) embedded into $D$ as an initial segment. This  embedding defines an extension $\map{\ovr f}{T'}{V_m}$ of $f$. Since $\fL$ has a minimal object, this shows that $\wek V$ is \fra\ and satisfies (E).

Finally, fix $T\in\trees$ with $\Ht(T)\loe\al+1$. Decompose $T$ into an inductive $\omega$-sequence, according to a fixed enumeration of $\max(T)$. Claims (1) and (2) say that $T$ can be embedded into $V$ as a closed initial subtree. This completes the proof.
\end{pf}

\begin{tw}\label{ofhqofjh}
Assume $U$ is a healthy binary tree of height $\omega_1$, whose all levels are countable. Let $\wek U$ be the natural $\trees$-decomposition of\/ $U$. Then $\wek U$ is a \fra\ sequence in $\trees$ which has the extension property. In particular, $\trees$ is directed and has the amalgamation property.
\end{tw}

\begin{pf}
Note that $\trees$ has a minimal object, namely the one-element tree. Clearly, such a tree embeds into $U_0$. It suffices to show that $\wek U$ satisfies (E), it will then follow that $\wek U$ is a \fra\ sequence. Since there exists a healthy binary tree of height $\omega_1$ with countable levels, we shall be able to conclude that $\trees$ has a \fra\ sequence satisfying (E) and consequently $\trees$ is directed and has the amalgamation property.
Thus, it remains to show that $\wek U$ has the extension property.

Fix $\al<\omega_1$ and fix an arrow $\map fT{U_\al}$ in $\trees$ and assume that $T$ is a closed initial subtree of $S$, i.e. the inclusion $T\subs S$ is an arrow of $\trees$. 
Fix $\al<\omega_1$ so that $\Ht(S)<\al$. Let $\setof{s_n}{\ntr}$ enumerate all minimal elements of $S\setminus T$. Let $t_n$ be the immediate predecessor of $s_n$. Then $t_n\in T$. Recall that $f(t_n)$ has exactly two immediate successors in $U_\al$ and at least one of them does not belong to $\img fT$, since otherwise $t_n$ would already have two immediate successors in $T$. Let $y_n$ be an immediate successor of $f(t_n)$ which does not belong to $\img fT$. Let $V_n=\setof{y\in U_\al}{y\goe y_n}$. Then $V_n$ is a healthy binary tree of countable height. By Lemma \ref{woehjqt} we can embed $G_n=\setof{s\in S}{s\goe s_n}$ onto a closed initial segment of $V_n$. Combining all these embeddings, we obtain an extension $\map{\ovr f}S{U_\al}$ of $f$. We claim that $\img{\ovr f}S$ is closed in $U_\al$. Indeed, if $C\subs S$ is a chain and $C\not\subs T$ then $c\goe s_n$ for some $c\in C$ and for some $n$. Thus $\sup C$ exists in $S$ and hence $\sup\img{\ovr f}C$ exists in $\img{\ovr f}{G_n}\subs\img{\ovr f}S$.
\end{pf}

\begin{tw}\label{efhfhqo}
Assume $\wek U$ and $\wek V$ are $\omega_1$-\fra\ sequences in \trees\, inducing healthy trees $U$ and $V$ respectively. Assume further that $\map F{\wek U}{\wek V}$ is an arrow of sequences. Then the trees $U$, $V$ are isomorphic.
\end{tw}

\begin{pf}
Let $\map fUV$ be the embedding induced by $F$, i.e. assuming both $\wek U$, $\wek V$ are chains of trees, $f$ is the union of arrows $\map{f_\al}{U_\al}{V_{\phi(\al)}}$ in \trees, where $\map \phi{\omega_1}{\omega_1}$ is an increasing function. We claim that $\img f{U}$ is closed in $V_{\phi(\al)}$. 
Suppose otherwise and fix a sequence $x_0<x_1<\dots$ in $U$ such that $y=\sup_{\ntr}f(x_n)\notin \img fU$. Find $\beta<\omega_1$ such that $\setof{x_n}{\ntr}\subs U_\beta$ and $y\in V_{\phi(\beta)}$. Then $f(x_n)=f_\beta(x_n)$ and $y\notin\img{f_\beta}{U_\beta}$, which shows that $f_\beta$ is not an arrow in \trees, a contradiction.

We finally claim that $V=\img fU$, which of course shows that $f$ is an isomorphism. Suppose $V\ne \img fU$ and fix a minimal element $y\in V\setminus \img fU$. Find $\al<\omega_1$ such that $y\in V_{\phi(\al)}$. Since $\img fU$ is closed in $V_{\phi(\al)}$, $y$ has an immediate predecessor, say $v=f(u)$. Let $a,b$ be the two immediate successors of $u$ in $U$ (which exist, because $U$ is healthy). Then either $y=f(a)$ or $y=f(b)$, because $V$ is binary and $\img fU$ is an initial segment of $U$. This is a contradiction.
\end{pf}

It can be easily shown that every tree induced by an $\omega_1$-\fra\ sequence in \trees\ is healthy, therefore this assumption can be removed from the statement above.

Recall that sequences $\wek a$ and $\wek b$ of the same length are \emph{comparable}
\index{sequence!-- (in)comparability}
if there exists an arrow of sequences ${\wek f}$ such that either $\map {\wek f}{\wek a}{\wek b}$ or $\map{\wek f}{\wek b}{\wek a}$. Otherwise, we say that $\wek a$ and $\wek b$ are \emph{incomparable}.

\begin{wn}
There exist two incomparable $\omega_1$-\fra\ sequences in \trees.
\end{wn}

\begin{pf}
Let $U=\setof{x\in 2^{<\omega_1}}{|x^{-1}(1)|<\aleph_0}$ and let $V$ be a healthy binary Aronszajn tree. 
Clearly, $U$ and $V$ are not isomorphic.
Both trees can be naturally decomposed into $\omega_1$-sequences $\wek U$ and $\wek V$ respectively. By Theorem \ref{ofhqofjh}, $\wek U$ and $\wek V$ are \fra\ sequences. By Theorem \ref{efhfhqo}, these sequences are incomparable.
\end{pf}

It is not hard to see that there are uncountably many pairwise non-isomorphic healthy binary trees of height $\omega_1$. 
Actually, it was proved by Abraham \& Shelah~\cite{AbrShe} that consistently there are $2^{\aleph_1}$ isomorphic types of Aronszajn trees.
In particular, the category $\trees$ has (at least consistently) $2^{\aleph_1}$ many pairwise incomparable $\omega_1$-\fra\ sequences.

\subsection{Universal trees}\label{SubXtreesfinerrw}

We now discuss certain classes of trees and associated categories, showing how \fra\ functors can help in finding universal trees.
We fix two cardinals $\lam > 0$, $\kappa \goe \aleph_0$ and we assume that $\kappa$ is regular.
We also fix a linearly ordered set $L$ with $0 = \min L$ and such that
\begin{enumerate}
	\item[(0)] $0 = \min L$ exists and $L$ is unbounded from above.
	\item[(1)] Every bounded increasing sequence of length $<\kappa$ has the supremum in $L$.
	\item[(2)] $L$ contains isomorphic copies of all ordinals $< \kappa$.
\end{enumerate}
Our aim is to find an $L$-embeddable $\lam$-branching tree of height $\kappa$, universal for this class of trees.

Recall that a tree $T$ is \emph{$L$-embeddable} if it admits a $<$-preserving function $\map \phi T L$.
A standard inductive argument using condition (1) shows that such a function can be ``corrected" to a continuous one.
More precisely, there exists a continuous $<$-preserving function $\map \psi T L$ such that $\psi \loe \phi$.
Furthermore, we may require that $\psi(0) = 0$.
Recall that $\psi$ is continuous if $\psi(\sup C) = \sup \img \psi C$ whenever $C$ is a bounded chain in $T$.
\index{tree!-- $\psi$-closed}
Given a pair $\pair T\psi$, where $\psi$ is as above, we will say that $T$ is \emph{$\psi$-closed} if 
for every unbounded chain $C \subs T$ the set $\img \psi C$ is unbounded in $L$.

Define the category $\fT_{L,\kappa,\lam}$ as follows.
The objects are pairs $\pair T \psi$ such that $T$ is a $\lam$-branching tree of height $< \kappa$, $\phi$ is a continuous strictly increasing function from $T$ to $L$ such that $\psi(0)=0$ and $T$ is $\phi$-closed.
Recall that $T$ is $\lam$-branching if $|t^+| \loe \lam$ for every $t \in T$.

\index{tree!-- end-extension}
A $\fT_{L,\kappa,\lam}$-arrow from $\pair T \psi$ to $\pair S \phi$ is a tree embedding $\map f T S$ such that $\phi(f(t)) = \psi(t)$ for every $t\in T$, $\img f T$ is closed in $S$ and $S$ is an \emph{end-extension} of $\img f T$.
This means, by definition, that for every $t\in T$ and $s\in S \setminus \img f T$ it holds that $\Lev_S(f(t)) < \Lev_S(s)$.
It is clear that $\fT_{L,\kappa,\lam}$ fails the amalgamation property.
Indeed, let $T = \sn 0$ be the trivial tree and let $\map f T S$, and $\map g T S$ be such that $S$ is the 2-element tree $\dn 01$ endowed with functions $\map {\phi_0, \phi_1} S L$ such that $\phi_0(1) \ne \phi_1(1)$.
Then it is not possible to find a pair $\pair R \psi$ such that $R$ is the end-extension of both $\pair S {\phi_0}$ and $\pair S {\phi_1}$ and $\map \psi R L$ extends both $\phi_0$ and $\phi_1$.

We shall now embed $\fT_{L,\kappa,\lam}$ in a bigger category, as follows.
Namely, let $\fS_{L,\kappa,\lam}$ be the category whose objects are the same as the objects of $\fT_{L,\kappa,\lam}$, while an $\fS_{L,\kappa,\lam}$-arrow from $\pair T \psi$ to $\pair S \phi$ is a tree embedding $\map f T S$ satisfying 
\begin{enumerate}
	\item[(3)] $\phi(f(t)) \loe \psi(t)$ for every $t \in T$.
	\item[(4)] $\img f T$ is an initial segment of $S$.
\end{enumerate}
Note that the trivial tree $\sn 0$ provides an initial object for both $\fT_{L,\kappa,\lam}$ and $\fS_{L,\kappa,\lam}$.
The following fact is straightforward.

\begin{lm}\label{LmMapmmer}
The pair $\pair {\fT_{L,\kappa,\lam}} {\fS_{L,\kappa,\lam}}$ has the mixed amalgamation property.
\end{lm}

We now specify an additional property of the linearly ordered set $L$, needed for the existence of a $\kappa$-\fra\ sequence over the inclusion $\fT_{L,\kappa,\lam} \subs \fS_{L,\kappa,\lam}$.
Namely, assume $L$ satisfies the following condition:
\begin{enumerate}
	\item[($\star$)] For every $x \in L$ there exists $C_x \subs L$ such that $\inf C_x = \inf(x,\rightarrow)$ and $|C_x| \loe \lam$.
\end{enumerate}
Let us fix the sets $C_x$ as above and let
$\map p {L \times \lam}L$ be the canonical projection.
Define $U$ to be the set of all functions $\map t {\al+1} {L \times \lam}$ such that
\begin{enumerate}
	\item[(i)] $\al < \kappa$ and $p \cmp t$ is continuous.
	\item[(ii)] $p(t(\xi+1)) \in C_{p(t(\xi))}$ for every $\xi < \al$.
\end{enumerate}
The order of $U$ is extension, that is, $t \loe s$ iff $s$ extends $t$.
Note that $U$ has no maximal elements, because the domain of each $t \in U$ is a successor ordinal and $L$ is unbounded from above.
By condition (2), the height of $U$ is $\kappa$.
We shall see that the canonical decomposition of $U$ provides a \fra\ sequence over $\fT_{L,\kappa,\lam} \subs \fS_{L,\kappa,\lam}$.
More precisely, let $\map \theta U L$ be defined by $\theta(t) = \max \dom(t)$.
Then $\theta$ is a continuous strictly increasing function.
Let $U_\al = \setof{t \in U}{\theta(t) \loe \al}$ and let $\theta_\al = \theta \rest U_\al$.
Clearly, $\pair {U_\al}{\theta_\al}$ is a $\fT_{L,\kappa,\lam}$-object and $\wek u = \sett{\pair {U_\al}{\theta_\al}}{\al < \kappa}$ with the obvious embeddings is a continuous sequence in $\fT_{L,\kappa,\lam}$.

\begin{lm}\label{LmFnirwfth}
The sequence $\wek u$ is \fra\ over the inclusion $\fT_{L,\kappa,\lam} \subs \fS_{L,\kappa,\lam}$.
\end{lm}

\begin{pf}
Fix $\al < \kappa$ and fix a $\fT_{L,\kappa,\lam}$-arrow $\map f {\pair {U_\al}{\theta_\al}}{\pair S \psi}$.
We may assume that $f$ is ``simple" in the sense that $S \setminus \img f {U_\al}$ consists of a single level, since arbitrary $\fT_{L,\kappa,\lam}$-arrows are compositions of ``simple" arrows and we can later use induction.

We shall find an $\fS_{L,\kappa,\lam}$-arrow $\map g {\pair S \psi}{\pair {U_{\al+1}}{\theta_{\al+1}}}$ so that $g \cmp f$ is the inclusion $U_\al \subs U_{\al+1}$.
Namely, let $T = \img f {U_\al}$ and start with $g_0 = f^{-1}$.
Let $K$ be the highest level of $T$ and let $T_s = T \cup s^+$ for $s \in K$.
We extend $g_0$ to each $T_s$ independently.
Namely, fix $s\in K$ and enumerate $s^+$ as $\sett{s_\xi}{\xi < \lam}$ (possibly with repetitions).
Let $t = g_0(s)$.
Then $\map t {\al+1}{L \times \lam}$ satisfies conditions (i), (ii) above.
Fix $\xi < \lam$.
Using ($\star$), we can extend  $t$ to a function $t_\xi \in U$ whose domain is $\al+2$ and $p(t_\xi(\al+1)) \loe \phi(s_\xi)$.
Define $g_s(s_\xi) = t_\xi$ and $g_s \rest T = g_0$.
This defines $\map {g_s}{T_s}{U_{\al+1}}$ extending $g_0$.
The union of all these maps defines the required $\fS_{L,\kappa,\lam}$-arrow $g$, showing that $\wek u$ is \fra\ over $\fT_{L,\kappa,\lam} \subs \fS_{L,\kappa,\lam}$.
\end{pf}

Using Lemmata~\ref{LmMapmmer}, \ref{LmFnirwfth} and Theorem~\ref{TmixedamsCof}, we obtain

\begin{tw}
Under the assumptions above, there exists a $\lam$-branching tree $U$ of height $\kappa$ together with a continuous strictly increasing function $\map \theta U L$, such that for every $\lam$-branching tree $T$ of height $\loe \kappa$, for every strictly increasing function $\map \phi T L$, there exists a tree embedding $\map f T U$ such that $\img f T$ is an initial segment of $U$ and $\theta(f(t)) \loe \phi(t)$ for every $t \in T$.
\end{tw}

Of course, the above result can be proved directly, without referring to \fra\ sequences.
The key inductive argument is contained in the proof of Lemma~\ref{LmFnirwfth}.
Our aim was rather showing how the existence and universal properties of certain trees can be explained by means of \fra\ sequences over a natural inclusion functor.

One of the special cases is $L = \kappa$; then $\lam$ could be even finite and ``being $L$-embeddable" is the same as ``being of height $\loe \kappa$".
When $L$ has the property that $|(x,\rightarrow)| = \lam$ for every $x \in L$, a natural example of a universal $\lam$-branching $L$-embeddable tree is the tree of all bounded closed well ordered subsets of $L$ endowed with the ``end-extension" ordering.
It is worth mentioning that the well-known tree $\sig\Qyu$ of bounded well ordered subsets of the rational numbers (again with the ``end-extension" ordering) is another example of a universal countably branching $\Err$-embeddable tree (see \cite{StevoRtrees} for a direct argument and more information).
It is a folklore fact that $\sig\Qyu$ is not $\Qyu$-embeddable, however its subtree consisting of all closed sets is, by the arguments above, a universal countably branching $\Qyu$-embeddable tree.

Finally, one can consider the category $\fK$ of finite trees with end-extensions and the category $\fL$ of finite trees with arbitrary embeddings.
In this case, sequences in $\fK$ lead to finitely-branching trees of height $\loe\omega$ and a sequence is \fra\ over the inclusion $\fK \subs \fL$ if and only if it leads to a perfect tree.
Recall that a tree $T$ is \emph{perfect} if for every $t\in T$ there are incomparable elements $s_0, s_1$ such that $s_0 > t$ and $s_1 > t$.
\index{tree!-- perfect}
As before, the category $\fK$ obviously fails the amalgamation property.
On the other hand, $\fL$ has the amalgamation property and since it is countable and directed, it has a unique \fra\ sequence.
This sequence, however, does not lead to a tree: Its natural co-limit is a certain tree-like partially ordered set, in the sense that all intervals $(a,b)$ with $a < b$ are linearly ordered (in fact, isomorphic to $\Qyu$) and for every $a < b$ there is $c > a$ not comparable with $b$.

\section{Embedding-projection pairs}\label{reterpeteairs}

In this section we describe a general construction on a given category, which is suitable for applications to the theory of Valdivia compacta and Banach spaces with Markushevich bases.
This construction had been first used by Dana Scott in order to get faithful models of untyped $\lambda$-calculus. It also appears in Droste \& G\"obel \cite{DrGoe92} in the context of Scott domains.

\subsection{Definitions and basic properties}

We fix a category $\fK$. Define $\rp\fK$ to be the category whose objects are the objects of $\fK$ and a morphism $\map fXY$ is a pair $\pair er$ of arrows in $\fK$ such that $\map eXY$, $\map rYX$ and $r\cmp e=\id X$.
We set $e(f):=e$ and $r(f):=r$, so $f=\pair{e(f)}{r(f)}$. Given morphisms $\map fXY$ and $\map gYZ$ in $\rp\fK$, we define their composition in the obvious way:
$$g\cmp f := \pair{e(g)\cmp e(f)}{r(f)\cmp r(g)}.$$
It is clear that this defines an associative operation on compatible arrows. Further, given an object $a\in \fK$, pair of the form $\pair{\id a}{\id a}$ is the identity morphism in $\rp\fK$. Thus, $\rp\fK$ is indeed a category. Note that $f\mapsto e(f)$ defines a covariant functor $\map e {\rp\fK}\fK$ and $f\mapsto r(f)$ defines a contravariant functor $\map r{\rp\fK}{\fK}$.
\index{EP-pairs}\index{embedding-projection pairs}
Following \cite{DrGoe92}, we shall call $\rp\fK$ the \emph{category of embedding-projection pairs} or briefly \emph{EP-pairs}.
The idea of considering EP-pairs is to obtain more special \fra\ sequences, namely we would like to obtain a unique cofinal object $u$ in the category of sequences such that every other $x$ can be both ``embedded" into $u$ and ``projected" from $u$. In other words, we would like to have arrows $\map jxu$ and $\map rux$ satisfying $r\cmp j = \id{x}$.
We shall demonstrate this in Subsection~\ref{SubSectCantorStDns} below, dealing with the category of finite nonempty sets.

We shall later need the following general facts about sequences of left-invertible arrows.

\begin{lm}\label{oiutiq}
Let $\fK$ be a category, let $\delta>0$ be a limit ordinal and let $\map a\delta{\rp\fK}$ be a sequence such that $\pair v{\sett{i_\al}{\al<\delta}}$ is the co-limit of $\img e{a}$. Then there exists a sequence $\map b{(\delta+1)}{\rp\fK}$ such that $b\rest\delta=a$, $b_\delta=v$ and $e(b_\al^\delta)=i_\al$ for every $\al<\delta$.
\end{lm}

\begin{pf}
Fix $\al<\delta$. Given $\al<\xi<\eta<\delta$, we have
$$r(a^\eta_\al) \cmp e(a_\xi^\eta) = r(a^\xi_\al) \cmp r(a^\eta_\xi) \cmp e(a_\xi^\eta) = r(a^\xi_\al).$$
Using the equation above and the fact that $v$ together with $\sett{i_\xi}{\al\loe\xi<\delta}$ is the co-limit of $a\rest [\al,\delta)$, we find a unique arrow $\map {f_\al}{v}{a_\al}$ satisfying $$f_\al \cmp i_\xi = r(a^\xi_\al)$$ for every $\al\loe\xi<\delta$. 
Setting $\xi:=\al$ in the equation above, we see that $f_\al \cmp i_\al = \id{a_\al}$.
By uniqueness, $f_\xi =  r(a_\xi^\eta) \cmp f_\eta$ for every $\xi<\eta<\delta$.
Finally, define $b_\delta := v$, $b^\delta_\al:=\pair {i_\al}{f_\al}$ and $b\rest \delta := a$.
\end{pf}

\begin{lm}\label{kwerjpit}
Let $\wek x$ be a continuous sequence in a category $\fK$ and assume that each bonding arrow $x_\al^\beta$ is left-invertible in $\fK$. Then there exists a sequence $\wek y$ in $\rp\fK$ such that $$\wek x= \img e{\wek y}.$$
\end{lm}

\begin{pf}
Let $\map{\wek x}\delta\fK$. Using induction, we shall construct arrows $\map {f_\al^\beta}{x_\beta}{x_\al}$ in $\fK$ so that
\begin{equation}
f_\al^\beta \cmp x_\al^\beta = \id{x_\al} \oraz
f_\al^\gamma = f_\beta^\gamma \cmp f_\al^\beta
\tag{*}\label{wetj0qwerj}
\end{equation}
holds for every $\al<\beta<\gamma$.
Fix $0<\beta<\delta$ and suppose $f_\xi^\eta$ have been constructed for every $\xi<\eta<\beta$.
Assume first that $\beta$ is a successor ordinal, say $\beta=\al+1$. Using the assumption, find $\map f{x_{\al+1}}{x_\al}$ such that $f\cmp x_\al^{\al+1}=\id{x_\al}$. Given $\xi<\beta$, define $f^{\al+1}_\xi = f^\al_\xi \cmp f$. It is clear that (\ref{wetj0qwerj}) holds.

Assume now that $\beta$ is a limit ordinal. Define $\map a\beta{\rp\fK}$ by setting $a_\al = x_\al$ and $a_\xi^\eta = \pair{x_\xi^\eta}{f_\xi^\eta}$ for $\xi<\eta<\beta$. By (\ref{wetj0qwerj}), $a$ is indeed a sequence in $\rp\fK$. By Lemma \ref{oiutiq}, this sequence has an extension $\map b{\beta+1}{\rp\fK}$ such that $e(b_\al^\beta)=x_\al^\beta$ for every $\al<\beta$. Setting $f_\al^\beta:=r(b_\al^\beta)$ for $\al<\beta$, it is clear that (\ref{wetj0qwerj}) holds.

Finally, setting $y_\al := x_\al$ and $y_\al^\beta := \pair{x_\al^\beta}{f_\al^\beta}$, we obtain the desired sequence $\wek y$ in the category $\rp\fK$.
\end{pf}

\subsection{Semi-continuity and the back-and-forth principle}

It turns out that usually $\rp\fK$ is not continuous, however one can consider the following weakening of continuity, which is good enough for applications.

\index{semi-continuous sequence}\index{sequence!-- semi-continuous}
We say that a sequence $\wek x$ in $\rp\fK$ is \emph{semi-continuous} if $\img e{\wek x}$ is continuous in $\fK$. The dual notion of semi-continuity with respect to the functor $r$ is obtained by considering $\rp{\dualcat\fK}$ instead of $\rp\fK$.

\begin{tw}\label{ierupqu}
Let $\kappa$ be an infinite regular cardinal and let $\fK$ be a $\kappa$-complete category such that $\rp\fK$ is directed and has the amalgamation property. Assume further that $\rp\fK$ has a dominating family consisting of at most $\kappa$ arrows. Then there exists a semi-continuous \fra\ sequence of length $\kappa$ in $\rp\fK$.
\end{tw}

\begin{pf}
By Lemma \ref{oiutiq}, $\rp\fK$ is $\kappa$-bounded. Thus, we may use Proposition~\ref{strangewef} applied to $\Phi:=e$.
\end{pf}

We shall prove that semi-continuous \fra\ sequences satisfy the back-and-forth principle. In order to do it, we need to show that the left-forgetful functor $e$ induces a quasi-limiting operator, as it is suggested by Lemma~\ref{oiutiq}.
It turns out that this is not completely trivial.
Note that we cannot use Proposition~\ref{Pkwasilimitz}, because the functor $e$ is not faithful.

\begin{lm}\label{LemQwasiLimOprtRPscot}
Let $\fK$ be a category.
Given a sequence $\wek x$ in $\rp \fK$ such that $\img e{\wek x}$ has the co-limit in $\fK$, let $\qlim \wek x$ be the co-cone in $\rp \fK$, existing by Lemma~\ref{oiutiq}, such that $e({\qlim \wek x}) = \lim \img e{\wek x}$.
Then $\qlim$ is a partial quasi-limiting operator on $\uzup \kappa \fK$.
\end{lm}

\begin{pf}
Let $\map {\wek f}{\wek x}{\wek y}$ be an isomorphism of sequences in $\rp \fK$ and let $\map {\wek g}{\wek y}{\wek x}$ be its inverse.
We may assume that $\wek f = \sett{f_\xi}{\xi < \lam}$, $\wek g = \sett{g_\xi}{\xi < \lam}$ so that $\map {f_\xi}{x_\xi}{y_\xi}$ and $\map {g_\xi}{y_\xi}{x_{\xi + 1}}$ for $\xi < \lam$.
Let the $\fK$-objects $\ovr x, \ovr y$ be the co-limits of $\img e {\wek x}$ and $\img e {\wek y}$, respectively.
The co-limiting co-cones will be denoted by $\sett{x_\xi^\infty}{\xi < \lam}$ and $\sett{y_\xi^\infty}{\xi < \lam}$, respectively.
Let $\map {f_\infty}{\ovr x}{\ovr y}$ be the unique arrow satisfying $e(y_\xi^\infty) \cmp e(f_\xi) = f_\infty \cmp e(x_\xi^\infty)$ for every $\xi < \lam$.
Similarly, let $\map {g_\infty}{\ovr y}{\ovr x}$ be the unique arrow satisfying $e(x_{\xi+1}^\infty) \cmp e(g_\xi) = g_\infty \cmp e(y_\xi^\infty)$ for every $\xi < \lam$.
By the universality of the co-limit, $f_\infty \cmp g_\infty = \id{\ovr y}$ and $g_\infty \cmp f_\infty = \id{\ovr x}$.
In particular, $\pair {f_\infty}{g_\infty}$ is an isomorphism in $\rp \fK$.
It remains to check that it commutes with $\wek f$ and $\wek g$.
By symmetry, it suffices to check that
\begin{equation}
r(f_\al) \cmp r(y_\al^\infty) = r(x_\al^\infty) \cmp g_\infty
\tag{1}\label{Eqneorgoerg}
\end{equation}
holds for all $\al < \lam$.

Fix $\al < \lam$.
Given $\xi \goe \al$, define $q_\xi = r(f_\al) \cmp r(y_\al^\xi)$.
If $\xi < \eta$ then
$$q_\eta \cmp e(y_\xi^\eta) = r(f_\al) \cmp r(y_\al^\xi) \cmp r(y_\xi^\eta) \cmp e(y_\xi^\eta) = q_\xi.$$
It follows that $\sett{q_\xi}{\al \loe \xi < \lam}$ commutes with the sequence $\sett{e(y_\xi)}{\al \loe \xi < \lam}$.
As $g_\infty$ is the co-limit of this sequence, there is a unique $\fK$-arrow $\map k{\ovr y}{x_\al}$ satisfying
\begin{equation}
k \cmp e(y_\xi^\infty) = q_\xi = r(f_\al) \cmp r(y_\al^\xi)
\tag{2}\label{Eqebtoert}
\end{equation}
for every $\xi \goe \al$.
Let
$$
k_0 = r(f_\al) \cmp r(y_\al^\infty) \oraz k_1 = r(x_\al^\infty) \cmp g_\infty.
$$
In order to show (\ref{Eqneorgoerg}), we check that both $k_0$ and $k_1$ satisfy (\ref{Eqebtoert}) in place of $k$.
Given $\xi \goe \al$, we have
$$
k_0 \cmp e(y_\xi^\infty) = r(f_\al) \cmp r(y_\al^\xi) \cmp r(y_\xi^\infty) \cmp e(y_\xi^\infty) = r(f_\al) \cmp r(y_\al^\xi).
$$
Hence $k_0 = k$.
Recalling that $g_\infty$ is the unique $\fK$-arrow satisfying $e(x_{\xi+1}^\infty) \cmp e(g_\xi) = g_\infty \cmp e(y_\xi^\infty)$ for every $\xi < \lam$, we obtain
\begin{align*}
k_1 \cmp e(y_\xi^\infty) &= r(x_\al^\infty) \cmp g_\infty \cmp e(y_\xi^\infty) = r(x_\al^\infty) \cmp e(x_{\xi+1}^\infty) \cmp e(g_\xi) \\
&= r(x_\al^{\xi+1}) \cmp r(x_{\xi+1}^\infty) \cmp e(x_{\xi+1}^\infty) \cmp e(g_\xi) = r(x_\al^{\xi+1}) \cmp e(g_\xi) \\
&= r(g_\xi \cmp y_\al^\xi \cmp f_\al) \cmp e(g_\xi) = r(f_\al) \cmp r(y_\al^\xi) \cmp r(g_\xi) \cmp e(g_\xi) \\
&=  r(f_\al) \cmp r(y_\al^\xi).
\end{align*}
Hence $k_1 = k$, which completes the proof.
\end{pf}

Combining the lemma above with Theorem~\ref{Tfjogosdg}, we obtain:

\begin{tw}\label{Thmhot55}
Let $\fK$ be a category.
Every two semi-continuous \fra\ sequences of the same regular length in $\rp\fK$ satisfy the back-and-forth principle.
\end{tw}

\subsection{Proper amalgamations and proper arrows}

We still have not proved that a semi-continuous \fra\ sequence in a category of EP-pairs is cofinal for uncountable sequences. In fact, even for the countable case, one can consider a better class of arrows between sequences, obtaining stronger cofinality results (see Subsection~\ref{SubSectCantorStDns} for a simple application in the category of sets). This leads to the following concept.

\index{proper amalgamation}\index{amalgamation!-- proper}
Let $\map fZX$ and $\map gZY$ be arrows in $\rp\fK$. We say that arrows $\map hXW$, $\map kYW$ provide a \emph{proper amalgamation}
of $f,g$ if $h\cmp f=k\cmp g$ and moreover $e(g)\cmp r(f)=r(k)\cmp e(h)$, $e(f)\cmp r(g)=r(h)\cmp e(k)$ hold.
Translating it back to the original category $\fK$, this means that the following four diagrams commute:
$$
\xymatrix{
W & Y\ar@{ >->}[l]_{e(k)}\\
X\ar@{ >->}[u]^{e(h)} & Z\ar@{ >->}[l]_{e(f)}\ar@{ >->}[u]_{e(g)}
}
\quad
\xymatrix{
W\ar@{->>}[r]^{r(k)}\ar@{->>}[d]_{r(h)} & Y\ar@{->>}[d]^{r(g)}\\
X\ar@{->>}[r]^{r(f)} & Z
}
\quad
\xymatrix{
W\ar@{->>}[r]^{r(k)} & Y\\
X\ar@{ >->}[u]^{e(h)}\ar@{->>}[r]^{r(f)} & Z\ar@{ >->}[u]_{e(g)}
}
\quad
\xymatrix{
W\ar@{->>}[d]_{r(h)} & Y\ar@{->>}[d]^{r(g)}\ar@{ >->}[l]_{e(k)}\\
X & Z\ar@{ >->}[l]_{e(f)}
}
$$
We draw arrows $\xymatrix{\ar@{ >->}[r] &}$ and $\xymatrix{\ar@{->>}[r] &}$ in order to indicate mono- and epimorphisms respectively.
We shall say that $\rp\fK$ \emph{has proper amalgamations} if every pair of arrows in $\rp\fK$ with common domain can be properly amalgamated in $\rp\fK$.

Below is a useful criterion for the existence of proper amalgamations.

\begin{lm}\label{aofihapsfj}
Let $\fK$ be a category and let $f,g$ be arrows in $\rp\fK$ with the same domain. If $e(f),e(g)$ have a pushout in $\fK$ then $f,g$ can be properly amalgamated in $\rp\fK$.
\end{lm}

\begin{pf} Let $\map hXW$ and $\map kYW$ form a pushout of $e(f)$, $e(g)$.
Consider the following diagram:
$$
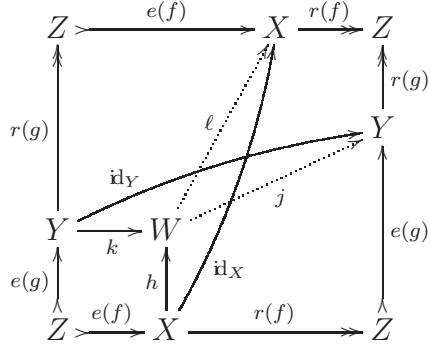
\xymatrix{
Z\ar@{ >->}[rr]^{e(f)} & & X\ar@{->>}[r]^{r(f)} & Z \\
 & & & Y\ar@{->>}[u]_{r(g)} \\
Y\ar@{->>}[uu]^{r(g)}\ar[r]_k\ar@/^/[urrr]^(.25){\id Y} & W\ar@{.>}[uur]^\ell \ar@{.>}[urr]_j & & \\
Z\ar@{ >->}[u]^{e(g)}\ar@{ >->}[r]^{e(f)} & X\ar[u]^h\ar@{->>}[rr]^{r(f)}\ar@/_/[uuur]_(.25){\id X} & & Z\ar@{ >->}[uu]_{e(g)}
}$$
The dotted arrows indicate unique morphisms completing appropriate diagrams, i.e. $j$ is the unique arrow satisfying equations $j\cmp h=e(g)\cmp r(f)$, $j\cmp k=\id Y$ and $\ell$ is the unique arrow satisfying equations $\ell\cmp k= e(f)\cmp r(g)$, $\ell\cmp h = \id X$. Consequently, $\pair kj$ and $\pair h\ell$ are morphisms in $\rp\fK$.
Set $s= r(f)\cmp \ell$. Then
\begin{equation}
s\cmp k = r(f)\cmp \ell\cmp k = r(f)\cmp e(f)\cmp r(g) = r(g) \quad\text{ and }\quad s\cmp h = r(f)\cmp \ell\cmp h  = r(f).
\tag{1}
\end{equation}
Recall that $r(f)\cmp e(f)=\id Z = r(g)\cmp e(g)$. Since $k,h$ is a pushout of $e(f)$, $e(g)$, we deduce that $s$ must be the unique arrow satisfying (1). Now let $t=r(g)\cmp j$. Similar computations show that $t\cmp k=r(g)$ and $t\cmp h=r(f)$, therefore by uniqueness we deduce that $s=t$ or, in other words, $r(f)\cmp \ell  = r(g)\cmp j$. This shows that the full diagram is commutative and hence $\pair kj$ and $\pair h\ell$ provide a proper amalgamation of $f,g$ in the category $\rp\fK$.
\end{pf}

As an example, if $\fK$ is the category of nonempty sets, then Lemma \ref{aofihapsfj} says that category $\rp\fK$ has proper amalgamations. We show below that not all amalgamations in $\rp\fK$ are proper.

\begin{ex}
Consider the category of nonempty finite sets $\sets^+$.
Let $a,b,c,d$ be pairwise distinct elements and set $Z=\sn a$, $X=\dn ab$, $Y=\dn ac$ and $W=\set{a,b,c}$.
We are going to define arrows $\map fZX$, $\map gZY$, $\map hXW$ and $\map kYW$ in the category $\rp\sets^+$.
Let $e(f)$, $e(g)$, $e(h)$ and $e(k)$ be the inclusion maps and let $r(f)$ and $r(g)$ be the obvious constant maps. Finally, let $r(h)(c)=a$ and $r(k)(b)=c$. 
This already defines $r(h)$ and $r(k)$, since these maps must be identity on the ranges of $e(h)$ and $e(k)$ respectively.
It is clear that $h\cmp f=k\cmp g$, i.e. $h,k$ amalgamate $f,g$ in the category $\rp\sets^+$.
On the other hand, $e(g)\cmp r(f)(b)=a$ and $r(k)\cmp e(h)(b)=c$, therefore $e(g)\cmp r(f)\ne r(k)\cmp e(h)$.
Note that actually $e(f)\cmp r(g)=r(h)\cmp e(k)$ holds, although redefining $r(h)(c)$ to $b$ we can even get $e(f)\cmp r(g)\ne r(h)\cmp e(k)$.
\end{ex}

\begin{lm}\label{Lemfjprprstb}
Let $\fK$ be a category such that sequences of length $< \kappa$ consisting of left-invertible arrows have co-limits, and  assume that $\rp \fK$ has proper amalgamations. Let $\qlim$ be the quasi-limiting operator on $\uzup \kappa {\rp \fK}$ associated with the left-forgetful functor $e$.
Let $\Ama$ denote the class of all proper amalgamations.
Then $\Ama$ is a stable amalgamation structure in $\rp \fK$ and $\triple {\rp\fK} \qlim \Ama$ is a \frajon\ category.
\end{lm}

\begin{pf}
The fact that $\Ama$ is an $\qlim$-continuous amalgamation structure is obvious (assuming of course that $\rp \fK$ admits proper amalgamations).
We only need to check that $\Ama$ is stable.

Fix a commutative diagram of the form
$$\xymatrix{
b \ar[r]^g & c \\
a \ar[r]_{\id a} \ar[u]^f & a \ar[u]_h
}$$
where $f,g,h$ are $\rp\fK$-arrows.
However, using the fact that $r(h) = r(f) \cmp r(g)$, we get $r(h) \cmp e(g) = r(f) \cmp r(g) \cmp e(g) = r(f)$ and, using the fact that $e(h) = e(g) \cmp e(f)$, we obtain
$r(g) \cmp e(h) = r(g) \cmp (g) \cmp e(f) = e(f)$.
\end{pf}

\index{proper arrow}
An arrow of sequences $\map{\wek f}{\wek x}{\wek y}$ that is admissible for the class of proper amalgamations will be called \emph{proper}.
In other words, $\wek f$ is proper if (up to equivalence) all squares of the form
$$\xymatrix{
x_\al \ar[rr]^{x_\al^\beta} \ar[d]_{f_\al} & & x_\beta \ar[d]^{f_\beta} \\
y_{\phi(\al)} \ar[rr]^{y_{\phi(\al)}^{\phi(\beta)}} & & y_{\phi(\beta)}
}$$
are proper amalgamations, where $\phi$ is the increasing ordinal function associated to $\wek f$.
As we shall see in Example~\ref{Exeimprprs}, proper arrows are rather restrictive even in the category of countable sequences of nonempty finite sets.

Using Lemma~\ref{Lemfjprprstb} and Theorem~\ref{Thmeotbnoerge}, we obtain:

\begin{tw}\label{lprpraf}
Assume $\fK$ is a category and $\wek u$ is a semi-continuous \fra\ sequence in $\rp \fK$ of regular length $\kappa\goe\aleph_0$. If $\rp\fK$ has proper amalgamations then for every semi-continuous sequence $\wek x\in\uzuple{\kappa}{\rp\fK}$ there exists a proper arrow of sequences $\map{\wek f}{\wek x}{\wek u}$.
\end{tw}

In Subsection~\ref{SubSectCantorStDns} we shall present a simple application of the result above (with $\kappa=\aleph_0$) in the category of nonempty finite sets, leading to a classical homogeneity property of the Cantor set.

\section{Applications}\label{SectAppsEPetsoon}

In this section we collect several applications of our results---mainly those from Section \ref{reterpeteairs}---to compact Hausdorff spaces, Banach spaces, the Cantor set, and linearly ordered sets.

\subsection{Compact spaces}

Let $\komp$ be the category of all nonempty compact Hausdorff spaces.
It is well-known that quotient maps have pullbacks in $\komp$.
More precisely, given quotient maps $\map f X Z$, $\map g Y Z$, the two projections from the set
$$P = \setof{\pair x y \in X \times Y}{f(x) = g(y)}$$
from the pullback of $f,g$.
It is obvious that $P$ is compact, as a closed subset of $X \times Y$ and both projections restricted to $P$ are onto.
Notice that, assuming the Continuum Hypothesis, the category of nonempty compact metric spaces with quotient maps has $\aleph_1$ many types of arrows, there it has an inverse \fra\ sequence of length $\omega_1$.
The limit of this sequence is the well-known space $\Nat^* = \beta \Nat \setminus \Nat$, i.e., the remainder of the \v Cech-Stone compactification of the set of natural numbers.
All these things clearly generalize to higher cardinals and were studied by Negrepontis~\cite{Negrep}, using \jon's work~\cite{Jon}.
In fact, when restricted to totally disconnected compact spaces, Stone duality allows moving to the category of Boolean algebras, where model-theoretic \frajon\ theory can be applied.
The fact that every compact Hausdorff space of weight $\loe\aleph_1$ is a continuous image of $\Nat^*$ was previously proved by Parovi\v cenko~\cite{Parov}.

\index{Valdivia compact}
\index{Corson compact}
\index{$\Sigma$-product}\index{support}
We now turn to a more special class of spaces, called Valdivia compacta.
A topological space $K$ is called \emph{Valdivia compact} if it is homeomorphic to a closed subspace $K'$ of some Tikhonov cube $\unii^\kappa$, such that $K'\cap \Sigma(\kappa)$ is dense in $K'$, where
$$\Sigma(\kappa) = \setof{x \in \unii^\kappa}{|\suppt(x)| \loe \aleph_0}$$
is the \emph{$\Sigma$-product} of $\kappa$ copies of the unit interval, and $\suppt(x) = \setof{\al}{x(\al)\ne 0}$ is the \emph{support} of $x$.
This class has been extensively studied by several authors; we refer to Kalenda's surveys~\cite{Kal2, Kal1} for further references.
Valdivia compacta are closely related to \emph{Corson compacta} which are defined by strengthening the requirement of $K'\cap \Sigma(\kappa)$ being dense to $K' \subs \Sigma(\kappa)$.
It is well-known that Corson compacta are precisely those Valdivia compacta that are countably tight.
This actually follows easily from the (not completely trivial) fact that $\Sigma$-products are countably tight.
\index{countably tight space}
Recall that a topological space is \emph{countably tight} if the closure of any set is the union of closures of its countable subsets.
One has to mention yet another class of compacta: Eberlein compact spaces. By definition,  a space is \emph{Eberlein compact} if it is homeomorphic to a weakly compact subset of some Banach space.
If this Banach space is Hilbert, the compact is called \emph{uniformly Eberlein}.
Every Eberlein compact is Corson and every Corson compact is Valdivia.
None of these implications can be reversed.

A result of Bell~\cite{Bell} says that, under the Continuum Hypothesis, there exists a universal uniform Eberlein compact of weight $\cont$.
Another positive result, due to Bell \& Marciszewski~\cite{BellMarciszewski} shows, under the usual assumption $2^{<\kappa}=\kappa$, the existence of a scattered Eberlein compact of weight $\kappa$ and of Cantor-Bendixson height 3 that is universal for this class in the sense of retractions.
Both of the results on Eberlein compacta are in fact \frajon\ constructions, by putting relevant structures on these spaces.

In contrast, Argyros \& Benyamini~\cite{ArgyrosBenyamini} proved that universal Eberlein compacta for weight $\kappa$ do not exist whenever $\kappa = \aleph_1$ or $\kappa^{\aleph_0}=\kappa$.
Later, Todor\v cevi\'c~\cite{To95} proved the same negative result for Corson compacta.
Both results are negative even when asking for a compact space $K$ such that every other compact in the given class is homeomorphic to a quotient of a closed subset of $K$.

The results of Section~\ref{reterpeteairs} allow us to get, under the Continuum Hypothesis, a Valdivia compact of weight $\aleph_1$ that is universal for this class in the sense of retractions.
More precisely:

\begin{tw}\label{ThmUnVldxyone}
Assume $2^{\aleph_0} = \aleph_1$.
There exists a compact space $V$ with the following properties.
\begin{enumerate}
	\item[$(1)$] $\w(V) = \aleph_1$ and $V$ is Valdivia compact.
	\item[$(2)$] Every Valdivia compact of weight $\loe \aleph_1$ is a retract of $V$.
	\item[$(3)$] Given a compact metrizable space $K$, given retractions $\map {f,g}V K$, there exists a homeomorphism $\map h V V$ such that $g = f \cmp h$.
\end{enumerate}
Furthermore, properties (1) and (3) describe the space $V$ uniquely, up to a homeomorphism, and every compact space satisfying (3) has weight at least $2^{\aleph_0}$.
\end{tw}

Before starting the proof, we need one important property of Valdivia compacta, explaining why we actually restrict attention to weight $\aleph_1$:

\begin{prop}[Kubi\'s \& Michalewski~\cite{KM}]\label{ThmKMhmm}
A nonempty compact space of weight $\loe \aleph_1$ is Valdivia compact if and only if it is homeomorphic to the limit of a continuous inverse sequence of compact metrizable spaces in which all bonding maps are retractions.
\end{prop}

\begin{pf}[Proof of Theorem~\ref{ThmUnVldxyone}.]
As we have mentioned above, quotient maps admit pullbacks in the category of nonempty compact spaces.
The same of course holds when restricted to the category $\fK$ of compact metric spaces.
It follows that $\rp \fK$ has proper amalgamations.
The assumption $2^{\aleph_0} = \aleph_1$ implies that $\rp \fK$ has a semi-continuous \fra\ sequence $\wek u$ of length $\omega_1$, where semi-continuity is with respect to quotient maps (obviously, every countable inverse sequence of quotient maps has a limit).
Thus, the sequence $\wek u$ leads to a continuous inverse sequence of compact metric spaces $U_\al$ ($\al < \omega_1$) in which all bonding maps are retractions.
Denote by $V$ the limit of this sequence.
By Proposition~\ref{ThmKMhmm}, $V$ is a Valdivia compact of weight $\aleph_1$, that is, $V$ satisfies (1).
Fix another nonempty Valdivia compact $K$ whose weight is $\loe \aleph_1$.
Again by Proposition~\ref{ThmKMhmm}, $K$ is the limit of a continuous sequence of retractions between metrizable compacta.
By Lemma~\ref{kwerjpit}, this sequence extends to a semi-continuous sequence $\wek x$ in $\rp \fK$.
Now Theorem~\ref{lprpraf} gives a proper arrow $\map {\wek f} {\wek x} {\wek u}$ which gives rise to a retraction $\map r V K$.
This shows (2).

Property (3) is just a translation of the back-and-forth principle, satisfied by the \fra\ sequence $\wek u$ by Theorem~\ref{Thmhot55}.
The relevant fact needed here is that every map $\map f V K$ onto a metrizable compact space factorizes through some $U_\al$, that is, it satisfies $f = f' \cmp p_\al$, where $\map {p_\al} V {U_\al}$ is the retraction obtained from the $\rp \fK$-arrow $\map {u^\infty_\al}{\wek u}{u_\al}$.

In order to show the uniqueness of $V$, let us consider the category $\fL$ of all retractions between nonempty compact metrizable spaces.
By Proposition~\ref{ThmKMhmm}, the space $V$ is the limit of a continuous inverse sequence $\wek v$ in $\fL$.
Property (3) implies that the sequence is \fra\ in $\fL$.
Even though countable inverse sequences of retractions may have no limits in $\fL$, the limits ``computed" in the bigger category of all continuous maps provide a quasi-limiting operator.
Thus, by Theorem~\ref{Tfjogosdg}, the sequence $\wek v$ is unique, up to an isomorphism of sequences.
This clearly translates to the fact that $V$ is unique up to a homeomorphism.

Let us finally see that $\w(V) \goe 2^{\aleph_0}$, whenever $V$ satisfies (3).
Clearly, there is a retraction $\map f V \unii$, where $\unii$ is the unit interval.
Fix a point $\infty \notin \unii$ and let $K = \unii \cup \sn \infty$, where $\infty$ is an isolated point.
Let $\map g V K$ be a fixed retraction.
For each $t \in \unii$, let $\map {r_t}{K}{\unii}$ be such that $r_t \rest \unii = \id{\unii}$ and $r_t(\infty) = t$.
Let $\map {h_t} V V$ be a homeomorphism such that $f = r_t \cmp g \cmp h_t$.
Let $U_t = f^{-1}(t)$.
Then $U_t = h_t^{-1}(g^{-1}(r_t^{-1}(t))) \sups h_t^{-1}(g_t^{-1}(\infty))$ and, since $\infty$ is isolated in $K$, we deduce that $U_t$ has nonempty interior.
Obviously, $U_t \cap U_s = \emptyset$ whenever $t \ne s$, which shows that $\w(V) \goe |\unii|=2^{\aleph_0}$.
\end{pf}

Using the full strength of categories of projection-embedding pairs, it is possible to make the result above more precise, by considering a pair of the form $\pair V D$, where $D = \Sigma(\omega_1) \cap V$, with respect to the appropriate embedding of $V$ into the Tikhonov cube $\unii^{\omega_1}$.
In particular, for every retractions $\map f V K$, $\map f V L$, where $K,L$ are metrizable compacta contained in $D$, every homeomorphism $\map h K L$ extends to an auto-homeomorphism of $V$.
Another fact would be that every Valdivia compact $K \subs \unii^{\omega_1}$ with a dense set $G = K \cap \unii^{\omega_1}$ is homeomorphic to a retract $K' \subs V$ of $V$ so that $G = D \cap K'$.
We have decided to present the simplified version which is more transparent and less technical.

A result from \cite{KM} says that every retract of a Valdivia compact is again Valdivia provided its weight does not exceed $\aleph_1$.
It is an open problem whether the same holds for arbitrary weight.
A category-theoretic description of the whole class of Valdivia compacta is given in~\cite{KM}. Unfortunately, it is not clear how to use it when the weight exceeds $\aleph_1$ and, in particular, it is not clear how to construct a retractively universal Valdivia compact space for larger weights.

\subsection{Banach spaces}

It is well-known that there exist separable Banach spaces that are isometrically universal for the class of all separable Banach spaces.
Perhaps the first example is $\C([0,1])$, the space of continuous functions on the unit interval, endowed with the maximum norm.
This space, however, does not have good homogeneity properties (see Proposition~\ref{progoelts} for the precise statement).
There is another separable Banach space $\G$, constructed by Gurari\u\i~\cite{Gurarii}, that is also isometrically universal for separable spaces and moreover has the following homogeneity property:
\begin{enumerate}
\item[(H)] Given an isometry $\map f X Y$ between finite-dimensional linear subspaces of $\G$, given $\eps > 0$, there exists a bijective isometry $\map h \G \G$ such that $\norm {f - h \rest X} \loe \eps$.
\end{enumerate}
Furthermore, $\G$ is the unique, up to isometry, separable Banach space containing isometric copies of all finite-dimensional spaces and satisfying (G).
For a direct and elementary proof we refer to~\cite{KubSol}.
Typically, the condition specifying the Gurari\u\i\ space is different from (G), it involves almost isometric embeddings of finite-dimensional spaces.
The space $\G$ can be viewed as some sort of an ``approximate \fra\ limit".
On the other hand, there seems to be no natural category whose \fra\ sequence would lead to $\G$.
In fact, Gurari\u\i\ had already observed that $\G$ does not satisfy the variant of condition (G) with $X$, $Y$ being one-dimensional and $\eps = 0$.

Below we show the existence, under the Continuum Hypothesis, of a natural variant of Gurari\u\i's space that comes from a \fra\ sequence in the category of separable Banach spaces.
Its existence was not known before\footnote{See footnote~\ref{Ffootone} in the last paragraph of Introduction.}, although it actually could be derived from the results of Droste \& G\"obel \cite{DrGoe92}, since the \fra\ sequence in this category can always be made continuous.

\index{\sbaniso}\index{category!-- \sbaniso}
Let $\sbaniso$ denote the category whose objects are separable Banach spaces and arrows are linear isometries.
The following fact is well-known.

\begin{lm}\label{dsgsdgerr}
$\sbaniso$ has the amalgamation property.
\end{lm}

\begin{pf}
Fix $X,Y,Z\in\sbaniso$ and fix linear isometric embeddings $\map fZX$ and $\map gZY$. Without loss of generality, we may assume that $f$ and $g$ are inclusions, i.e. $Z\subs X$ and $Z\subs Y$.
We may also assume that $X\cap Y=Z$. Now let $W$ be the formal algebraic sum of $X$ and $Y$, i.e. $W=\setof{x+y}{x\in X,\;y\in Y}$ and $x+y=x'+y'$ whenever $x-x'=y'-y \in Z$.
The formula
$$\norm {w} = \inf \setof{ \norm {x}_X + \norm {y}_Y }{w = x + y}$$
defines a norm on $W$ such that the canonical embeddings are isometric.
Finally, $W$ is a separable Banach space: It can be seen as a suitable quotient of $X\oplus Y$ endowed with the $\ell_1$-norm.
\end{pf}

Clearly, $\sbaniso$ has an initial object, the zero space. Thus, directedness follows from amalgamation.

\begin{lm}\label{ewtijwpeitj}
$\sbaniso$ is $\sig$-complete.
\end{lm}

\begin{pf}
The completion of a countable chain of separable Banach spaces is a separable Banach space.
\end{pf}

\begin{tw}\label{qijwqwrfp}
Assume $2^{\aleph_0} = \aleph_1$. There exists a Banach space $V$ of density $\aleph_1$ such that every Banach space of density $\loe\aleph_1$ is linearly isometric to a subspace of $V$ and every linear isometry $\map TXY$ between separable subspaces of $V$ can be extended to a linear isometry of $V$. Moreover, the space $V$ is unique, up to a linear isometry.
\end{tw}

\begin{pf}
Assuming the Continuum Hypothesis, there are only $\aleph_1$ many isometric types of separable Banach spaces and there are only $\aleph_1$ many types of linear isometries. Thus, by Lemmata~\ref{dsgsdgerr} and \ref{ewtijwpeitj} and by Theorem \ref{jofjaiopf}, $\sbaniso$ has a \fra\ sequence $\wek u$ of length $\omega_1$. We may further assume that this sequence is continuous. Let $V$ be the co-limit of $\wek u$ in the category of all Banach spaces. 

Fix a Banach space $X$ of density $\loe\aleph_1$. We can write $X=\bigcup_{\al<\omega_1}X_\al$, where $\sett{X_\al}{\al<\omega_1}$ is an increasing chain of closed separable subspaces of $X$ such that $X_\delta=\cl(\bigcup_{\xi<\delta}X_\xi)$ for every limit ordinal $\xi<\omega_1$. Translating it to the language of category theory, we obtain a continuous $\omega_1$-sequence in $\sbaniso$ whose co-limit, in the category of all Banach spaces, is $X$.
By Theorem \ref{continussdfsdfarfqarw}, there is an arrow of sequences $\map F{\wek x}{\wek u}$. This arrow has a co-limit in the category of all Banach spaces, which is just a linear isometric embedding of $X$ into $V$.

The second statement is obtained by the back-and-forth principle, using the continuity of $\wek u$.
\end{pf}

It has been shown (using non-trivial arguments) in~\cite{ACCGM} that the space $V$ from the theorem above is not isomorphic to any $\C(K)$ space.
The next simple statement provides a short and elementary argument that $V$ cannot be linearly isometric to any $\C(K)$ space.

\begin{prop}\label{progoelts}
Let $K$ be a compact space which contains at least two points. Then there exists a linear isometry $\map TXY$ between $1$-dimensional subspaces of $\C(K)$, which cannot be extended to a linear isometry of\/ $\C(K)$.
\end{prop}

\begin{pf}
Fix $a\ne b$ in $K$. Let $X$ consist of all constant functions on $K$. Let
$$\map R{\C(\dn ab)}{\C(K)}$$ be a regular extension operator for the inclusion $\dn ab\subs K$. That is, $R$ is a linear operator which assigns to each $f\in\C(\dn ab)$ its extension $Rf\in\C(K)$ so that $R1=1$ and $Rf\goe0$ whenever $f\goe0$. For example, let $(Rf)(t)=\phi(t)f(a) + (1-\phi(t))f(b)$, where $\map\phi K{[0,1]}$ is a continuous function such that $\phi(a)=1$ and $\phi(b)=0$ (which exists by Urysohn's Lemma).
Note that $R$ is an isometric embedding of $\C(\dn ab)$ into $\C(K)$.

Now define $\map TX{\C(K)}$ by $T1 = R1_{\sn a}$, where $1_{\sn a}$ is the function which takes value $1$ at $a$ and value $0$ at $b$.
Let $Y=\img TX=\setof{\lam R1_{\sn a}}{\lam \in \Err}$. Then $T$ is an isometry.

Suppose $\map{\ovr T}{\C(K)}{\C(K)}$ is a linear isometry extending $T$.
Let $$v=\inv{(\ovr T)}{R1_{\sn b}}.$$
Then $\norm v = 1$. By compactness, there exists $t\in K$ such that $\abs{v(t)}=1$. Let $\al=v(t)$ and consider $u=R(\al 1_{\sn a}+1_{\sn b})=\al R1_{\sn a}+R1_{\sn b}$. Notice that $\norm u=1$, while $\norm{(\ovr T)^{-1}(u)}=\norm{\al 1_K+v}\goe\abs{\al+v(t)}=2$. Hence $\norm{(\ovr T)^{-1}}\goe2$, a contradiction.
\end{pf}

As we have mentioned above, under the Continuum Hypothesis, the \v Cech-Stone remainder of the natural numbers $\Nat^*$ is the inverse limit of the $\omega_1$-\fra\ sequence in the category of nonempty compact metric spaces with quotient maps.
Thus it is natural to expect that our space $V$ is isometric (or at least isomorphic) to $\ell_\infty \by c_0 = \C(\Nat^*)$.
This is not the case, because of Proposition~\ref{progoelts} (and its isomorphic version in \cite{ACCGM}).

\separator

\index{\sban}\index{category!-- \sban}
We now turn to a more special class, namely Banach spaces with projectional resolutions.
From this point on, we consider the category $\sban$ whose objects are again all separable Banach spaces and arrows are linear operators of norm $\loe1$. We shall apply the results of Section \ref{reterpeteairs}.

\index{PRI}\index{projectional resolution of the identity}
A \emph{projectional resolution of the identity} (briefly: \emph{PRI}) on a Banach space $X$ is a sequence of norm-one projections $\sett{P_\al}{\al < \kappa}$ of $X$, satisfying the following conditions:
\begin{enumerate}
	\item[(P1)] $P_\al \cmp P_\beta = P_\al = P_\beta \cmp P_\al$ for every $\al, \beta < \kappa$.
	\item[(P2)] $\img {P_\delta}X = \cl{\bigcup_{\xi < \delta}\img{P_\xi}X}$ for every limit ordinal $\delta < \kappa$.
	\item[(P3)] $\dens(\img {P_\al} X) \loe \al + \aleph_0$.
	\item[(P4)] $X = \cl{\bigcup_{\al<\kappa} \img{P_\al}X}$.
\end{enumerate}
We are interested in the case of $\kappa = \aleph_1$, where the existence of a PRI is equivalent to the existence of the so-called countably 1-norming Markushevich basis, a natural generalization of a Schauder basis.
Banach spaces with norming Markushevich bases are often called \emph{Plichko spaces}.
We refer to \cite{Kal1} or \cite{KKLP} for details.
In view of Lemma~\ref{kwerjpit}, it is clear that a PRI in a space $X$ of density $\aleph_1$ is induced by a semi-continuous sequence in $\rp\sban$.

Our aim is to obtain a complementably universal Banach space with a projectional resolution of the identity, that is, we want that every other Banach space with a PRI is linearly isometric to its 1-complemented subspace.
\index{complemented space}
Recall that a Banach space $X$ is \emph{1-complemented} in $Y$ if there exists a projection $\map P Y Y$ such that $\norm P \loe 1$ and $\img P Y = X$.

\begin{lm}\label{ghwoieqtrzyr}
Let $\map fZX$, $\map gZY$ be left-invertible arrows in $\sban$. Then $f,g$ have a pushout in $\sban$.
\end{lm}

\begin{pf}
It is well-known and standard to check that the amalgamation described in the proof of Lemma~\ref{dsgsdgerr} is actually the pushout of $f, g$ in the category of Banach spaces with linear operators of norm $\loe1$.
\end{pf}

We now have all ingredients needed to show the existence of a universal Banach space with a PRI.

\begin{tw}\label{ThmPliCoUniv1}
Assume the Continuum Hypothesis. There exists a Banach space $U$ with a projectional resolution of the identity and of density $\aleph_1$, which has the following properties:
\begin{enumerate}
	\item[$(1)$] The family $\setof{X\subs U}{X\text{ is $1$-complemented in }U}$ is, modulo linear isometries, the class of all Banach spaces of density $\loe\aleph_1$ that admit a PRI.
	\item[$(2)$] Given separable 1-complemented subspaces $X,Y\subs U$, given a bijective linear isometry $\map T X Y$, there exists a bijective linear isometry $\map H U U$ such that $H \rest X = T$.
\end{enumerate}
Moreover, properties (1) and (2) together with the existence of a PRI describe the space $U$ uniquely, up to a linear isometry.
\end{tw}

\begin{pf}
The proof is completely analogous to that of Theorem~\ref{ThmUnVldxyone}, using the category $\sban$ instead of the category of compact metric spaces, and inductive sequences instead of inverse sequences.
Regarding part of property (1), we need to invoke a result from \cite{K_lin} saying that the existence of a PRI in Banach spaces of density $\aleph_1$ is preserved under 1-complemented subspaces.
\end{pf}

One has to admit that Theorem~\ref{ThmPliCoUniv1} (except for the uniqueness part) is already contained in the last chapter of the monograph~\cite{KKLP}, however the construction and arguments were entirely based on the early draft of this work, which existed in a preprint form before the aforementioned monograph appeared.

It is well-known that the Gurari\u\i\ space has a monotone Schauder basis.
On the contrary, the space $V$ from Theorem~\ref{qijwqwrfp} cannot have any kind of bases (of length $\omega_1$), because this would imply the existence of an isomorphic variant of a PRI which in turn would imply some isomorphic properties that are not shared by all Banach spaces of density $\loe 2^{\aleph_0} = \aleph_1$.
An example of such a property is the existence of a strictly convex renorming.
It is well-known that the space $\ell_\infty \by c_0$ fails this property, yet it is obviously isometric to a subspace of $V$.
In particular, the space $V$ is not isomorphic to the space $U$ from Theorem~\ref{ThmPliCoUniv1}, since the latter one admits a strictly convex renorming (implied by the PRI).
We refer to the book~\cite{DGZ} for details on renorming theory of Banach spaces.

\separator

One can argue that Markushevich bases in non-separable Banach spaces are not the proper generalization of Schauder bases.
In fact, especially when viewed from category-theoretic perspective, the notion of a PRI is the natural generalization of \emph{monote finite-dimensional Schauder decompositions}, not Schauder bases.
A monotone finite-dimensional Schauder decomposition in a Banach space $X$ is a sequence of norm one projections $\ciag P$ satisfying conditions (P1) and (P4) with $\kappa = \omega$ and such that $\img {P_n}X$ is finite-dimensional for each $\ntr$.
On the other hand, the existence of a monotone Schauder basis is equivalent to the existence of a monotone finite-dimensional decomposition $\ciag P$ such that $\dim \img {(P_{n+1}-P_n)}X = 1$ for every $\ntr$.
The last property obviously generalizes to higher cardinals.
In fact, a Banach space $X$ of density $\aleph_1$ has a \emph{monotone Schauder basis} (of length $\omega_1$) if and only if it has a PRI $\sett{P_\al}{\al < \omega_1}$ such that the image of $P_{\al+1} - P_\al$ is 1-dimensional for every $\al < \omega_1$.
\index{Schauder basis}
Formally, a \emph{Schauder basis} of type $\delta$ (where $\delta$ is an ordinal) is a sequence $\sett{e_\xi}{\xi < \delta}$ of vectors in a Banach space $X$ such that for every $x \in X$ there are uniquely determined scalars $\sett{t_\xi}{\xi < \delta}$ such that
$$x = \sum_{\xi < \delta}t_\xi e_\xi,$$
where the convergence of the series is taken with respect to the norm.
Once this happens, for each $\al < \delta$ there is a canonical projection $P_\al$ defined by
$$P_\al\left( \sum_{\xi < \delta} t_\xi e_\xi \right) = \sum_{\xi < \al} t_\xi e_\xi.$$
The basis is \emph{monotone} if $\norm{P_\al} \loe 1$ for every $\al < \delta$.
It is clear that $\sett{P_\al}{\al < \delta}$ satisfies conditions (P1), (P2), (P4) and
\begin{enumerate}
	\item[(P5)] For each $\al < \delta$ the image of $P_{\al+1} - P_\al$ is 1-dimensional.
\end{enumerate}
On the other hand, given a sequence $\sett{P_\al}{\al < \delta}$ satisfying (P1), (P2), (P4) and (P5), we can choose for each $\al$ a vector $e_\al$ from the image of $P_{\al+1} - P_\al$, obtaining a Schauder basis of type $\delta$.

\index{initial subspace}
In order to formulate the result on universal Schauder bases, we need the following notion.
Given Banach spaces $Y \subs X$, we say that $Y$ is an \emph{initial subspace} of $X$ if there is a sequence of norm one projections $\sett{P_\al}{\al < \delta}$, where $\delta \loe \omega_1$, satisfying conditions (P1), (P2), (P4), (P5) and such that $Y = \img{P_0}X$.
Typical examples of initial subspaces are linear spans of initial parts of a Schauder basis.
Note that an initial subspace is 1-complemented and the trivial space is initial in $X$ if and only if $X$ has a monotone Schauder basis of type $\loe \omega_1$.
Given a Schauder basis $\sett{e_\xi}{\xi < \delta}$ in $X$, given a subset $S \subs \delta$, we say that $\sett{e_\xi}{\xi \in S}$ is a \emph{canonically 1-complemented subbasis} if the linear operator $\map {P_S} X X$, defined by conditions $P_S e_\xi = e_\xi$ for $\xi \in S$ and $P e_\xi = 0$ for $\xi \notin S$, has norm $\loe1$.
Finally, we say that a basis $\sett{v_\xi}{\xi < \delta}$ is \emph{isometric to} a subbasis of $\sett{e_\al}{\al < \eta}$ if there is an increasing function $\map \phi \delta \eta$ such that the linear operator $f$ defined by equations $f(v_\xi) = e_{\phi(\xi)}$ ($\xi < \delta$) is a linear isometric embedding.

\begin{tw}\label{ThmtwSchauderbasics}
Assume the Continuum Hypothesis.
Then there exists a Banach space $\Pel$ with a monotone Schauder basis $\sett{e_\al}{\al < \omega_1}$ and with the following properties:
\begin{enumerate}
	\item[(i)] Every monotone Schauder basis of type $\loe \omega_1$ is isometric to a canonically 1-complemented subbasis of $\sett{e_\al}{\al < \omega_1}$.
	\item[(ii)] Given separable initial subspaces $X, Y \subs \Pel$, each with a monotone Schauder basis, given a bijective linear isometry $\map h X Y$, there exists a bijective linear isometry $\map H \Pel \Pel$ such that $H \rest X = h$.
\end{enumerate}
Furthermore, conditions (i) and (ii) determine the space $\Pel$ uniquely, up to a linear isometry.
\end{tw}

\begin{pf}
We shall first consider a pair of categories $\fB \subs \fB'$ which has the mixed amalgamation property.
Namely, the objects of $\fB'$ are pairs of the form $\pair X {\sett{x_\xi}{\xi<\delta}}$, where $\delta < \omega_1$ and $\sett{x_\xi}{\xi<\delta}$ is a monotone Schauder basis in $X$.
A $\fB'$-arrow from $\pair X {\sett{x_\xi}{\xi<\delta}}$ into $\pair Y {\sett{y_\xi}{\xi<\eta}}$ is an isometric embedding of the basis $\sett{x_\xi}{\xi < \delta}$ onto a canonically 1-complemented subbasis of $\sett{y_\xi}{\xi < \eta}$.
The category $\fB$ has the same objects as $\fB'$, the arrows are more restrictive.
Namely, a $\fB$-arrow from $\pair X {\sett{x_\xi}{\xi<\delta}}$ into $\pair Y {\sett{y_\xi}{\xi<\eta}}$ is an isometric embedding $\map f X Y$ such that $f(x_\xi) = y_\xi$ for $\xi < \delta$ (note that $\sett{y_\xi}{\xi < \delta}$ is canonically 1-complemented in $\sett{y_\xi}{\xi < \eta}$ because the basis is monotone).

\begin{claim}\label{Clkltggiek}
The inclusion $\fB \subs \fB'$ has the mixed amalgamation property.
\end{claim}

\begin{pf}
It suffices to prove the mixed amalgamation involving a $\fB$-arrow $\map f Z X$ such that $f$ is the inclusion $Z \subs X$ and $X \by \img f Z$ is 1-dimensional.
For general $\fB$-arrows, the mixed amalgamation can later be proved by an easy induction, using continuity at limit steps.

So fix Banach spaces $Z, X, Y$ such that $Z \subs X$, $\sett{x_\xi}{\xi\loe\delta}$ is a monotone Schauder basis in $X$ such that $\sett{x_\xi}{\xi < \delta}$ is a Schauder basis for $Z$.
Fix also a $\fB'$-arrow $\map i Z Y$, where $\sett{y_\xi}{\xi<\rho}$ is a monotone Schauder basis and $\sett{i(x_\xi)}{\xi < \delta}$ is a canonically 1-complemented subbasis of $\sett{y_\xi}{\xi<\rho}$.
Let $i(x_\xi) = y_{\phi(\xi)}$, where $\map \phi \delta \rho$ is strictly increasing.
Let $\map P Y Z$ and $\map T X Z$ be the canonical projections.

Forgetting the structure induced by Schauder bases, we move into the category of separable Banach spaces $\sban$, where by Lemmata~\ref{ghwoieqtrzyr} and \ref{aofihapsfj} we conclude that the pair $\pair i P$, $\pair \subs T$ of $\rp\sban$-arrows has a proper amalgamation.
In other words, there exist a Banach space $W \sups Y$ and an isometric embedding $\map j X W$, together with norm one projections $\map S W Y$ and $\map R W X$, satisfying
$$j \rest Z = i, \; T \cmp R = P \cmp S, \; S \cmp j = i \cmp T \text{ and } R\rest Z = P.$$
Replacing $W$ by $Y + \img j X$, we may assume that $W \by Y$ is 1-dimensional (in fact, the pushout of the embeddings already has this property).
Define $w_\xi = y_\xi$ for $\xi < \rho$ and define $w_\rho = j(x_\delta)$.
Observe that $S j(x_\delta) = i T(x_\delta) = 0$, that is, $x_\delta \in \ker S$.
It follows that $\sett{w_\xi}{\xi\loe \rho}$ is a monotone Schauder basis of $W$ and $S$ is the canonical projection, that is, the inclusion $Y \subs W$ becomes a $\fB'$-arrow.
Finally, $j$ is a $\fB$-arrow, which is witnessed by the projection $R$.
\end{pf}

Now, since $2^{\aleph_0} = \aleph_1$, by Theorem~\ref{ThmtwIstnienief} we conclude that there exists a continuous $\omega_1$-sequence $\wek U$ in $\fB$ that is \fra\ over the inclusion $\fB \subs \fB'$.
It is clear that every monotone Schauder basis of length $\omega_1$ is the result of a continuous chain in $\fB$.
Thus, (i) follows from Theorem~\ref{TmixedamsCof}.

In order to show (ii), we shall move to another category.
Namely, let $\fK$ be the subcategory of $\sban$ consisting of all isometric embeddings $\map f X Y$ such that $\img f X$ is an initial subspace of $Y$ and $\sn 0$ is an initial subspace of $X$ (i.e., $X$ has a monotone Schauder basis of a countable type).

Adapting the arguments from Claim~\ref{Clkltggiek}, it is clear that $\fK$ has the amalgamation property.
Furthermore, co-limits in the category of Banach spaces provide a quasi-limiting operator on $\fK$, therefore by Theorem~\ref{Tfjogosdg} there is at most one continuous $\omega_1$-\fra\ sequence $\wek u$ in $\fK$ and it satisfies the back-and-forth principle.
This translates to the fact that the Banach space obtained as the co-limit of $\wek u$ in the category of Banach spaces satisfies (ii) and is unique up to a linear isometry.

In order to complete the proof, it suffices to show that $\wek u$ is isomorphic to the image of $\wek U$ under the forgetful functor from $\fB'$ into $\fK$.
Represent $\wek U$ as a suitable continuous chain $\sett{U_\al}{\al<\omega_1}$ of separable spaces.
Fix $\al < \omega_1$ and fix an isometric embedding $\map f {U_\al} Y$ such that $\img f {U_\al}$ is an initial subspace of $Y$.
We need to find an isometric embedding $\map g Y{U_\beta}$ with $\beta > \al$, such that $\img g Y$ is an initial subspace of $U_\beta$ and $g \cmp f$ is the inclusion $U_\al \subs U_\beta$.
Notice that, just by the definition of an initial subspace, $Y$ has a monotone Schauder basis extending the image of the fixed basis of $U_\al$.
In other words, $f$ is actually a $\fB$-arrow.
Using the fact that $\wek U$ is \fra\ over the inclusion $\fB \subs \fB'$, we find $\beta > \al$ and a $\fB'$-arrow $\map g Y {U_\beta}$ such that $g\cmp f$ is the inclusion $U_\al \subs U_\beta$.
Finally, note that a Banach space $E$ is initial in $F$ whenever $E \subs F$ and a fixed monotone Schauder basis of $E$ is canonically 1-complemented in a fixed monotone Schauder basis of $F$.
In particular, $\img g Y$ is an initial subspace of $U_\beta$, which says that $g \in \fK$. This shows that $\wek u$ is a \fra\ sequence in $\fK$.
\end{pf}

Let us mention that in the separable case, Pe\l czy\'nski~\cite{Pelczynski} constructed a complementably universal (in the sense of isomorphic embeddings) Banach space for Schauder bases of type $\omega$.
In the same paper, he also constructed a complementably universal unconditional basis of type $\omega$.
Its uncountable counterpart can be obtained by a suitable modification of the proof above.
We skip the details.

\separator

We finally recall a rather standard argument, going back to Ostrovki\u\i~\cite{Ostrov}, showing that the Continuum Hypothesis is really needed for the existence of the Banach spaces from Theorems~\ref{qijwqwrfp}, \ref{ThmPliCoUniv1} and~\ref{ThmtwSchauderbasics}.
The fact that a space satisfying Theorem~\ref{qijwqwrfp} cannot have density less than the continuum has been already noted in~\cite[Prop. 4.1]{ACCGM}.

\begin{prop}
Let $V$ be a Banach space containing a fixed isometric copy of $c_0$ and satisfying the following condition:
\begin{enumerate}
	\item[$(\times)$] For every Banach space $E \sups c_0$ such that $\dim(E\by c_0) = 1$ and $c_0$ is 1-comple\-men\-ted in $E$, there exists an isometric embedding $\map f E V$ such that $f \rest c_0 = \id{c_0}$.
\end{enumerate}
Then the density of $V$ is at least the continuum.
\end{prop}

\begin{pf}
We assume that $c_0 \subs V$.
Given an infinite set $A \subs \nat$, denote by $\chi_A$ its characteristic function, treated as an element of $\ell_\infty$.
Let $E_A = c_0 + \Err \chi_A \subs \ell_\infty$ endowed with the norm $\norm{x + \lam \chi_A} = \max \{\norm{x + \lam \chi_A}_\infty, |\lam|\}$.
By this way, the canonical projection onto $c_0$ that maps $\chi_A$ to $0$ has norm one and we may use condition ($\times$).
Let $\map {f_A}{E_A}V$ be an embedding such that $f_A \rest c_0 = \id{c_0}$ and let $v_A = f_A(\chi_A)$.
Now fix another infinite set $B \subs \nat$ such that $A \not\subs B$ and fix $m \in A \setminus B$.
Then $\norm{v_A - v_B} \goe \norm{v_B - 2 \chi_{\sn m}} - \norm{2\chi_{\sn m} - v_A} = 3 - 2=1$.
Finally, taking an independent family of cardinality $2^{\aleph_0}$ we obtain a 1-discrete subset of $V$ of the same cardinality.
\end{pf}

\begin{wn}
A Banach space satisfying the assertion either of Theorem~\ref{qijwqwrfp} or \ref{ThmPliCoUniv1} or \ref{ThmtwSchauderbasics} must have density at least $2^{\aleph_0}$.
\end{wn}

\subsection{The Cantor set}\label{SubSectCantorStDns}

Let $\sets^+$ denote the category of nonempty finite sets.
It is clear that one-to-one mappings have pushouts in $\sets^+$, therefore it has proper amalgamations.
We claim that $\ciagi{(\rp\sets^+)}$ is isomorphic to the following category $\fL$.
The objects of $\fL$ are pairs of the form $\pair K D$, where $K$ is a totally disconnected compact metric space and $D\subs K$ is a countable dense set.
An arrow from $\pair K D$ to $\pair LE$ is a pair of functions $\pair j f$, where $\map j D E$, $\map f L K$, $f\cmp j = \id D$ and $f$ is continuous.

Given a sequence $\wek x$ in $\rp \sets^+$, we take $D$ to be the co-limit of $\img e{\wek x}$ in the category of sets. Clearly, $D$ is a countable set.
Now take $K$ to be the inverse limit of $\img r{\wek x}$ in the category of topological spaces.
Why do we mention topology here? As we shall see in a moment, arrows correspond to continuous quotients.
Clearly, $K$ is compact metrizable and totally disconnected. In view of Lemma~\ref{oiutiq}, there is a canonical embedding of $D$ into $K$.
Thus, we may think of $D$ as a subset of $K$.
Notice that every projection from $K$ to an element of the sequence $\img r{\wek x}$ is a continuous quotient; in other words, it corresponds to a partition into clopen sets.

Fix two sequences $\wek x, \wek y$ in $\rp \sets^+$ and fix an arrow $\map {\wek f}{\wek x}{\wek y}$.
Let $\pair K D$ and $\pair L E$ be the $\fL$-objects corresponding to $\wek x$ and $\wek y$, respectively.
Notice that the sequence $e(f_n)$ ``converges" to a one-to-one map $\map j D E$, and the sequence $r(f_n)$ ``converges" to a continuous quotient $\map f L K$. Clearly, $f \cmp j = \id D$.
By this way we have described a functor from $\ciagi{(\rp\sets^+)}$ into $\fL$.

Now fix an $\fL$-object $\pair K D$ and write $D = \sett{d_n}{\ntr}$. Assume $D$ is infinite and its enumeration is one-to-one.
Construct inductively partitions $\Yu_n$ of $K$ into clopen sets, in such a way that for each $U\in \Yu_n$ there is a unique index $i \loe n$ such that $d_i \in U$.
Let $D_n = \{d_0, \dots, d_n\}$ and define $\map {r_n}{K}{D_n}$ so that $r_n^{-1}(d_i) \in \Yu_n$ and $r_n(d_i)=d_i$ for $i \loe n$.
Let $\map {e_n}{D_n}D$ be the inclusion.
Set $x_n^\infty = \pair {e_n}{r_n}$. Then $x_n^\infty$ is a $\rp\sets^+$-arrow and for each $n<m$ there is a unique $\rp\sets^+$-arrow $x_n^m$ satisfying $x_n^\infty = x_m^\infty \cmp x_n^m$.
It is clear that this defines a sequence $\wek x$ in $\rp\sets^+$ which induces the $\fL$-object $\pair K D$.
Call an $\fL$-arrow \emph{proper} if it corresponds to a proper arrow of sequences in $\rp\sets^+$.

\begin{lm}\label{Lmerigrrito}
An $\fL$-arrow $\map f {\pair K D}{\pair L E}$ is proper if and only if $f = \pair j r$, where $\map j K L$, $\map r L K$ are continuous maps such that $\img j D \subs E$, $\img r E \subs D$, and $r \cmp j = \id K$.
\end{lm}

The last three conditions actually imply that $\img r E = D$.

\begin{pf}
It is clear that a proper arrow of sequences induces a pair $\pair j r$ satisfying the conditions above.
Fix continuous maps $\map j K L$, $\map r L K$ satisfying $r \cmp j = \id K$ and $\img j D \subs E$, $\img r E = D$.
Without loss of generality, we may assume that $L \subs K$, $j$ is the inclusion, and $D = E \cap L$.
Choose a chain $\sett{D_n}{\ntr}$ of finite subsets of $D$ such that $\bigcup_{\ntr}D_n = D$.
Choose inductively finite sets $E_n\subs E$ so that $D_n \subs E_n$, $\img r{E_n} = D_n$, and $E_n \subs E_{n+1}$ for $\ntr$.
We can do it in such a way that $\bigcup_{\ntr}E_n = E$, because $\inv r D \sups E$.
Now observe that $f_n := \pair {e_n}{r_n}$ is a $\rp\sets^+$-arrow, where $e_n$ is the inclusion $D_n \subs E_n$ and $r_n = r\rest E_n$.
Finally, $\ciag f$ is a proper arrow from the sequence $\ciag D$ to the sequence $\ciag E$.
\end{pf}

The category $\rp\sets^+$ is certainly countable, therefore it has a \fra\ sequence.
It is easy to check that a sequence $\wek u$ is \fra\ in $\rp\sets^+$ if and only if the inverse limit of $\img r{\wek u}$ is the Cantor set $\Cantor$, together with the canonical projections.
Moreover, for every countable dense set $D$ of $\Cantor$ there exists a sequence $\wek u$ in $\rp\sets^+$ such that $\Cantor$ is the limit of $\img r{\wek u}$ and $D$ is the co-limit of $\img e{\wek u}$ in the category of sets.
In particular, we obtain the following properties of the Cantor set which belong to the folklore:

\begin{tw}
Let $Q$ denote the set of rational numbers in the standard representation of the Cantor set.
Then
\begin{enumerate}
	\item[(1)] Given a countable dense set $D \subs \Cantor$, there exists a homeomorphism $\map h \Cantor \Cantor$ such that $\img h D = Q$.
	\item[(2)] Given a compact totally disconnected metric space $K$ together with a countable dense set $D \subs K$, there exist a continuous surjection $\map f \Cantor K$ and a topological embedding $\map i K \Cantor$ such that $f \cmp i = \id K$, $\img i D \subs Q$ and $\img f Q = D$.
\end{enumerate}
\end{tw}

Homogeneity translates to the following statement, again belonging to the folklore:

\begin{tw}
Let $Q$ be a countable dense subset of $\Cantor$.
Given a nonempty finite set $s$, given continuous mappings $\map {p_i}{\Cantor}{s}$ and embeddings $\map {e_i}{s}{Q}$ for $i=0,1$, such that $p_0 \cmp e_0 = \id s = p_1 \cmp e_1$, there exists a homeomorphism $\map h \Cantor \Cantor$ such that $\img h Q = Q$, $p_1 \cmp h = p_0$ and $h \cmp e_0 = e_1$.
\end{tw}

Below we give two examples of improper arrows in $\ciagi{(\rp\sets^+)}$.

\begin{ex}\label{Exeimprprs}
(a) Let $x_n = \{0, \dots, n \}$, $y_n = \{0, \dots, n, \infty \}$.
Then $x_0\subs x_1 \subs \dots$ and $y_0 \subs y_1 \subs \dots$.
Let $K = L = \omega \cup \sn \infty$, where $\omega$ is discrete and $\infty=\lim_{n\to\infty}n$.
Given $n<m$, let $x_n^m\in {(\rp\sets^+)}$ be such that $e(x^n_m)$ is inclusion and $r(x_n^m)$ is the retraction mapping all elements of the set $\{n+1,\dots,m\}$ to $n$.
Let $y_n^m\in {(\rp\sets^+)}$ be such that $e(y_n^m)$ is again inclusion, while $r(y_n^m)$ maps the set $\{n+1,\dots,m,\infty\}$ to $\infty$.
This defines $\omega$-sequences $\wek x$ and $\wek y$ in ${\rp\sets^+}$.
Now let $\map {f_n}{x_n}{y_n}$ be defined by $f_n = \pair {e_n}{r_n}$, where $e_n$ is inclusion and $r_n$ maps $\infty$ to $n$.
Notice that, given $n<m$, the equation $y_n^m \cmp f_n = f_m \cmp x_n^m$ holds in $\rp\sets^+$, although this amalgamation is not proper.

The sequences $\wek x$ and $\wek y$ correspond to pairs $\pair K\omega$ and $\pair LL$ in $\fL$, respectively.
The sequence $\wek f = \sett{f_n}{\ntr}$ induces maps $\map e\omega L$ and $\map rLK$, where $e$ is the inclusion of $\omega$ into $L$ and $r$ is the identity. Note that $\img r L \not \subs \omega$, therefore $\pair e r$ is not proper.

Another example of an improper arrow can be obtained as follows. Consider the compact space $L = \Z \cup \dn{-\infty}{\infty}$, where $\Z$ is the (discrete!) set of the integers and $\lim_{n\to\infty}\infty$ and $\lim_{n\to \infty}-n = -\infty$. Let $K = \Z \cup \sn \infty$ be a quotient of $L$ obtained by identifying $-\infty$ with $\infty$. Let $\map r L K$ be the quotient map.
Furthermore, let $D = K$, $E = L$, and let $\map e D E$ be the inclusion.
We claim that $\pair e r$ is an arrow from $\pair K D$ into $\pair L E$.
Of course, it cannot be proper, because $e$ is even not continuous.

Let $x_n = \setof{ i \in \Z}{|i|\loe n} \cup \sn \infty$, $y_n = \setof{ i \in \Z}{|i|\loe n} \cup \dn {-\infty}\infty$.
Let $x_n^m$ be such that $e(x_n^m)$ is the inclusion $x_n \subs x_m$ and $r(x_n^m)$ is identity on $x_n$ and maps all $j\in x_m \setminus x_n$ to $\infty$.
Define $y_n^m$ in a similar manner, with the difference that $r(y_n^m)(j) = -\infty$ whenever $j < -n$ and $r(y_n^m)(j) = \infty$ whenever $j > n$.
We have defined sequences $\wek x$ and $\wek y$ corresponding to $\pair K K$ and $\pair L L$, respectively.
Finally, let $f_n = \pair {e_n}{r_n}$, where $e_n$ is the inclusion $x_n \subs y_n$ and $r_n$ is a quotient map of $y_n$ onto $x_n$ that maps $-\infty$ onto $\infty$ and satisfies $r_n \cmp e_n = \id {x_n}$.
Then $\wek f = \sett{f_n}{\ntr}$ is an improper arrow of sequences inducing $\pair e r$.
\end{ex}

\separator

We have just seen above that the Cantor set $\Cantor$ together with its fixed countable dense subset $Q$ can be seen as the co-limit of a \fra\ sequence in $\rp\sets^+$.
Using the concepts from Section~\ref{reterpeteairs},
we shall now find a universal family of continuous self-functions on $\Cantor$ that preserve $Q$.

Namely, we define the following category $\fK$.
The objects of $\fK$ are triples of the form $\triple a b {\sett{f_i}{i \in s}}$, where $a,b, s$ are nonempty finite sets, and $\map {f_i} a b$ is a mapping for each $i \in s$.
Given $\fK$-objects $p = \triple a b {\sett{f_i}{i \in s}}$, $q = \triple c d {\sett{g_i}{i \in t}}$,
a $\fK$-arrow from $p$ to $q$ is a triple $\triple k \ell \phi$ such that $\map k a c$, $\map \ell b d$ are $\rp \sets^+$-arrows and $\map \phi s t$ is a one-to-one function such that for each $i\in s$ the following diagrams commute:
$$
\xymatrix{
c \ar[rr]^{g_{\phi(i)}} & & d \\
a \ar@{ >->}[u]^{e(k)} \ar[rr]_{f_i} & & b \ar@{ >->}[u]_{e(\ell)}
}
\qquad \qquad
\xymatrix{
c \ar@{->>}[d]_{r(k)} \ar[rr]^{g_{\phi(i)}} & & d \ar@{->>}[d]^{r(\ell)} \\
a \ar[rr]_{f_i} & & b
}
$$
Recall that also $r(k) \cmp e(k) = \id a$ and $r(\ell) \cmp e(\ell) = \id b$.
It is obvious that $\fK$ is countable and has a weakly initial object (a single function between one-element sets).
The key point is to show that $\fK$ has the amalgamation property.
For this aim, we need the following extension property, whose proof is rather trivial.

\begin{lm}\label{Lmoouwghow}
Let $\map f a b$ be a $\sets^+$-arrow.
Let $\map i a {a'}$, $\map j b {b'}$ be $\rp\sets^+$-arrows.
Then there exists a $\sets^+$-arrow $\map {f'} {a'} {b'}$ for which the diagrams
$$
\xymatrix{
a' \ar[r]^{f'} & b' \\
a \ar@{ >->}[u]^{e(i)} \ar[r]_f & b \ar@{ >->}[u]_{e(j)}
}
\qquad \qquad
\xymatrix{
a' \ar@{->>}[d]_{r(i)} \ar[r]^{f'} & b' \ar@{->>}[d]^{r(j)} \\
a \ar[r]_f & b
}
$$
commute.
\end{lm}

An amalgamation of two $\fK$-arrows will be called \emph{proper} if the two squares involving $\rp\sets^+$ are proper amalgamations.

\begin{lm}
$\fK$ has the proper amalgamation property.
\end{lm}

\begin{pf}
Fix $\fK$-objects $p = \triple a b {\sett{f_i}{i \in s}}$, $q^j = \triple {c^j}{d^j}{\sett{g^j_i}{i \in t^j}}$, where $j = 0,1$, and let $\triple {k_j}{\ell_j}{\phi_j}$ be $\fK$-arrows from $p$ to $q^j$.
By (the proof of) Lemma~\ref{aofihapsfj}, there exist $\rp \sets^+$-arrows $\map {k'_j}{c^j}{\ovr c}$ and $\map {\ell'_j}{d^j}{\ovr d}$,
where $j=0,1$, such that
$$
\xymatrix{
c^0 \ar@{ >->}[rr]^{e(k_0')} & & \ovr c \\
a \ar@{ >->}[u]^{e(k_0)} \ar@{ >->}[rr]_{e(k_1)} & & c^1 \ar@{ >->}[u]_{e(k_1')}
}
\qquad \qquad
\xymatrix{
d^0 \ar@{ >->}[rr]^{e(\ell_0')} & & \ovr d \\
b \ar@{ >->}[u]^{e(\ell_0)} \ar@{ >->}[rr]_{e(\ell_1)} & & d^1 \ar@{ >->}[u]_{e(\ell_1')}
}
$$
are pushout squares in $\sets$.
In particular, $\pair {k_0'}{k_1'}$ is a proper amalgamation of $\pair {k_0}{k_1}$ and $\pair {\ell_0'}{\ell_1'}$ is a proper amalgamation of $\pair {\ell_0}{\ell_1}$.
Without loss of generality, we may assume that $s = t^0 \cap t^1$ and $\phi_0$, $\phi_1$ are inclusions.
Let $\ovr t = t^0 \cup t^1$ and let $\map {\phi_j'}{t^j}{\ovr t}$ be the inclusion for $j=0,1$.

Given $i \in t^0 \setminus s$, using Lemma~\ref{Lmoouwghow}, choose $\map {\ovr g_i}{\ovr c}{\ovr d}$ satisfying 
$$\ovr g_i \cmp e(k_0') = e(\ell_0') \cmp g_i^0 \oraz r(\ell_0') \cmp \ovr g_i = g_i^0 \cmp r(k_0').$$
Similarly, given $i \in t^1 \setminus s$, find $\ovr g_i$ satisfying analogous conditions.
In order to complete the proof, we still need to define $\ovr g_i$ for $i \in s$, so that $w = \triple {\ovr c}{\ovr d}{ \sett{\ovr g_i}{i \in \ovr t} }$ will become a $\fK$-object which together with $k_0', k_1',\ell_0'\ell_1'$ witnesses the proper amalgamation property.

Fix $i \in s$ and define
$$\al = g_i^0 \cmp r(k_0'), \qquad \beta = g_i^1 \cmp r(k_1').$$
Notice that
\begin{align*}
r(\ell_0) \cmp \al
&= r(\ell_0) \cmp g_i^0 \cmp r(k_0') = f_i \cmp r(k_0) \cmp r(k_0') \\
&= f_i \cmp r(k_1) \cmp r(k_1') = r(\ell_1) \cmp g_i^1 \cmp r(k_1')
= r(\ell_1) \cmp \beta.
\end{align*}
By the property of a pushout, there exists a unique map $\map {\ovr g_i}{\ovr c}{\ovr d}$ satisfying
$$r(\ell_0') \cmp \ovr g_i = \al \oraz r(\ell_1') \cmp \ovr g_i = \beta.$$
We still need to check that $g_i^0 \cmp r(k_0') = r(\ell_0') \cmp \ovr g_i$ and $g_i^1 \cmp r(k_1') = r(\ell_1') \cmp \ovr g_i$.
By symmetry, we shall prove the first equality only.
Define
$$\delta = r(\ell_0') \cmp \ovr g_i \cmp e(k_1').$$
We have
\begin{align*}
\delta \cmp e(k_1)
&= r(\ell_0') \cmp \ovr g_i \cmp e(k_1' \cmp k_1) = r(\ell_0') \cmp e(\ell_1') \cmp g_i^1 \cmp e(k_1) \\
&= r(\ell_0') \cmp e(\ell_1') \cmp e(\ell_1) \cmp f_i = r(\ell_0') \cmp e(\ell_0') \cmp e(\ell_0) \cmp f_i \\
&= e(\ell_0) \cmp f_i = g_i^0 \cmp e(k_0).
\end{align*}
This shows that $\pair {g_i^0}{\delta}$ is an amalgamation of $\pair {e(k_0)}{e(k_1)}$.
By the property of a pushout, there is a unique mapping $\map \gamma {\ovr c}{d^0}$ satisfying
\begin{equation}
\gamma \cmp e(k_0') = g_i^0 \oraz \gamma \cmp e(k_1') = \delta.
\tag{$*$}\label{Equjborgt}
\end{equation}
Finally, it remains to check that both $\gamma_1 = g_i^0 \cmp r(k_0')$ and $\gamma_2 = r(\ell_0') \cmp \ovr g_i$ satisfy (\ref{Equjborgt}).
That $\gamma_2$ satisfies (\ref{Equjborgt}) is clear.
It is also obvious that $\gamma_1 \cmp e(k_0') = g_i^0$.
Finally,
\begin{align*}
\gamma_1 \cmp e(k_1')
&= g_i^0 \cmp r(k_0') \cmp e(k_1') = g_i^0 \cmp e(k_0) \cmp r(k_1) \\
&= e(\ell_0) \cmp f_i \cmp r(k_1) = e(\ell_0) \cmp r(\ell_1) \cmp g_i^1 \\
&= r(\ell_0') \cmp e(\ell_1') \cmp g_i^1 = r(\ell_0') \cmp \ovr g_i \cmp e(k_1') = \delta,
\end{align*}
which completes the proof.
\end{pf}

The proof above could be made simpler and shorter, using the point-structure of sets and the explicit form of pushouts of sets.
On the other hand, purely category theoretic arguments have the advantage of possible applicability to other categories.

Once we have proper amalgamations, it is natural to consider the category $\ciagi \fK$ with \emph{proper arrows}, where the meaning of being proper is rather obvious.

Let $\wek u$ be a \fra\ sequence in $\fK$, where $u_n = \triple {v_n}{w_n}{\sett{z_i^n}{i \in s_n}}$.
We may assume that $s_0 \subs s_1 \subs \dots$ and that the bonding maps $u_m^n$ are of the form $\triple {v_m^n}{w_m^n}{\psi_m^n}$, where $\psi_m^n$ is the inclusion of $s_m$ into $s_n$.

\begin{lm}
$\wek v$ and $\wek w$ are \fra\ sequences in $\rp \sets^+$.
\end{lm}

\begin{pf}
This is an immediate consequence of Lemma~\ref{Lmoouwghow}.
\end{pf}

It follows that both $\wek v$ and $\wek w$ determine the Cantor set together with a fixed countable dense (as we know, it can be chosen arbitrarily).
Finally, in order to show that $\wek u$ provides a sequence of continuous self-mappings of the Cantor set that is universal both in the sense of quotients and in the sense of embeddings, we need to decode the information contained in the category $\ciagi \fK$.

Namely, let $\fL_p$ be the subcategory of $\fL$ consisting of all proper arrows (see Lemma~\ref{Lmerigrrito}).
Recall that the objects of $\fL_p$ (as well as of $\fL$) are pairs $\pair K D$, where $K$ is a 0-dimensional compact metrizable space and $D$ is a countable dense subset of $K$.
Let $\fB$ be the category whose objects are the same as the objects of $\fL$, while a $\fB$-arrow from $\pair K D$ to $\pair L E$ is a continuous map $\map f K L$ satisfying $\img f D \subs E$.

We now define the following category $\fF$, in the same spirit as the category $\fK$ defined above.
The objects of $\fF$ are triples of the form $\triple A B { \sett{f_n}{n \in S} }$, where $A$, $B$ are $\fL$-objects, $S$ is a nonempty countable set and $\map {f_n} A B$ is a $\fB$-arrow.

An $\fF$-arrow from $X = \triple A B { \sett{f_n}{n \in S} }$ to $X' = \triple {A'} {B'} { \sett{f'_n}{n \in S'} }$ is a triple $f = \triple k \ell \phi$, where $\map k A {A'}$, $\map \ell B {B'}$ are $\fL_p$-arrows, $\map \phi S {S'}$ is a one-to-one mapping such that
$$f'_{\phi(n)} \cmp e(k) = e(\ell) \cmp f_n \oraz r(\ell) \cmp f'_{\phi(n)} = f_n \cmp r(k)$$
holds for every $n \in S$.

\begin{lm}\label{Lmueghr}
The canonical functor from $\ciagi \fK$ to $\fF$, identifying the objects of $\fK$ with constant sequences, is surjective on the class of objects.
\end{lm}

\begin{pf}
As usual, we identify $\fK$ with constant sequences.
It is obvious how to define a suitable functor $\map F {\ciagi \fK} \fF$ which is the identity on $\fK$.
It remains to check that every $\fF$-arrow is induced by a proper arrow of the corresponding sequences in $\fK$.

Fix an $\fF$-object $X = \triple A B { \sett{f_n}{n \in S} }$, where $A = \pair K D$, $B = \pair L E$.
For convenience, we shall assume that $S = \nat$ (if $S$ is finite, the argument becomes simpler).
Let $\wek a$, $\wek b$ be two sequences in $\rp\sets^+$ such that $K = \liminv \img r {\wek a}$, $L = \liminv \img r {\wek b}$, $D = \bigcup_{\ntr} a_n$, $E = \bigcup_{\ntr} b_n$, assuming $\ciag a$ and $\ciag b$ are chains of finite sets such that $e(a_m^n)$ and $e(b_m^n)$ are inclusions.
For each $\ntr$ we can find a strictly increasing function $\map {\phi_n} \nat \nat$ such that $f_n$ is the limit of $\sett{f_n^k}{k\in\nat}$, or more precisely, for each $k, n \in \nat$ the squares
$$
\xymatrix{
D \ar[rr]^{f_n} & & E \\
a_k \ar@{ >->}[u] \ar[rr]_{f_n^k} & & b_{\phi_n(k)} \ar@{ >->}[u]
}
\qquad \qquad
\xymatrix{
K \ar@{->>}[d] \ar[rr]^{f_n} & & L \ar@{->>}[d] \\
a_k \ar[rr]_{f_n^k} & & b_{\phi_n(k)}
}
$$
are commutative, where the vertical arrows are induced by the sequences $\wek a$, $\wek b$.
Define
$$\psi(n) = \max \{ \phi_0(n), \phi_1(n), \dots, \phi_n(n) \}$$
and $s_n = \{ 0, 1, \dots, n \}$.
Given $n \loe k$, replace $f_n^k$ by $b_{\phi(n)}^{\psi(n)} \cmp f_n^k$.
Now
$$x_k = \triple {a_k} {b_k} { \sett{f_n^k}{k \in s_n} }$$
is a $\fK$-object and the bonding maps of $\wek a$ and $\wek b$ can be used to compose a sequence $\wek x$ converging to $X$ in the sense that each $f_n$ is the co-limit of $\sett{f_n^k}{k\goe n}$.
\end{pf}

It is not hard to check that the functor mentioned in Lemma~\ref{Lmueghr} witnesses the equivalence of $\fF$ to the category of sequences in $\fK$ with proper arrows.

Summarizing, we obtain the following result concerning continuous functions on the Cantor set.

\begin{tw}
Let $Q$ be a fixed countable dense set in the Cantor set $\Cantor$.
Then there exists a sequence $\sett{u_n}{\ntr}$ of continuous self-maps of the Cantor set such that $\img {u_n}Q \subs Q$ for every $\ntr$ and the following condition is satisfied.
\begin{enumerate}
	\item[($\ddagger$)] Given a sequence $\sett{\map {f_n}K L}{\ntr}$ of continuous functions between 0-dim\-en\-sio\-nal compact metric spaces, given countable dense sets $D \subs K$, $E\subs L$ such that $\img {f_n}D \subs E$ for every $\ntr$, there exist topological embeddings $\map i K \Cantor$, $\map j L \Cantor$ together with continuous surjections $\map p \Cantor K$, $\map q \Cantor L$, satisfying $p \cmp i = \id K$, $q \cmp j = \id L$, $\img i D \subs Q$, $\img p Q = D$, $\img j E \subs Q$, $\img p Q = E$, and there exists a one-to-one function $\map \theta \nat \nat$ such that for every $\ntr$ the following diagrams commute.
$$
\xymatrix{
\Cantor \ar[rr]^{u_{\theta(n)}} & & \Cantor \\
K \ar@{ >->}[u]^i \ar[rr]_{f_n} & & L \ar@{ >->}[u]_j
}
\qquad \qquad
\xymatrix{
\Cantor \ar@{->>}[d]_p \ar[rr]^{u_{\theta(n)}} & & \Cantor \ar@{->>}[d]^q \\
K \ar[rr]_{f_n} & & L
}
$$
\end{enumerate}
\end{tw}

On could also formulate the homogeneity part, describing the sequence $\ciag u$ uniquely up to isomorphism.
The precise statement seems to be too technical, therefore we have decided to omit it.

\subsection{Linear orders}

Historically, the set $\Qyu$ of rational numbers was the first universal homogeneous object, discovered in this context by Cantor at the end of 19th century.
In the realm of \fra\ theory, its properties are rather obvious and easily explained.
Since the category of linear orders with embeddings has the amalgamation property, uncountable versions of the set of rational numbers may exist, as usual, subject to a ``proper" cardinal arithmetic.
For example, the Continuum Hypothesis is equivalent to the statement ``there exists a linearly ordered set $\Qyu_{\omega_1}$ of cardinality $\aleph_1$ that contains isomorphic copies of all linearly ordered sets of cardinality $\loe \aleph_1$ and is homogeneous with respect to its countable subsets".
Indeed, assuming CH, this is just the consequence of the \frajon\ theory. On the other hand, if such a set $\Qyu_{\omega_1}$ exists then, assuming $\Qyu \subs \Qyu_{\omega_1}$, every extension $\Qyu \subs \Qyu \cup \sn t$ with $t \in \Err$ must be realized in $\Qyu_{\omega_1}$, showing that $|\Qyu_{\omega_1}| \goe 2^{\aleph_0}$.

Notice that every embedding between finite linearly ordered sets is left-invertible.
This may serve as a motivation for studying uncountable versions of the set of rational numbers, where homogeneity is considered with respect to left-invertible embeddings.
It turns out that no cardinal arithmetic assumptions are needed for the existence of such a set of size $\aleph_1$.
In particular, we shall see a natural example of a directed $\sig$-complete category with amalgamations that is dominated by a single arrow.

Namely, consider the category $\LO$ of nonempty countable linearly ordered sets with increasing maps.
By ``increasing" we mean ``preserving the non-strict ordering".
In particular, constant maps are increasing.

Let $\map f12$ be the inclusion map, where, as usual, $1=\sn\emptyset$ and $2=\dn{\emptyset}{1}$. It is easy to see that the pair $\pair ff$ has no pushout in $\LO$.
However, we shall see that the category $\rp\LO$ has proper amalgamations.

Given linearly ordered sets $A,B$ we denote by $A \cdot B$ the set $A \times B$ endowed with the lexicographic ordering.
Note that, fixing $p \in B$, the pair $\pair e r$ is a $\rp\LO$-arrow from $A$ to $A \cdot B$, where $e(a) = \pair a p$ and $r(a,b) = a$ for $a \in A$, $b \in B$.

In particular, when $B = \Qyu$ and $p=0$, the $\rp\LO$-arrow $\map {\ell_A} A {A \cdot \Qyu}$ just described above will be called \emph{canonical}.

\begin{lm}
The category $\rp\LO$ has proper amalgamations.
\end{lm}

\begin{pf}
Let $\map f Z X$, $\map g Z Y$ be $\rp \LO$-arrows.
Without loss of generality, we may assume that $X \cup Y \subs \Qyu$ and $X \cap Y = Z$.
Let us identify $X$ and $Y$ with suitable subsets of $Z \cdot \Qyu$.
Since $\Qyu$ is homogeneous, we may actually assume that the embedding of $Z$ into $Z \cdot \Qyu$ is of given by the formula $z \mapsto \pair z 0$.
Now, consider the $\rp \LO$-arrows $\map p {Z \cdot \Qyu}{Z \cdot \Qyu \cdot \Qyu}$ and $\map q {Z \cdot \Qyu}{Z \cdot \Qyu \cdot \Qyu}$ defined by $e(p)(z,s) = \triple z s 0$, $r(p)(z,s,t) = \pair z s$ and $e(q)(z,s) = \triple z 0 s$, $r(q)(z,t,u) = \pair z u$.
We shall actually consider $p$ restricted to $X$ and $q$ restricted to $Y$.
We claim that
$$\xymatrix{
Y \ar[rr]^q  & & Z \cdot \Qyu \cdot \Qyu \\
Z \ar[rr]_f \ar[u]^g & & X \ar[u]_p
}$$
is a proper amalgamation.
Fix $x = \pair t s \in X$.
Then $r(f)(x) = t$ and $e(g)(r(f)(x)) = \pair t 0$.
On the other hand, $e(p)(x) = \triple t s 0$ and $r(q)(e(p)(x)) = \pair t 0$.

Now fix $y = \pair t s \in Y$.
Then $r(g)(y) = t$ and $e(f)(r(g)(y)) = \pair t 0$.
On the other hand, $e(q)(y) = \triple t 0 s$ and $r(p)(e(q)(y)) = \pair t 0$.
This shows that the amalgamation provided by $p$ and $q$ is proper.
\end{pf}

Recall that $\ell_\Qyu$ is the canonical $\rp\LO$-arrow from $\Qyu$ to $\Qyu \cdot \Qyu$.
Note that, up to isomorphism, $\ell_\Qyu$ is a unique $\rp \LO$-arrow $h$ with domain $\Qyu$ and such that all fibers of $r(h)$ are isomorphic to $\Qyu$.

\begin{lm}
The arrow $\ell_\Qyu$ is dominating in $\rp \LO$.
\end{lm}

\begin{pf}
Condition (D1) in the definition of a dominating family of arrows is obvious, because given a countable linear order $X$, the lexicographic product $X \cdot \Qyu$ is isomorphic to $\Qyu$.
Fix a $\rp \LO$-arrow $\map f \Qyu Y$.
Let $\map {\ell_Y} Y {Y\cdot \Qyu}$ be the canonical $\rp\LO$-arrow and let $h = \ell_Y \cmp f$.
Notice that all fibers of $r(h)$ are isomorphic to $\Qyu$, therefore $h$ is isomorphic to $\ell_\Qyu$.
\end{pf}

Now we know that the category $\rp \LO$ has an $\omega$-\fra\ sequence, which we do not find of any particular interest.
Since $\LO$ is $\sig$-complete, we can also talk about $\omega_1$-\fra\ sequences.
By the results of Section~\ref{reterpeteairs}, there exists a unique semi-continuous $\omega_1$-sequence in $\rp \LO$.
Its natural co-limit $\Quna$ can be regarded as a natural uncountable variant of the set of rational numbers, existing without any extra set-theoretic assumptions.
Note that the $\omega_1$-\fra\ sequence has the property that each bonding arrow is isomorphic to $\ell_\Qyu$.
With this information at hand, it is straightforward to see that
$$\Quna = \setof{x \in \Qyu^{\omega_1}}{ |\suppt(x)| < \aleph_0},$$
endowed with the lexicographic ordering (recall that $\suppt(x) = \setof{\al}{x(\al) \ne0}$).

We now come back to the theory of Valdivia compact spaces, in the context of linearly ordered spaces.
Namely, every linearly ordered set induces a natural interval topology.
This topology is compact if and only if every subset has the least upper bound (the supremum of the empty set is the minimal element).
There is a natural duality (formally: dual equivalence) between the category of linearly ordered sets with increasing maps and the category of compact 0-dimensional linearly ordered spaces with continuous increasing maps preserving the minimum and the maximum.
To be more precise, by a \emph{compact line} we mean a compact linearly ordered space $K$ with distinguished elements $0_K = \min K$ and $1_K = \max K$.
The arrows of our category are continuous increasing maps preserving the distinguished elements.
In particular, compact lines are supposed to be nonempty.
The duality can be described briefly as follows.
Given a linearly ordered set $X$, we let $K(X)$ to be the set of all increasing maps from $X$ into $2 = \dn 01$.
This is clearly a compact 0-dimensional line.
Conversely, given a compact 0-dimensional line $K$, we define $L(K)$ to be the set of all arrows into $2 =\dn 01$ treated as a compact line.
Specifically, $L(K)$ consists of all continuous increasing maps $\map f K 2$ satisfying $f(0) = 0$ and $f(1) = 1$.
Again, this is a linearly ordered set and we consider it without topology.
Both operations canonically extend to contravariant functors, providing the dual equivalence of the two categories.
The dual equivalence described above is actually a special case of Priestley duality~\cite{DaveyPriestley} between partially ordered sets and distributive lattices.

It turns out that Valdivia compact lines cannot be too big.
For this aim, we quote the following result from~\cite{K_classR}:

\begin{prop}[{\cite[Props. 5.5, 5.7]{K_classR}}]\label{ThmKlinesRetsV}
Valdivia compact lines are precisely those compact spaces that can be represented as limits of continuous inverse sequences of compact metrizable lines whose bonding maps are increasing and right-invertible (in the category of continuous maps).
In particular, Valdivia compact lines have weight $\loe \aleph_1$.
\end{prop}

It has been furthermore proved in~\cite{K_classR} that the class of (nonempty!) Valdivia compact lines has a universal pre-image $V_{\omega_1}$ which is 0-dimensional.
More precisely, every Valdivia compact is a continuous increasing quotient of $V_{\omega_1}$.
As a byproduct, we get the fact that every Valdivia compact line is a continuous increasing quotient of a 0-dimensional one (this can be also proved more directly).

For a moment, let us consider the class of 0-dimensional Valdivia compact lines.
In view of Proposition~\ref{ThmKlinesRetsV}, this class can be analyzed by looking at continuous inverse sequences in the category $\fL$ of topologically right-invertible increasing quotient maps between 0-dimensional compact metric lines.
It has been proved in~\cite{K_classR}, without referring to \fra\ sequences, that there is a single increasing quotient that dominates this category.
Namely, consider the Cantor set $\Cantor$ endowed with the usual linear order.
Actually, all linear orderings of $\Cantor$ compatible with its topology are isomorphic, so by the ``usual one" we may consider the lexicographic ordering.
Call a point $t \in \Cantor$ \emph{rational} if it is eventually constant when viewed as a function from $\omega$ to $2$.
In the language of orderings, rational points are precisely those that are isolated from one side (including 0 and 1).
Let the map $\map \pi {\Cantor}{\Cantor}$ be defined by the following conditions: $\pi^{-1}(t)$ is the Cantor set, if $t$ is rational and $|\pi^{-1}(t)| = 1$, otherwise.
One can easily ``extract" from the proof of~\cite[Thm. 5.8]{K_classR} that this map dominates the category $\fL$.
The Valdivia compact line $V_{\omega_1}$ is the inverse limit of an $\omega_1$-sequence $\wek v$ consisting of such maps.
Thus, $\wek v$ is a \fra\ sequence in $\fL$ and by this way we obtain the following fact, not discovered in~\cite{K_classR}:

\begin{tw}
Given an increasing isomorphism $\map h A B$, where $A, B$ are closed metrizable subsets of\/ $V_{\omega_1}$, there exists an increasing isomorphism $\map H {V_{\omega_1}}{V_{\omega_1}}$ such that $H \rest A = h$.
\end{tw}

\begin{pf}
Note that every increasing isomorphism is a homeomorphism and every closed subset of a 0-dimensional compact line is its increasing retract.
Thus, after moving to a suitable category of embedding-projection pairs, the statement follows from Theorem~\ref{Thmhot55}.
\end{pf}

One has to admit that in the proof of~\cite[Thm. 5.8]{K_classR} it had been shown that the canonical functor from $\fL$ into the category of all compact metric lines has the mixed amalgamation property and the sequence $\wek v$ mentioned above is \fra\ over this functor.
Thus, the argument leading to the fact that every Valdivia compact line is an increasing continuous image of $V_{\omega_1}$ can be explained by Theorem~\ref{TmixedamsCof}.

Coming back to linearly ordered sets and to the \fra\ limit $\Quna$ of the category $\LO$, let us mention a result from~\cite{KalKub}, characterizing linearly ordered sets $X$ for which $K(X)$ is Valdivia compact:

\begin{prop}[{\cite[Thm. 3.1]{KalKub}}]
Given a linearly ordered set $X$, the compact line $K(X)$ is Valdivia if and only if $X$ has the following properties:
\begin{enumerate}
	\item[$(1)$] $|X| \loe \aleph_1$.
	\item[$(2)$] Every bounded monotone $\omega_1$-sequence in $X$ is convergent.
	\item[$(3)$] Given a stationary set $S \subs \omega_1$, given a function $\map f S X$, there exists a stationary set $T \subs S$ such that $f \rest T$ is monotone.
\end{enumerate}
\end{prop}

Condition (3) perhaps looks somewhat artificial and not natural, however it appears to be one of the proper uncountable versions of the Bolzano-Weierstrass principle, saying that every sequence in a linearly ordered set contains a monotone subsequence.

It is important to note that for a continuous increasing map the property of being right-invertible in the category of topological spaces does not imply being right-invertible in the category of compact lines, because a right inverse may not preserve 0 and 1, although it is always increasing.
In fact, it is well-known that every closed set (not necessarily containing 0, 1) of a 0-dimensional compact line is its topological (even increasing) retract.
It is easy to check that if $\map f {K(Y)}{K(X)}$ is a topologically right-invertible increasing quotient map, then $L(f)$ corresponds to an increasing embedding $\map e X Y$ such that $\img e X$ is an increasing retract of its convex hull
$$\conv(\img e X) = \setof{y \in Y}{(\exists\;x_0, x_1 \in X) \;\; e(x_0) \loe y \loe e(x_1)}.$$
On the other hand, if  $\map f X Y$ is left-invertible in the category of linearly ordered sets then $K(f)$ is right-invertible in the category of compact lines with continuous increasing maps.
One can easily adapt the arguments of~\cite{KalKub} for obtaining the following characterization:

\begin{tw}
Let $V$ be a linearly ordered set of cardinality $\aleph_1$.
Then $V$ is order isomorphic to $\Quna$ if and only if it satisfies the following conditions:
\begin{enumerate}
	\item[$(1)$] Every monotone $\omega_1$-sequence in $V$ is convergent.
	\item[$(2)$] Given countable sets $A, B \subs V$ such that $a < b$ for every $a \in A$, $b \in B$, there exists $p \in V$ such that $a < p$ for $a \in A$ and $p < b$ for $b \in B$ (we allow the possibility that one of the sets $A,B$ is empty).
	\item[$(3)$] Given a stationary set $S \subs \omega_1$, given a function $\map f S X$, there exists a stationary set $T \subs S$ such that $f \rest T$ is monotone.
\end{enumerate}
\end{tw}

We refer to \cite{KalKub} for details, where one can find all arguments showing that (1) and (3) are necessary and sufficient for the fact that $V$ is the co-limit of a sequence in $\LO$.
Finally, condition (2) says that any $\omega_1$-sequence with co-limit $V$ is \fra\ in $\LO$.

\subsection*{Acknowledgments}

A significant part of this research was done during the author's longer research visits at the following institutions: The Fields Institute for Research in Mathematical Sciences (January--March 2007 and October 2012), Department of Mathematics, University of Murcia (November 2011) and the Department of Mathematical Analysis, University of Valencia (April--October 2010).
The author wishes to thank all these institutions for their warm hospitality.

\printindex
\end{document}